\documentclass{elsart}

\usepackage{amssymb}
\usepackage{amsfonts}
\usepackage{amsmath}
\usepackage{rotating}
\usepackage{graphicx,color}
\usepackage{subfigure}
\usepackage{bm}
\usepackage{color}
\usepackage{multirow}
\usepackage{enumitem}
\usepackage{tikz}
\usepackage{breakcites}
\usepackage{epstopdf}
\usepackage{morefloats}
\usepackage[pdfpagelabels, hypertexnames=true, plainpages=false, naturalnames=false, breaklinks=true]{hyperref}

\usetikzlibrary{decorations.markings}
\usetikzlibrary{shapes,arrows}
\usetikzlibrary{shapes.geometric}
\usetikzlibrary{fit}	
\renewcommand*{\today}{20 July 2015}

\newcommand\change[1]{#1}

\begin{document}

\setcounter{page}{1}
\allowdisplaybreaks[4]

\begin{frontmatter}

\title{Adaptive Algebraic Multiscale Solver for Compressible Flow in Heterogeneous Porous Media}

\author[TuDelft]{Matei \c{T}ene\thanksref{label1}\corauthref{cor1}}
\author[Stanford]{Yixuan Wang\thanksref{label2}}
\author[TuDelft]{Hadi Hajibeygi\thanksref{label3}}

\address[TuDelft]{Department of Geoscience and Engineering, Faculty of Civil Engineering and Geosciences, Delft University of Technology, P.O. Box 5048, 2600 GA Delft, The Netherlands.}
\address[Stanford]{Department of Energy Resources Engineering, Stanford University, 367 Panama St., Rm. 065, Stanford, CA 94305-2220, USA.}

\thanks[label1]{m.tene@tudelft.nl}
\thanks[label2]{yixuanw@stanford.edu}
\thanks[label3]{h.hajibeygi@tudelft.nl}

\corauth[cor1]{Corresponding author.}

\begin{abstract}
This paper presents the development of an Adaptive Algebraic Multiscale Solver for Compressible flow (C-AMS) in heterogeneous porous media. Similar to the recently developed AMS for incompressible (linear) flows [Wang et al., JCP, 2014], C-AMS operates by defining primal and dual-coarse blocks on top of the fine-scale grid. These coarse grids facilitate the construction of a conservative (finite volume) coarse-scale system and the computation of local basis functions, respectively. However, unlike the incompressible (elliptic) case, the choice of equations to solve for basis functions in compressible problems is not trivial. Therefore, several basis function formulations (incompressible and compressible, with and without accumulation) are considered in order to construct an efficient multiscale prolongation operator. As for the restriction operator, C-AMS allows for both multiscale finite volume (MSFV) and finite element (MSFE) methods. Finally, in order to resolve high-frequency errors, fine-scale (pre- and post-) smoother stages are employed. In order to reduce computational expense, the C-AMS operators (prolongation, restriction, and smoothers) are updated adaptively. In addition to this, the linear system in the Newton-Raphson loop is infrequently updated. Systematic numerical experiments are performed to determine the effect of the various options, outlined above, on the C-AMS convergence behaviour. An efficient C-AMS strategy for heterogeneous 3D compressible problems is developed based on overall CPU times. Finally, C-AMS is compared against an industrial-grade Algebraic MultiGrid (AMG) solver. Results of this comparison illustrate that the C-AMS is quite efficient as a nonlinear solver, even when iterated to machine accuracy.

\end{abstract}


\end{frontmatter}
\endNoHyper

\section{Introduction}

Accurate and efficient simulation of multiphase flow in large-scale heterogeneous natural formations is crucial for a wide range of applications, including hydrocarbon production optimization, risk management of Carbon Capture and Storage, water resource utilizations and geothermal power extractions. Unfortunately, considering the size of the domain along with the high resolution heterogeneity of the geological properties, such numerical simulation is often beyond the computational capacity of traditional reservoir simulators. Therefore, Multiscale Finite Element (MSFE) \cite{Hou97,Hou99,EffendievBookSiam,Efendiev00,AarnesHou02} and Finite Volume (MSFV) \cite{Jenny03,Jenny06} methods and their extensions have been developed to resolve this challenge.  \change{A comparison of different multiscale methods, based on their original descriptions, has been studied in the literature \cite{Kippe06}.}
MSFV and MSFE methods map a discrete fine-scale system to a much coarser space. In MultiGrid (MG) terminology \cite{MG}, this map is considered as a special prolongation operator, represented by locally-supported (and adaptively updated) basis functions \cite{Zhou}. The restriction operator is then defined based on either a Finite Element (MSFE), Finite Volume (MSFV), or a combination of both.

MSFV has been applied to a wide range of applications (see, e.g., \cite{Zhou,Zhou-tams,Ivan-gravity,SeongTrJcp,imsfv-jcp,hadi-frac-jcp,hadi-erest-spej,Wolfsteiner-well,hadi-compositional-spej,yixuan-ams,Davide14,Olav_unstruct}), thus recommending multiscale as a very promising framework for the next-generation reservoir simulators. However, most of these developments, including the state-of-the-art algebraic multiscale formulation (AMS) \cite{yixuan-ams}, have focused on the incompressible (linear) flow equations. 

When compressibility effects are considered, the pressure equation becomes nonlinear, and its solution requires an iterative procedure involving a parabolic-type linear system of equations \cite{Aziz79}. Therefore, the development of an efficient and general algebraic formulation for compressible nonlinear flows is crucial in order to advance the applicability of multiscale methods towards more realistic problems. 

The present study introduces the first algebraic multiscale iterative solver for compressible flows in heterogeneous porous media (C-AMS), along with a thorough study of its computational efficiency (CPU time) and convergence behaviour (number of iterations). 

In contrast to cases with incompressible flows, the construction of basis functions for compressible flow problems is not straightforward. In the past, incompressible elliptic \cite{Ivan-compressible}, compressible elliptic \cite{Zhou,SEON-BO}, and pressure-independent parabolic \cite{hadi-comp-jcp} basis functions have been considered. However, the literature lacks a systematic study to reveal the benefit of using one option over the other, especially when combined with a fine-scale smoother stage. Moreover, no study of the overall efficiency of the multiscale methods (based on the CPU time measurements) has been done so far for compressible three-dimensional  problems.

In order to develop an efficient prolongation operator, in this work, several formulations for basis functions are considered. These basis functions differ from each other in the amount of compressibility involved in their formulation, ranging from incompressible elliptic to compressible parabolic types. In terms of the restriction operator, both MSFE and MSFV are considered, along with the possibility of mixing iterations of the former with those of the latter, allowing C-AMS to benefit from the Symmetric Positive Definite (SPD) property of MSFE and the conservative physically correct solutions of MSFV. 

The low-frequency errors are resolved in the global (multiscale) stage of C-AMS, while high-frequency errors are tackled using a second-stage smoother at fine-scale. In this paper, we consider two options for the smoothing stage: the widely used local correction functions with different types of compressibility involved (i.e., more general than the specific pressure-independent operator \cite{hadi-comp-jcp}), as well as ILU(0) \cite{Saad}. The best C-AMS procedure is determined among these various strategies, on the basis of the CPU time for 3D heterogeneous problems. It is important to note that the setup and linear system population are measured alongside the solve time - a study which has so far not appeared in the previously published compressible multiscale works.  

Though C-AMS is a conservative method (i.e., only a few iterations are enough in order to obtain a high-quality approximation of the fine-scale solution), in the benchmark studies of this work, it is iterated until machine accuracy is reached. And, thus, its performance as an exact solver is compared against an industrial-grade Algebraic MultiGrid (AMG) method, SAMG \cite{SAMG}.
This comparative study for compressible problems is the first of its kind, and is made possible through the presented algebraic formulation, which allows for easy integration of C-AMS in existing advanced simulation platforms. Numerical results, presented for a wide range of heterogeneous 3D domains, illustrate that the C-AMS is quite efficient for  simulation of nonlinear compressible flow problems. 

The paper is structured as follows. First, the Compressible Algebraic Multiscale Solver (C-AMS) is presented, where several options for the prolongation, restriction operators as well as the second-stage solver are considered. Then the adaptive updating of the C-AMS operators are studied, along with the possibility of infrequent linear system updates in the Newton-Raphson loop. Numerical results are subsequently presented for a wide range of 3D heterogeneous test cases, aimed at determining the optimum strategy. Finally, the C-AMS is compared with an Algebraic MultiGrid Solver (i.e., SAMG) both in terms of the number of iterations and overall CPU time. 


\section{Compressible Flow in Heterogeneous Porous Media}

Single phase compressible flow in porous media, using Darcy's law (without gravity and capillary effects), can be stated as:
\begin{align}\label{mass}
 \frac{\partial}{\partial t} (\phi \rho) - \nabla~ .~ \big(\rho \ \bm \lambda \cdot \nabla p\big)=\rho {q},
\end{align}
where $\phi$, $\rho$, and $q$ are the porosity, density, and the source terms, respectively. Moreover, $\bm \lambda=\bm{K}/\mu$ is the fluid mobility with positive-definite permeability tensor, $\bm{K}$, while $\mu$ is the fluid viscosity.

The semi-discrete form of this nonlinear flow equation using implicit (Euler-backward) time integration reads
\begin{equation}\label{p-nonlinear}
\frac{\phi^{n+1}}{\Delta t}-\frac{\phi^{n} {\rho^n}}{\Delta t \rho^{n+1}} - \frac {1}{\rho^{n+1}}{\nabla~ .~ \big(\rho^{n+1} \bm \lambda \cdot \nabla p^{n+1} \big)} = q,
\end{equation}
which is linearized as 
\begin{align}\label{lin_p}
{c^{\nu}}(p^{\nu+1}-p^{\nu}) - \frac{1}{\rho^{\nu}}&\nabla~ .~ (\rho^{\nu}\bm \lambda \cdot \nabla p^{\nu+1})= b^\nu,
\end{align}
where
\begin{align}\label{p_lhs}
c^{\nu}=\frac{1}{\Delta t} \bigg[
\frac{\partial \phi}{\partial p} \bigg\arrowvert ^{\nu}-\phi^n \frac{\partial }{\partial p} \bigg(\frac{1}{\rho}\bigg){{\bigg\arrowvert}}^{\nu}\rho^n \bigg]
\end{align}
and
\begin{align}\label{p_rhs}
b^{\nu} = -\frac{\phi^{\nu}}{\Delta t}+\frac{\phi^n \ \rho^{n}}{\Delta t \ \rho^{\nu}} + q.
\end{align}
The superscripts $(\nu)$ and $(\nu+1)$ denote the old and new Newton-Raphson iteration levels, respectively. Note that, as $(\nu\rightarrow\infty)$, Eq.~\eqref{lin_p} converges to the nonlinear Eq.~\eqref{p-nonlinear}, and $\big( p^{\nu+1} - p^{\nu} \big) \rightarrow 0$. Therefore, the coefficient $c$, which is a by-product of the linearization lemma, plays a role only during iterations. This fact opens up the possibility to alter $c$ by computing it based on either $p^{\nu}$ (resulting in $c^{\nu}$) or $p^n$ (corresponding to $c^n$) -the pressure at the previous time-step. Each choice can potentially lead to a different convergence behaviour and, thus, computational efficiency. 

Algebraically, Eq.~\eqref{lin_p} can be written for the unknown pressure vector $\bm{p}$ as
\begin{align}\label{ap=r}
\big(\bm{C}^{\nu} + \bm{\tilde{A}}^{\nu}\big){\bm{p}}^{\nu+1} \equiv {\bm{A}}^{\nu} {\bm{p}}^{\nu+1} =  {\bm{f}}^{\nu}  \equiv {\bm{b}}^{\nu} + \bm{C}^{\nu} \bm{p}^{\nu},
\end{align}
where $\bm{C}^{\nu}$ is a diagonal matrix having $(c^{\nu}_i \ dV_i)$ at cell $i$ in its $i$-th diagonal entry, where $dV_i$ is the volume of cell $i$. Also, $\bm{\tilde{A}}$ is the convective compressible flow matrix, having fine-scale transmissibilities computed on the basis of a finite-volume scheme as entries. Moreover, the vector $\bm{b}^{\nu}$ contains the integrated source terms in the fine-scale volumes, i.e., $(q_i \ dV_i)$. The total Right-Hand-Side (RHS) terms are denoted by the vector $\bm{f}^{\nu}$.


\section{Compressible Algebraic Multiscale Solver (C-AMS)}

\begin{figure} [htb!]
\centering
\includegraphics [width=0.9\linewidth]{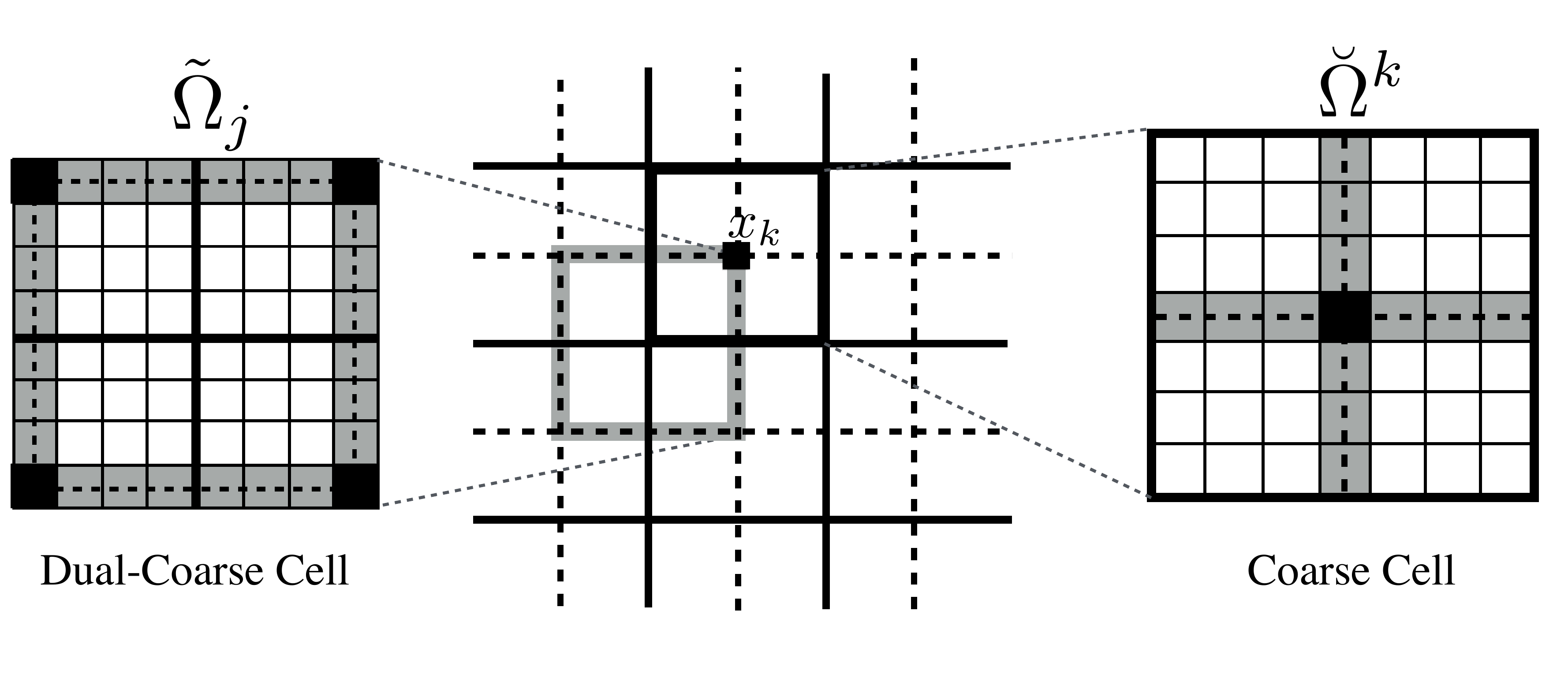}
\caption{Multiscale grids imposed on the given fine grid. A primal- and a dual-coarse block are highlighted on the right and left sides, respectively.}
\label{msfv-grid}
\end{figure}

The C-AMS relies on the primal- and dual-coarse grids, which are superimposed on the fine-scale grid (See Fig.~\ref{msfv-grid}). There are $N_p$ and $N_d$ coarse and dual-coarse grid cells in a domain with $N_f$ fine-grid cells.

The transfer operators between fine-scale and coarse-scale are defined as the multiscale Restriction ($\bm{R}$) and Prolongation ($\bm{P}$). The former is defined based on either Finite Element (MSFE), i.e., $\bm{R}^{FE} = \bm{P}^T$, or Finite Volume (MSFV), for which $\bm{R}^{FV}$ corresponds to the integral over primal-coarse blocks, i.e.,
\begin{equation}
\bm{R}^{FV}(i,j) = \begin{cases} 
	1 & \text{, if fine-cell } j \text{ is contained in primal-coarse block } i \\
	0 & \text{, otherwise.}
\end{cases}
\end{equation}

The columns of $\bm{P}$ are the basis functions, which are computed on dual-coarse cells (see Fig.~\ref{msfv-grid}), subject to simplified boundary conditions (the localization assumption).

In contrast to the (incompressible) AMS, C-AMS can be formulated based on different choices of the basis functions, depending on the level of compressibility involved. The first two types read
\begin{align}
{c^\nu} \Phi^{\nu+1}_{k,h} - \frac{1}{\rho^{\nu}} \nabla(\rho^{\nu}\bm \lambda \cdot \nabla \Phi^{\nu+1}_{k,h})= 0 \label{B1}
\end{align}
and
\begin{align}
 - \frac{1}{\rho^{\nu}} \nabla\cdot(\rho^{\nu}\bm \lambda\cdot \nabla \Phi^{\nu+1}_{k,h})= 0, \label{B2}
\end{align}
both being pressure dependent (through $c$ and $\rho$), but different in the sense of the consideration of the accumulation term, $c$. Alternatively, one can also formulate basis functions using 
\begin{align}
{c^{n}} \Phi_{k,h} - \nabla\cdot(\bm \lambda \cdot \nabla \Phi_{k,h})&=0 \label{B3},
\end{align}
or
\begin{align}
- \nabla\cdot(\bm \lambda \cdot \nabla \Phi_{k,h})&=0 \label{B4}
\end{align}
which are both pressure independent (since $c$ is now based on the pressure from the previous time step).
All of these equations are subject to reduced-problem boundary conditions along dual-coarse cell boundaries $\partial \tilde \Omega^h$ \cite{yixuan-ams}. One can also obtain the equations for the corresponding four types of local correction functions, $\Psi_h$, by substituting the corresponding RHS term in Eqs.~\eqref{B1}-\eqref{B4}.
As mentioned before, in this work, systematic studies on the basis of the CPU time as well as the number of iterations are performed in order to find the optimum formulation for basis function (i.e., prolongation operator). 

\begin{figure} [htb!]
\centering
\includegraphics[width=0.35\textwidth]{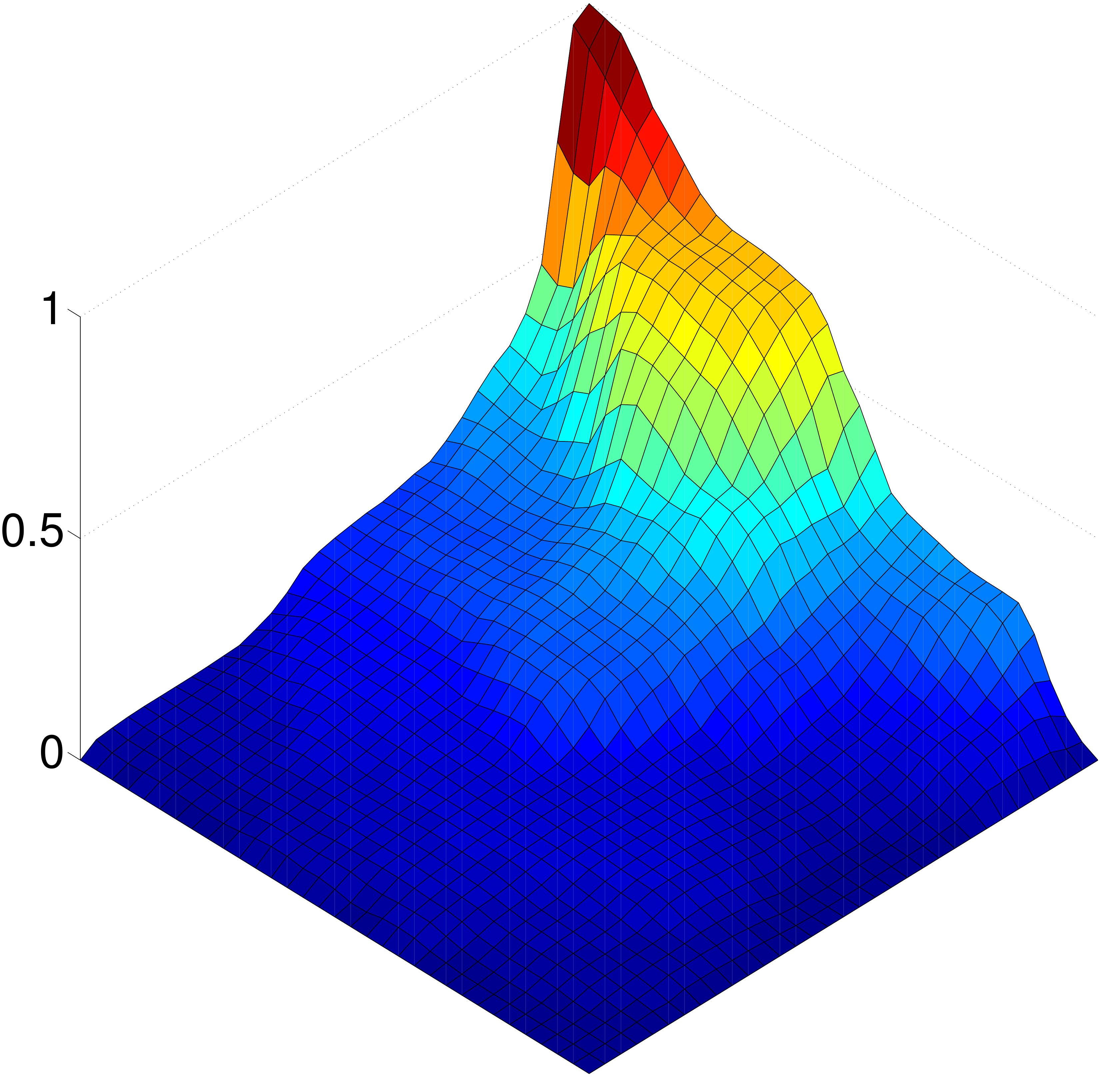}\hspace{1cm}
\includegraphics[width=0.35\textwidth]{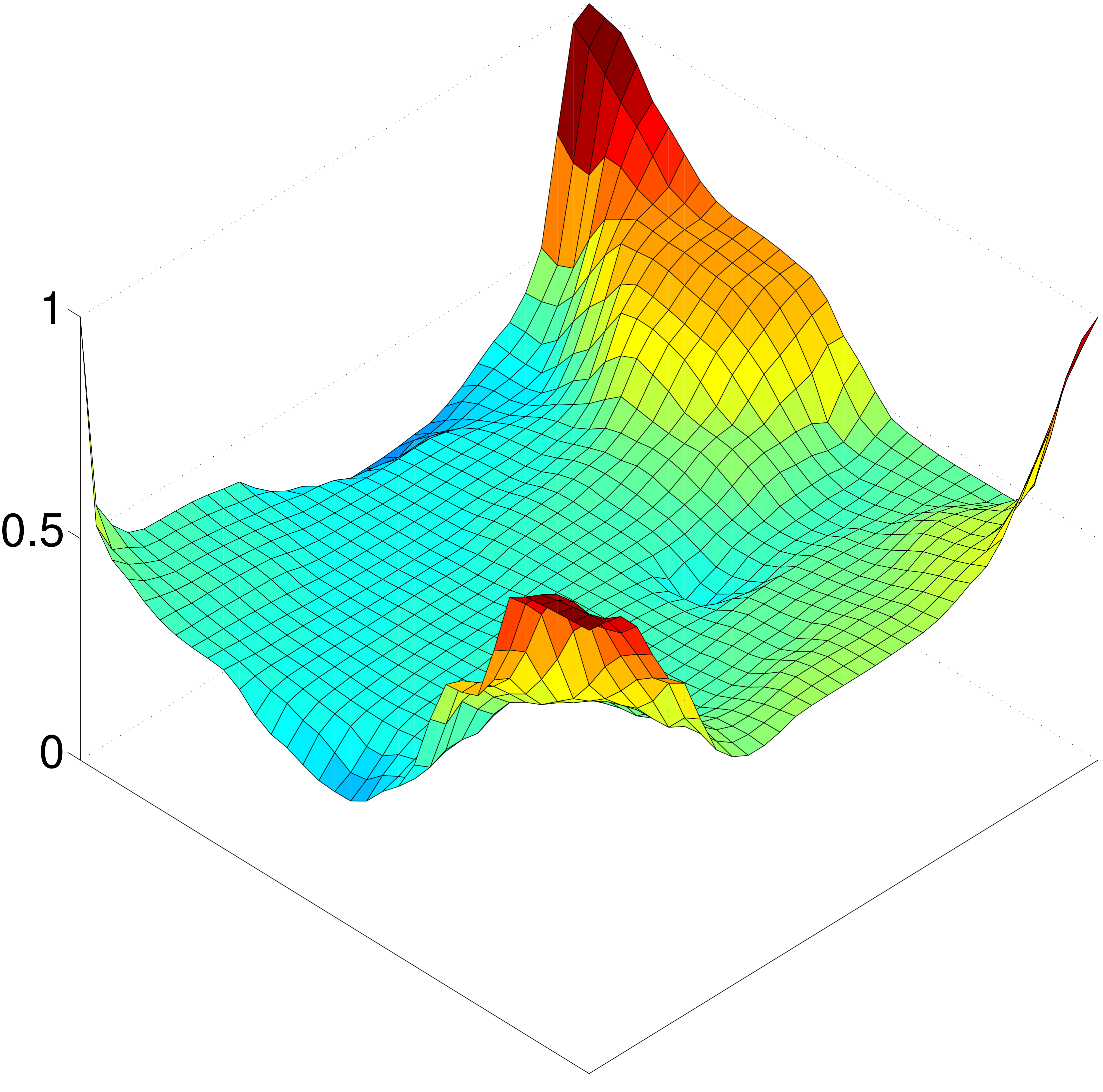}\\
(a): $c^\nu  \Phi^{\nu+1} -  \frac{1}{\rho^{\nu}} \nabla~ .~ (\rho^{\nu} {\bm{\lambda}}  \cdot \nabla \Phi^{\nu+1}) = 0$ \\
\vspace{3mm}
\includegraphics[width=0.35\textwidth]{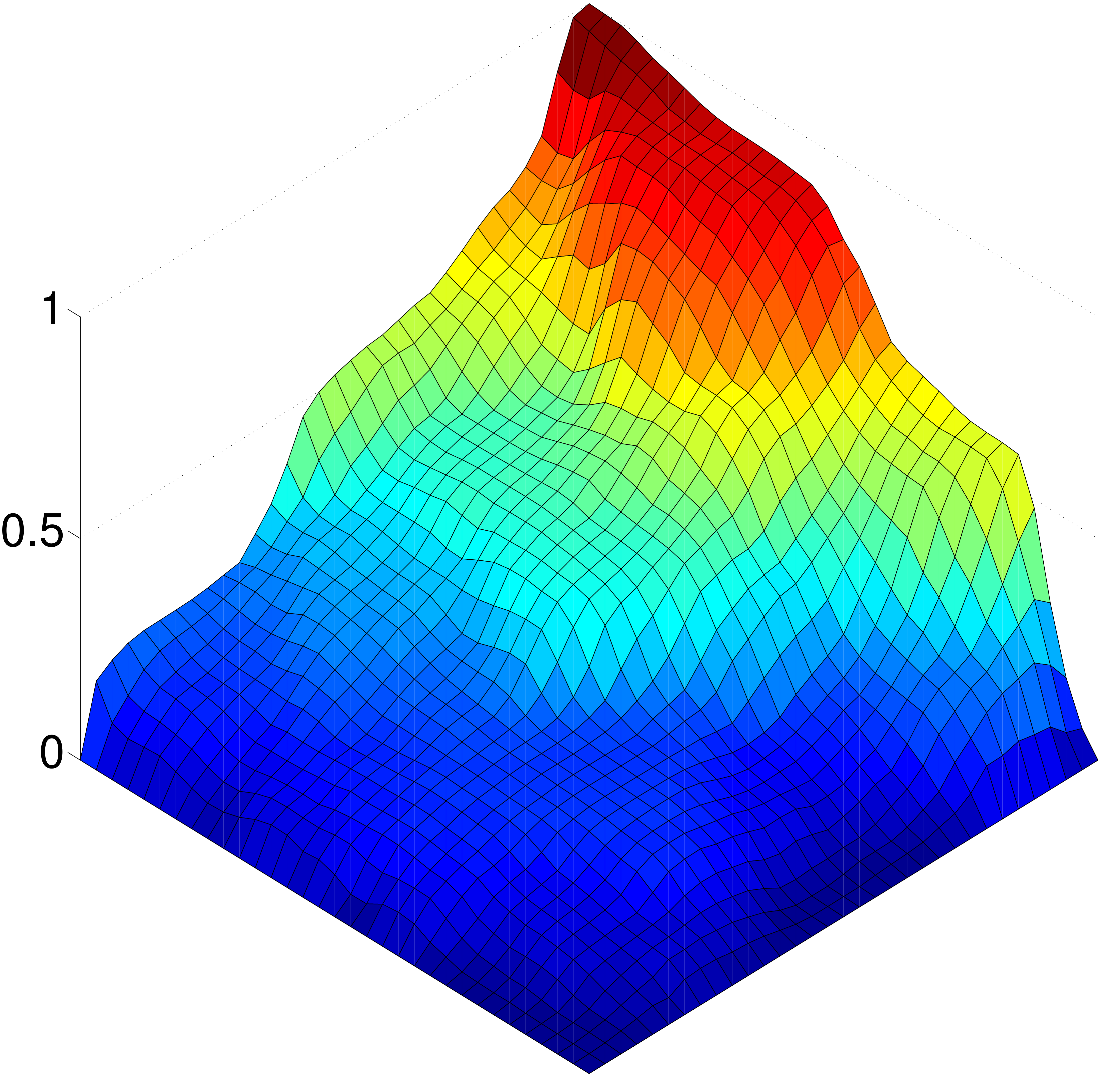}\hspace{1cm}
\includegraphics[width=0.35\textwidth]{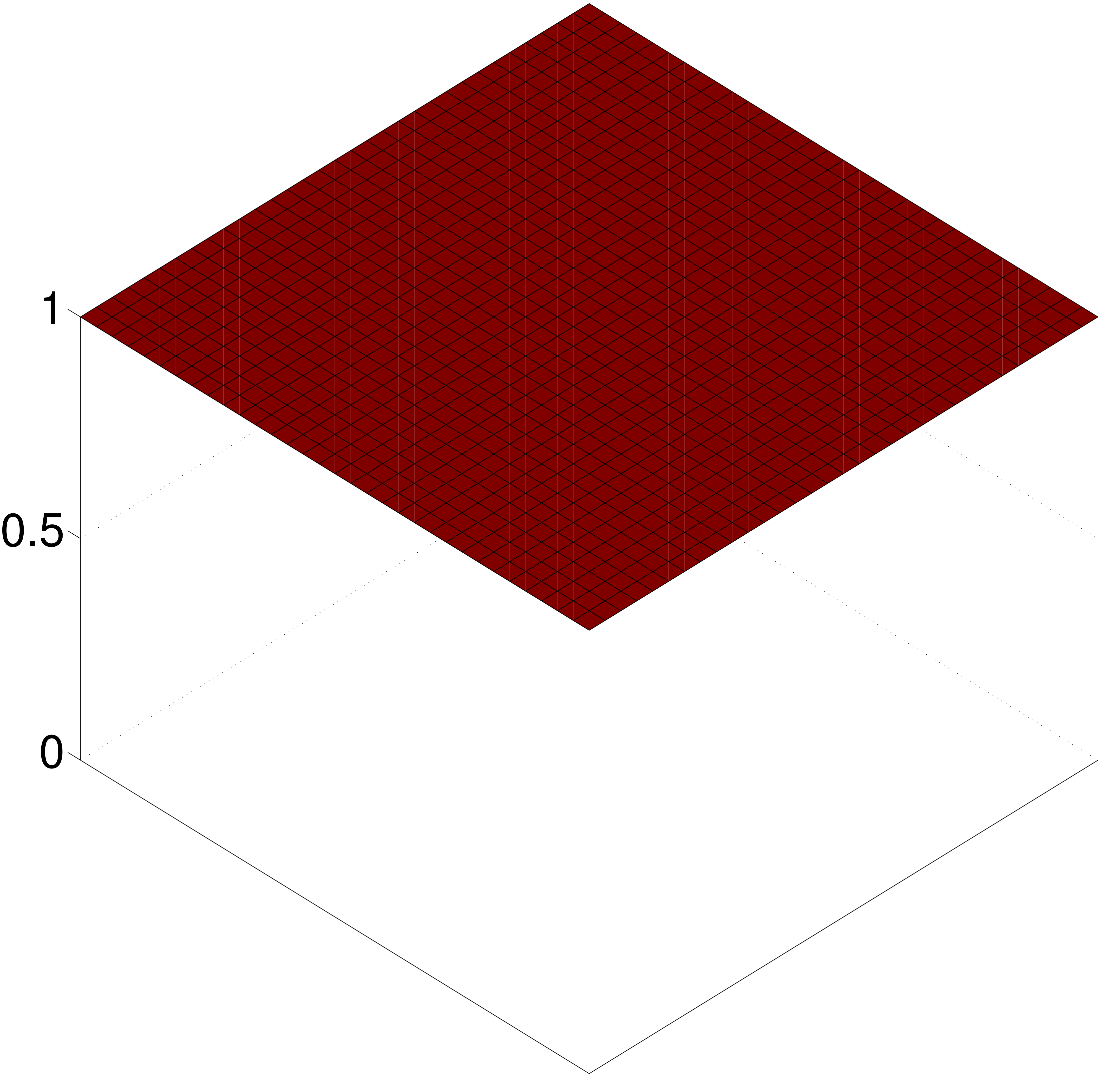}\\
(b): $- \nabla~ .~ ({\bm{\lambda}}  \cdot  \nabla \Phi) = 0$ \\
\caption{Two choices of multiscale basis functions in a reference dual-coarse block (left column), Summation of the basis functions over the dual-coarse block (right column), i.e., partition of unity check.}
\label{Basis-plot}
\end{figure}

The basis functions $\Phi_k$ are assembled over dual-coarse cells $\tilde \Omega^h, \forall h\in\{1, ..., N_d\}$, i.e., $\Phi_k = \bigcup_{h=1}^{N_d} \Phi_{k,h}$, and, if used, the correction functions are also assembled as $\Psi = \bigcup_{h=1}^{N_d} \Psi_{h}^{\nu}$.
Fig.~\ref{Basis-plot} illustrates that the basis functions do not form a partition of unity when compressibility effects are included, \change{which is the intrinsic nature of the parabolic compressible equation}.

The choices formulated above affect computational efficiency of constructing and updating the multiscale operators. More precisely, while basis functions of Eqs.~\eqref{B1} and \eqref{B2} depend on pressure (hence, updated adaptively when pressure changes), Eqs.~\eqref{B3} and \eqref{B4} are pressure independent; thus, they only need to be computed once for single-phase problems (for multi-phase flows, they need to be adaptively updated when local transmissibility changes beyond a prescribed threshold value).
While the basis and correction functions from Eq.~\eqref{B3} were previously used \cite{hadi-comp-jcp}, the other options are, as of yet, have not been studied.

The C-AMS approximates the fine-scale solution $\bm{p}^{\nu}$ by $\bm{p'}^{\nu}$ using the Prolongation operator $\bm{P}$, which is a matrix of size $N_f\times N_c$, having basis function $\Phi_k$ in its $k-th$ column. 
The map between the coarse ($\bm{\breve{p}}$) and fine-scale solution ($\bm{p'}$) reads
\begin{align}\label{superposition}
\bm{p'} = \bm{P \breve{p}}.
\end{align}
The coarse-scale system is obtained using the restriction operator, $\bm{R}$, as 
\begin{align}\label{MS}
\bm{\breve{A}}^{\nu} \bm{\breve{p}}^{\nu+1} \equiv (\bm{R} \bm{A}^\nu \bm{P})  \bm{\breve{p}}^{\nu+1} = \bm{R} \bm{f}^\nu,
\end{align}
and its solution is prolonged to the fine-scale using Eq.~\eqref{superposition}, i.e.,
\begin{align}\label{p-MS}
\bm{p'}^{\nu+1} = \bm{P} (\bm{R} \bm{A}^\nu \bm{P})^{-1} \bm{R} \bm{f}^\nu.
\end{align}
In residual form, it reads
\begin{align}\label{dp-MS}
\bm{\delta p'}^{\nu+1} = \bm{P} (\bm{R} \bm{A}^\nu \bm{P})^{-1} \bm{R} \bm{r}^\nu.
\end{align}
Here $\bm{p'}^{\nu+1} = (\bm{p}^\nu + \bm{\delta p'}^{\nu+1})$, while $\bm{r}^\nu = (\bm{f}^\nu - \bm{A}^\nu \bm{p}^\nu)$ is the fine-scale residual. Note that all the different options for basis functions can be considered in construction of the prolongation operator.

The C-AMS employs Eq.~\eqref{p-MS} as the global solver (for resolving low-frequency errors). In addition to the coarse-scale solver, an efficient convergent multiscale solver needs to include a second-stage smoother at fine  scale. The  smoother accounts for the high-frequency errors, arising from simplified localization conditions, the nonlinearity of the operator, and the complex RHS term.
Among the choices for this smoother (block-, line-, or point-wise solvers), the correction functions (CF) and ILU(0) are considered in this work. The C-AMS procedure is finally summarized in Table \ref{cams}.

\begin{table}[htp]
\begin{tabular}{ll}
Do until convergence ($\| \bm{\epsilon} \|_2^\nu < e$) achieved (See Eq.~\eqref{epsilon_nu}) \{ \\
  \hspace{1mm} 1. \textbf{Initialize:} update linear system components, $\mathbf{A}^{\nu}$ and $\bm{f}^\nu$, based on $\bm{p}^\nu$\\
  \hspace{1mm} 2. \textbf{Update residual:} $\bm{r}^\nu = \bm{f}^\nu - \bm{A}^\nu \bm{p}^\nu$\\
  \hspace{1mm} 3. \textbf{Adaptively compute Basis Functions}: use either of Eqs.~\eqref{B1}-\eqref{B4} \\
  \hspace{1mm} 4. \textbf{Pre-smoothing Stage}: only if CF is used, apply CF on  $\bm{r}^\nu$ and update residual\\
  \hspace{1mm} 5. \textbf{Multiscale Stage:} Solve \eqref{dp-MS} for $\bm{\delta p'}^{\nu+1/2}$\\
  \hspace{1mm} 6. \textbf{Post-smoothing Stage:} smooth $\bm{\delta p'}^{\nu+1/2}$ for $n_s$ times using a fine-scale iterative\\\hspace{51mm}solver (here, ILU(0) is used), obtaining $\bm{\delta p'}^{\nu+1}$ \\
  \hspace{1mm} 7. \textbf{Update solution:} $\bm{p}^{\nu+1} = (\bm{p}^\nu + \bm{\delta p'}^ {\nu +1/2} + \bm{\delta p'}^{\nu+1}) $\\
  \hspace{1mm} 8. \textbf{Update error:} compute $\bm{\epsilon}^\nu$, and assign $\bm{p}^\nu \leftarrow \bm{p}^{\nu+1} $\\
  \}
\end{tabular}
\caption{C-AMS iteration procedure, converging to $\bm{p}^{n+1}$ with tolerance $e$.}
\label{cams}
\end{table}
 
In the next section, numerical results for 3D heterogeneous test cases are presented, in order to provide a thorough assessment of the applicability of C-AMS to large-scale problems.

\section{Numerical Results}

The numerical experiments presented in this section are divided into: (1) finding a proper iterative procedure and multi-stage multiscale components for efficiently capturing the nonlinearity within the flow equation, and (2)  systematic performance study by comparing against a commercial algebraic multigrid solver, i.e., SAMG \cite{SAMG}.
Note that the second aspect is mainly to provide the computational physics community with an accurate assessment of the convergence properties of the state-of-the-art compressible multiscale solver (i.e., C-AMS).
{\change{As an advantage over many advanced linear solvers, C-AMS allows for construction of locally conservative velocity after any MSFV stage. Therefore, for multiphase flow scenarios, only a few C-AMS iterations are necessary to obtain accurate solutions \cite{hadi-erest-spej}.}}

\change{
For the studied numerical experiments of this paper, sets of log-normally distributed permeability fields with spherical variograms are generated by using  sequential Gaussian simulations \cite{SGeMS}. The variance and mean of natural logarithm of the permeability, i.e., $ln(k)$, for all test cases are 4 and -1, respectively, unless otherwise is mentioned. Furthermore, the fine-scale grid size and dimensionless correlation lengths in the principle directions, i.e., $\psi_1$, $\psi_2$ and $\psi_3$, are provided in Table~\ref{tab:permeability set}. Each set has 20 statistically-equivalent realizations. The sets with orientation angle of $15^\circ$ are referred to as the layered fields. Also, the grid aspect ratio $\alpha$ is 1, i.e., $\Delta x/\alpha = \Delta y = \Delta z = 1$ m, unless otherwise is specified.
\begin{table}[htp!] 
\begin{center}
\begin{tabular}{c|c|c|c||c|c|c}
\hline
\hline 
\change{Permeability Set} & \change{1}  & \change{2} & \change{3} & \change{4} & \change{5} & \change{6}\\
\hline 
\change{Fine-scale grid} & \change{$64^3$} & \change{$128^3$}  & \change{$256^3$} & \change{$64^3$} & \change{$128^3$} & \change{$256^3$}\\
\hline
\change{$\psi_1$} & \change{0.125} & \change{0.125} & \change{0.125} & \change{0.5} & \change{0.5}  & \change{0.5}\\
 \hline
\change{$\psi_2$} &\change{0.125} & \change{0.125} & \change{0.125}& \change{0.03} & \change{0.03} & \change{0.03}\\
\hline
\change{$\psi_3$} &\change{0.125} & \change{0.125} & \change{0.125} & \change{0.01} & \change{0.06} & \change{0.01}\\
\hline
\change{Angle between $\psi_1$ and $y$ direction }& \multicolumn{3}{|c||}{\change{patchy}} & \multicolumn{3}{|c}{\change{$15^\circ$}} \\
\hline
\change{Variance of ln(k)} & \multicolumn{5}{|c}{\change{4}} \\
\hline
\change{Mean of ln(k)} & \multicolumn{5}{|c}{\change{-1}} \\
\hline
\hline
\end{tabular}
\caption{\change{Permeability sets (each with 20 statistically-equivalent realizations) used for numerical experiments of this paper. Layered fields refer to the sets 4-6, in which the orientation angle between $\psi_1$ and y direction is $15^\circ$.}}
\label{tab:permeability set}
\end{center}
\end{table}
}

\change{Phase properties and simulation time are described as non-domensional numbers. The non-dimensional pressure and density are introduced as
\begin{equation}
p^* = \frac{p-p_{east}}{p_{west}-p_{east}},
\end{equation}
and
\begin{equation}
\rho^* = \frac{\rho}{\rho_0} = 1 +  \eta \ p^*,
\end{equation}
respectively, where the coefficient $\eta$ is set to $1$ for all subsequent test cases in this paper. }

\change{The $p_{west}$ and $p_{east}$ values of $10^6$ and $0$ Pa, relative to the Standard (Atmospheric) condition, are considered. These correspond to non-dimensional pressure values of $1$ and $0$, which are set as Dirichlet conditions at the west and east boundaries, respectively, for all the cases unless otherwise is mentioned.} Also, all the other surfaces are subject to no-flow Neumann conditions.

\change{The non-dimensional time is introduced as $t^* = t/\tau$, where
\begin{equation}
\tau = \frac{\mu \phi L^2}{\bar{K}(p_{west} - p_{east})}.
\end{equation}
Here,} $\bar{K}$ is the average permeability, and $L$ is a length scale of the domain. \change{With the values of $10^6$ Pa pressure difference, in-situ viscosity of $2\times 10^{-6}$ Pa.s, $\Delta x = 1$ m, $\phi = 0.1$, and $\bar{K}$ value of $10^{-12}$ m$^2$ for homogeneous cases, the $\tau$ will be $128$ s for problem size of $L =64$ m in SI units.}

The implementation used to obtain the results presented in this paper consists of a single-threaded {\change{object-oriented C++ code}}, and the CPU times were measured on an Intel Xeon E5-1620 v2 quad-core system with 64GB RAM.

\subsection{C-AMS: determining the most effective iterative procedure and multi-stage multiscale components}

\begin{figure} [htb!]
\centering
\includegraphics[width=0.4\textwidth]{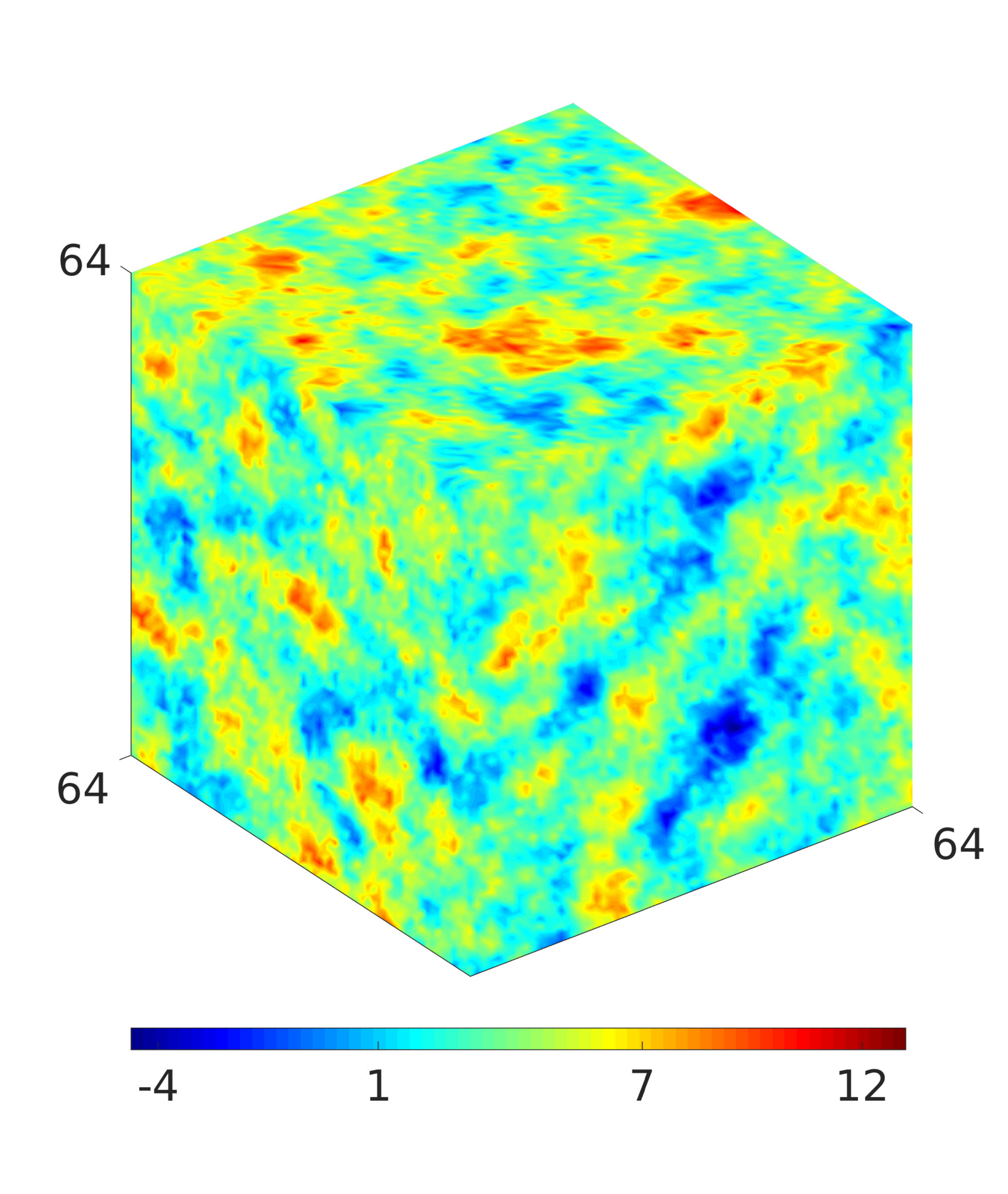}\hspace{2cm}
\includegraphics[width=0.4\textwidth]{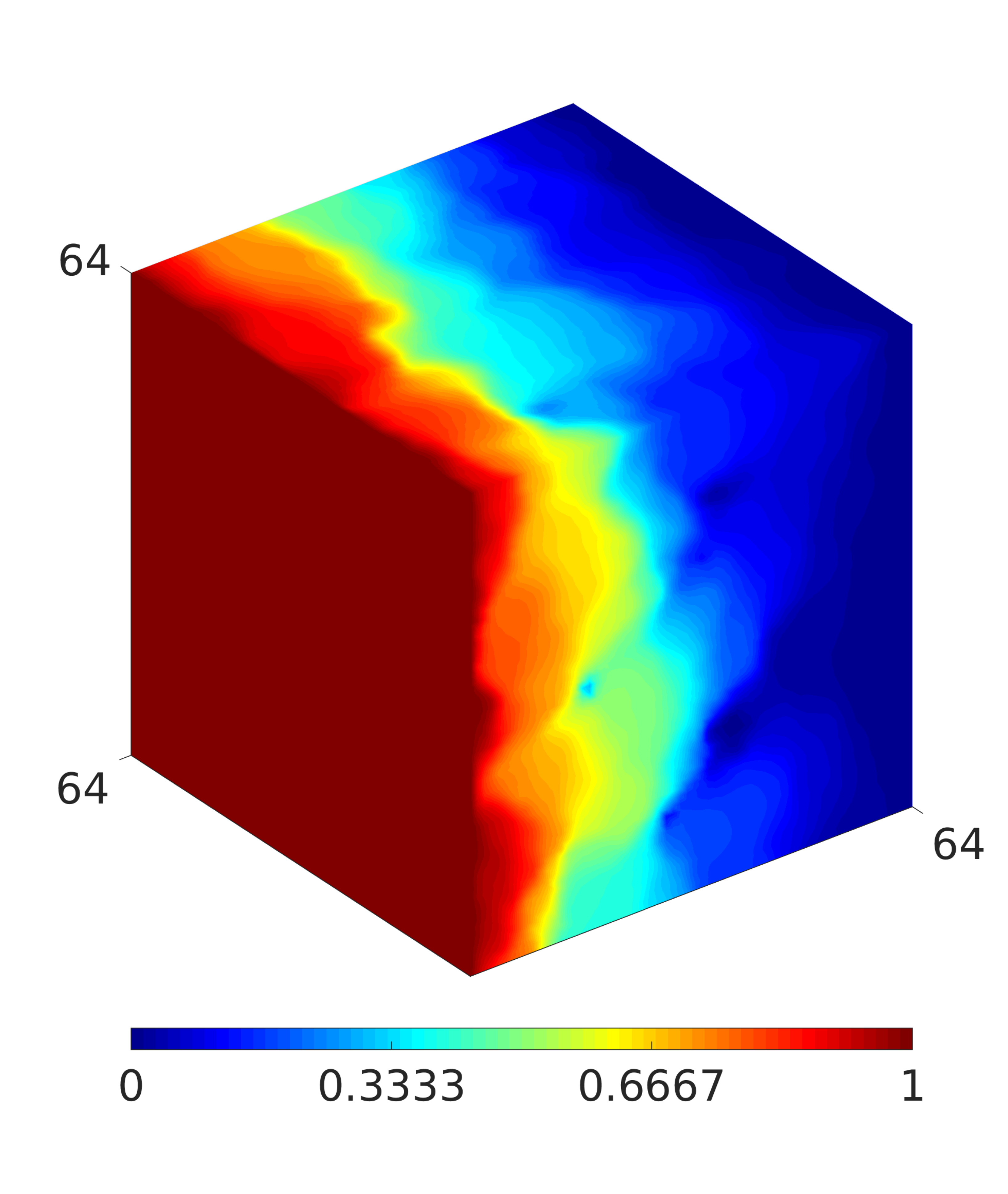}
\caption{Natural-log  of the permeability (left) and pressure solution after \change{$t^* = 0.4$} (right) corresponding \change{to one of the realization of permeability Set 1 from Table \ref{tab:permeability set}}.}
\label{fig:realiz}
\end{figure}

The efficient capturing of the nonlinearity within the iterations is important in designing an efficient multiscale strategy. For the purposes of a conclusive result, in this section, a set of 20 statistically\change{-equivalent} patchy fields\change{, i.e., permeability Set 1 from Table \ref{tab:permeability set},} is considered. One of the realizations and its corresponding solution at $t^* = 0.4$ are shown in Fig.~\ref{fig:realiz}.

\subsubsection{Nonlinear and linear level updates}
\label{sec:sol_strat}

In formulating a convergence criterion for the C-AMS, one can express the error of the approximate solution at step $\nu$ on the basis of either the linear or nonlinear expressions. According to Eq.~\eqref{p-nonlinear}, the nonlinear error in each grid cell reads
\begin{align}\label{epsilon_nu}
\epsilon^{\nu} = q -  \frac{\phi^{n+1}}{\Delta t} + \frac{\phi^{n}}{\Delta t} \frac {{\rho^n
S_{\alpha}^n}}{\rho^{n+1}} + \frac {1}{\rho^{n+1}}{\nabla \cdot \big(\rho^{n+1} \bm \lambda \cdot \nabla p^{n+1} \big)},
\end{align}
and is assembled in the vector $\bm{\epsilon}^\nu$, which allows the computation of the \textit{error norm}, $\| \bm{\epsilon}^{\nu} \|_2$. On the other hand, the linear-level error is based on the linearized equation \eqref{lin_p}, which leads to the computation of the \textit{residual norm}, $\| \bm{r}^{\nu} \|_2$.

In order to determine a suitable sequence of the linear and nonlinear stages, the same patchy domain of $64 \times 64 \times 64$ grid cells is considered (Fig.~\ref{realiz}), for which the pressure equation is solved using the following solution strategy:
\begin{table}[htp]
\begin{tabular}{ll}
Do until ($\| \bm{\epsilon} \|_2 < 10^{-6} $) is reached \{ \\
  \hspace{1mm} 0. Update parameters, linear system matrix and RHS vector based on $\bm{p}^\nu$\\
  \hspace{1mm} 1. Solve linear system using the Richardson iterative scheme, preconditioned\\
    \hspace{5mm}  with one multigrid V-cycle until $\| \bm{r} \|_2 < 10^{-6}$ \\
  \}
\end{tabular}
\caption{Solution strategy used to determine a suitable stopping criterion}
\label{tab:conv}
\end{table}

\begin{figure} [htb!]
\centering
\includegraphics[width=0.47\textwidth]{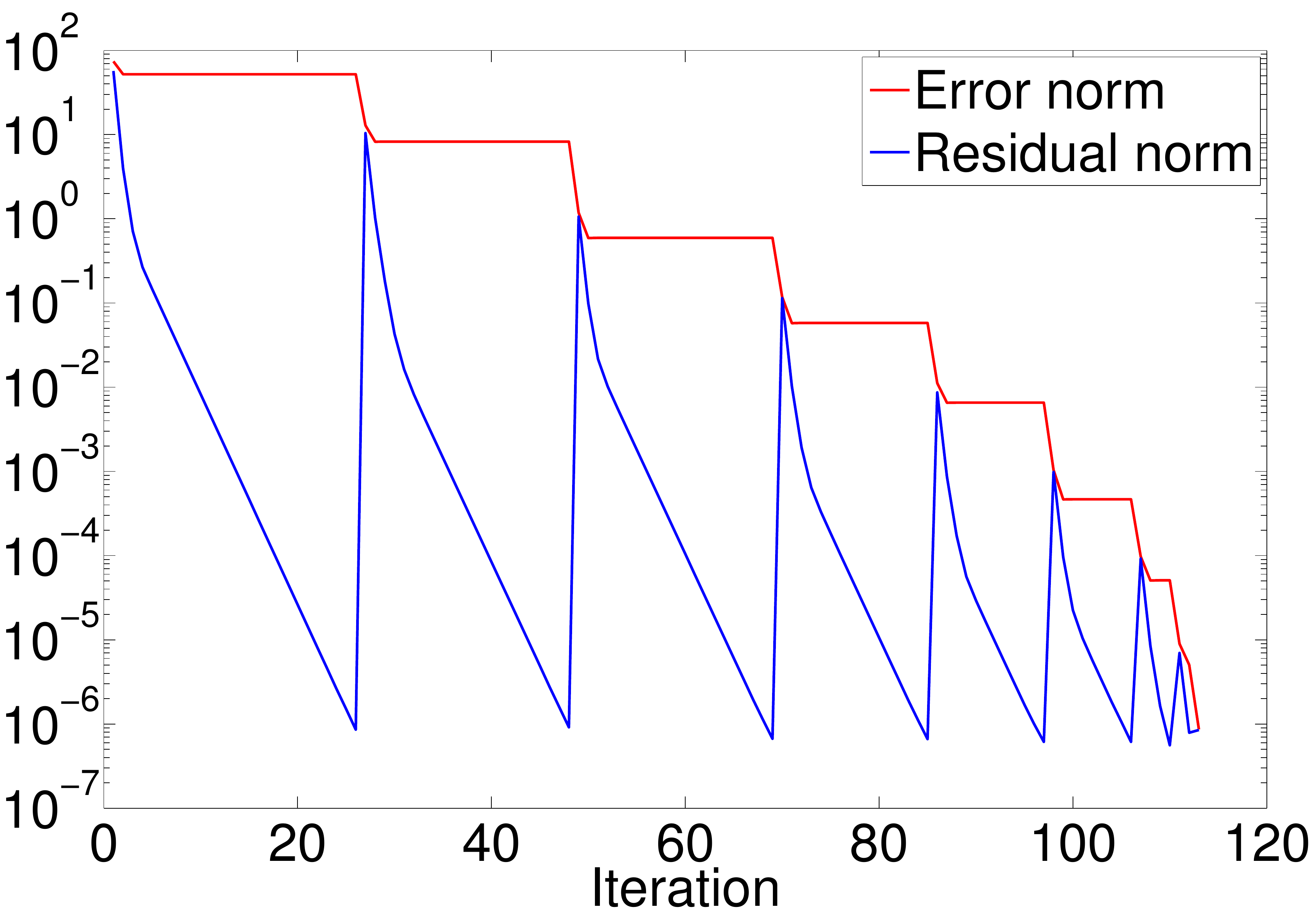} \hspace{0.5cm}
\includegraphics[width=0.47\textwidth]{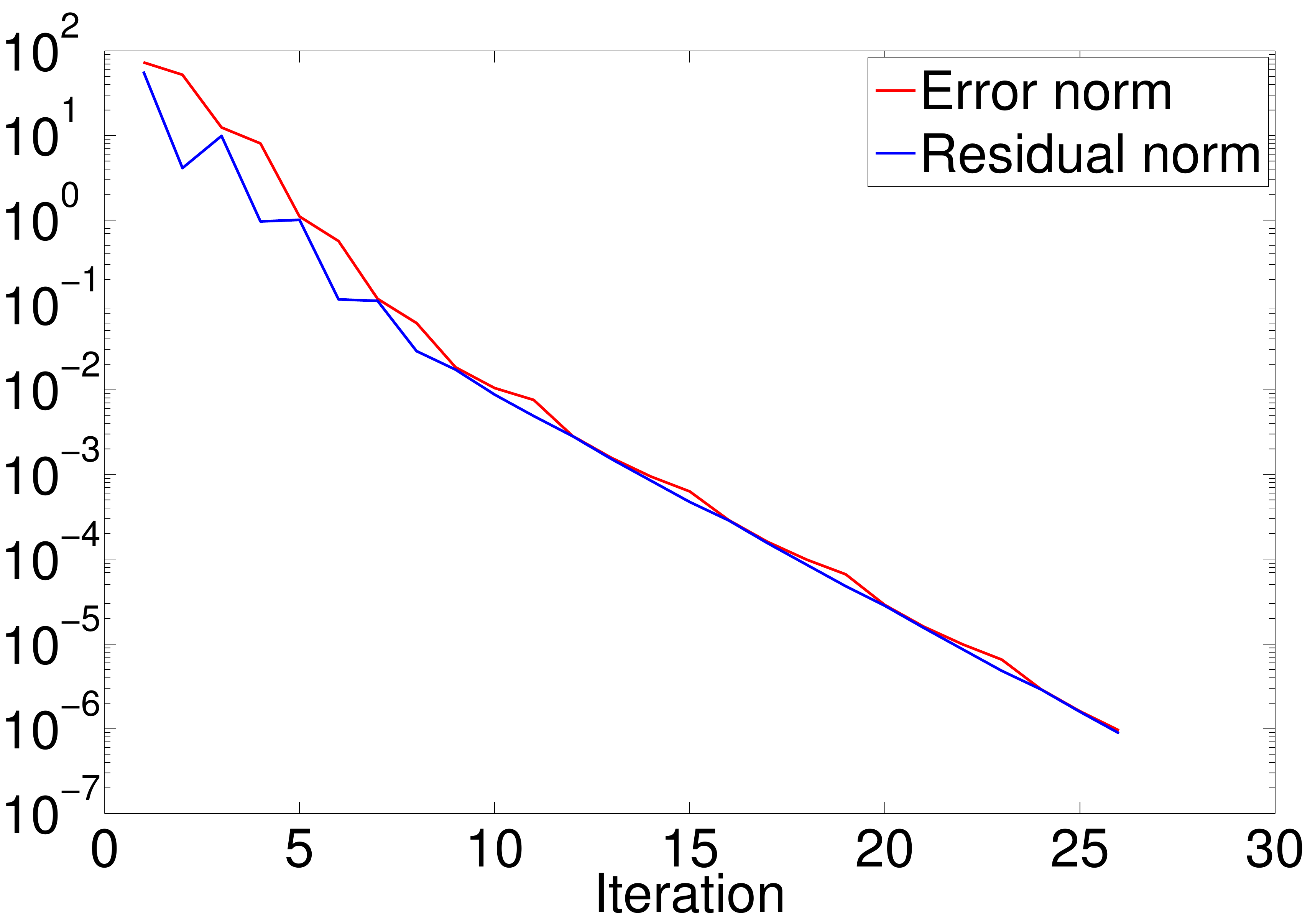}
\begin{picture}(360,0)
\put(70, 0){\change{(\bf a)}}
\put(280, 0){\change{(\bf b)}}
\end{picture}
\caption{Error and residual norm histories for one of the realizations of \change{permeability Set 1 from Table \ref{tab:permeability set}} over a single time step of \change{$t^* = 0.4$}. Shown on the left is the strategy where at each nonlinear stage, the fully converged linear solution is obtained. Shown on the right is the strategy where in each outer (nonlinear) loop the residual is reduced only by one order of magnitude.}\label{fig:conv_norm}
\end{figure}
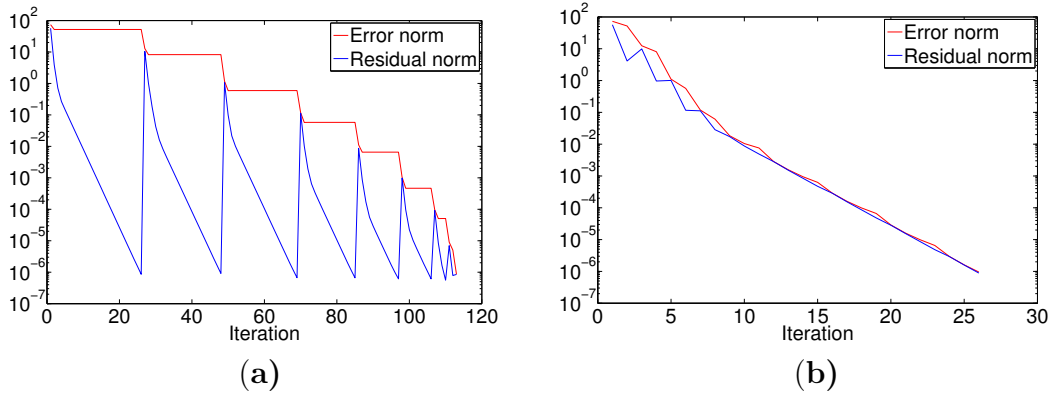

The error and residual norms were recorded after each iteration of the Richardson loop and are presented in Fig.~\ref{fig:conv_norm}. Note that the reduction of the residual norm beyond the first few iterations does not contribute to the reduction in the (nonlinear) error norm. Therefore, one could ideally speed up the solution scheme by monitoring the error norm and updating the linear system after its decrease starts to stagnate. However, the computational cost of evaluating the nonlinear equation is roughly the same as that of a linear system update and, thus, much more expensive than the evaluation of the residual norm.

Fig.~\ref{fig:conv_norm}{\change{{(a)}}} also reveals that the stagnation of the error norm happens roughly after the residual norm has been approximately reduced by $1/10$ of its initial value (i.e., immediately after the linear system update). Fig.~\ref{fig:conv_norm}{\change{{(b)}}} shows the convergence behaviour after implementing this heuristic strategy, which is deemed quite efficient, since the two norms are in agreement. Hence, in the following experiments the same strategy is employed, i.e., for linear level, $\displaystyle \frac{\| \bm{r}^i \|_2}{\| \bm{r}^0 \|_2} < 10^{-1}$ after iteration $i$ of the inner (linear) loop and, for nonlinear level, $\| \epsilon^{\nu} \|_2 < 10^{-6}$ after iteration $\nu$ of the outer (nonlinear) loop are set (see Table~\ref{tab:conv}).

\subsubsection{Adaptive updating of multiscale operators}

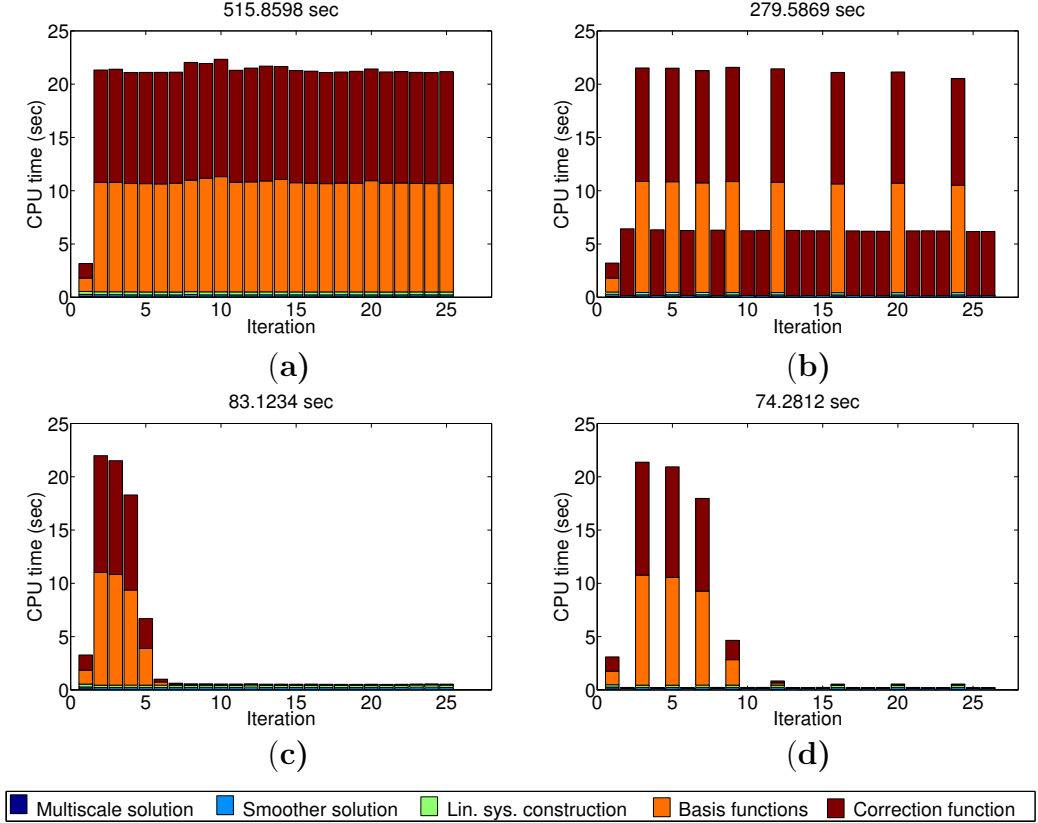
\begin{figure} [htb!]
\centering
\includegraphics[width=0.45\textwidth]{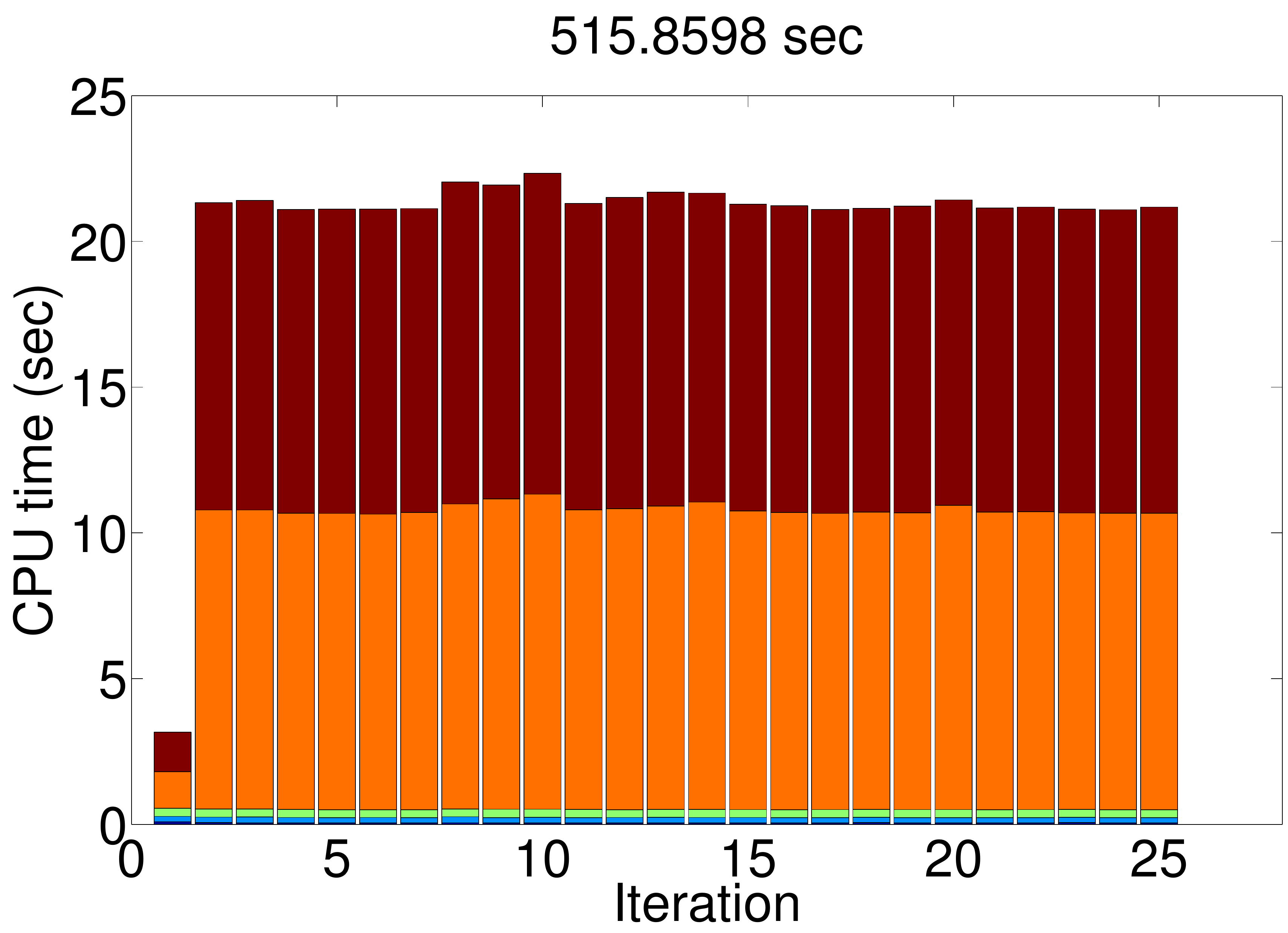}\hspace{0.6cm}
\includegraphics[width=0.45\textwidth]{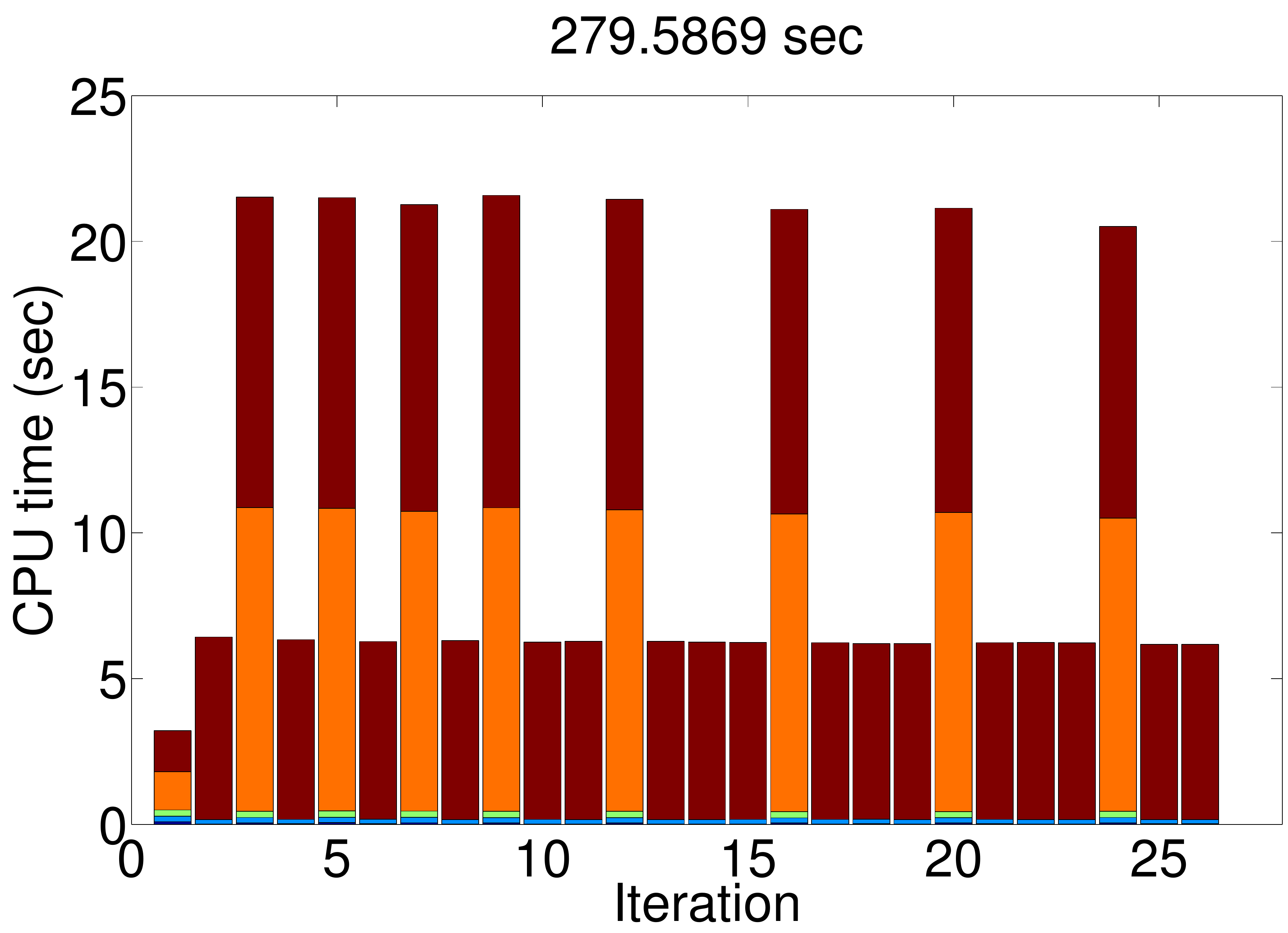}\vspace{0.65cm}
\includegraphics[width=0.45\textwidth]{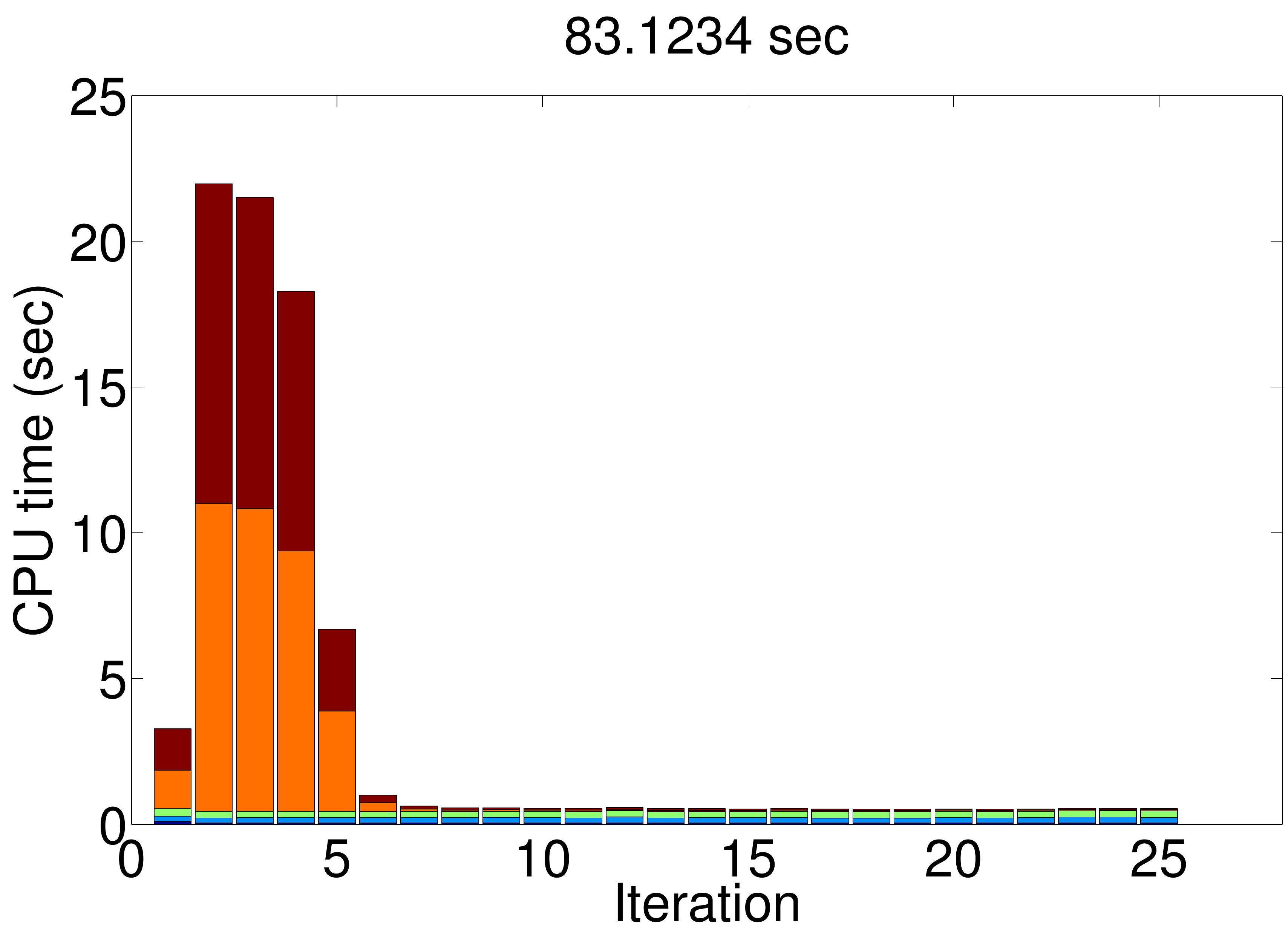}\hspace{0.6cm}
\includegraphics[width=0.45\textwidth]{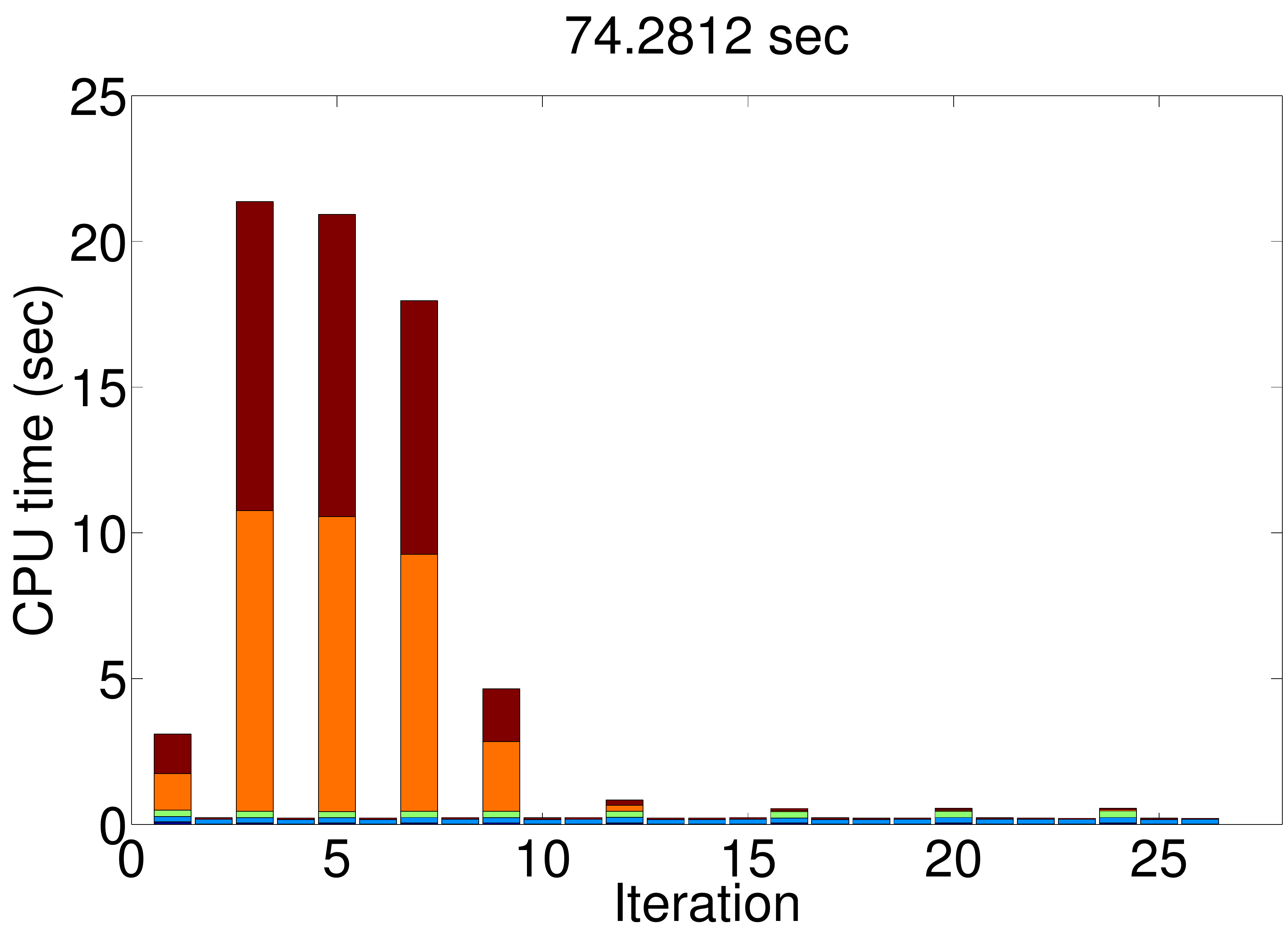}\vspace{0.65cm}
\includegraphics[width=\textwidth]{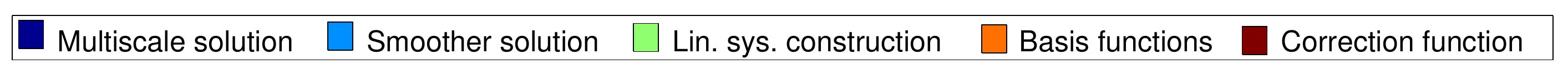}
\begin{picture}(360,0)
\put(84,187){(\bf a)}
\put(280,187){(\bf b)}
\put(84,40){(\bf c)}
\put(280,40){(\bf d)}
\end{picture}
\caption{Effect of different types of adaptivity on the C-AMS performance for \change{the permeability Set 1 from Table \ref{tab:permeability set}} after a time step of \change{$t^* = 0.4$}: (a): No adaptivity, (b): Linear system update adaptivity only, (c): Multiscale operator update adaptivity only, (d):  Fully adaptive, i.e in terms of both linear system and multiscale operator updates.} \label{fig:adapt}
\end{figure}

The previous study described the first adaptive aspect considered in this work, namely, updating the linear system only after the residual norm drops by an order of magnitude. The C-AMS procedure can be further optimized by employing adaptive updates of its multiscale components, i.e., the basis and (if considered) the correction functions. To this end, one has to monitor the changes in the entries of the transmissibility matrix $\bm{A}$ and RHS $\bm{f}$ between the  iteration steps. Fig.~\ref{fig:adapt}(c) shows that the adaptive update of the C-AMS basis functions leads to a significant speed-up in terms of CPU time.

Furthermore, the two adaptivity methods (for linear system and local function updates) are combined and shown in Fig. \ref{fig:adapt}(d).
Hence, C-AMS will perform its iterations such that it exploits all adaptivity within the multiscale components and the nonlinearity within the flow equation. Note that for this case, the compressible variant from Eq.~\eqref{B1} was used for both basis and correction functions. However, if the incompressible Eqs~\eqref{B3} and \eqref{B4} are used, then the basis functions do not require updates during iterations. Finally, for this and all the following results (unless otherwise stated), the C-AMS coarsening ratio was taken as $8 \times 8 \times 8$ \change{, because it was found efficient (see Subsection~\ref{sec:coarsening})}.

\subsubsection{C-AMS global stage: choice of basis functions}

The aim of this study is to determine an optimum choice for the type of basis functions for the C-AMS algorithm. The correction function is computed based on Eq.~\eqref{B1} in all cases (and, hence, updated adaptively with pressure), 20 iterations of ILU(0) are used for smoothing and all possibilities for the basis functions, i.e., Eqs.~\eqref{B1}-\eqref{B4}, are considered. Finally, there is a single time step in the simulation, which takes the initial solution at time 0 ($p^*_0 = 0$ everywhere) to the solution at time \change{$t^* = 0.4$}.

\begin{figure} [htb!]
\centering
\includegraphics[width=1\textwidth]{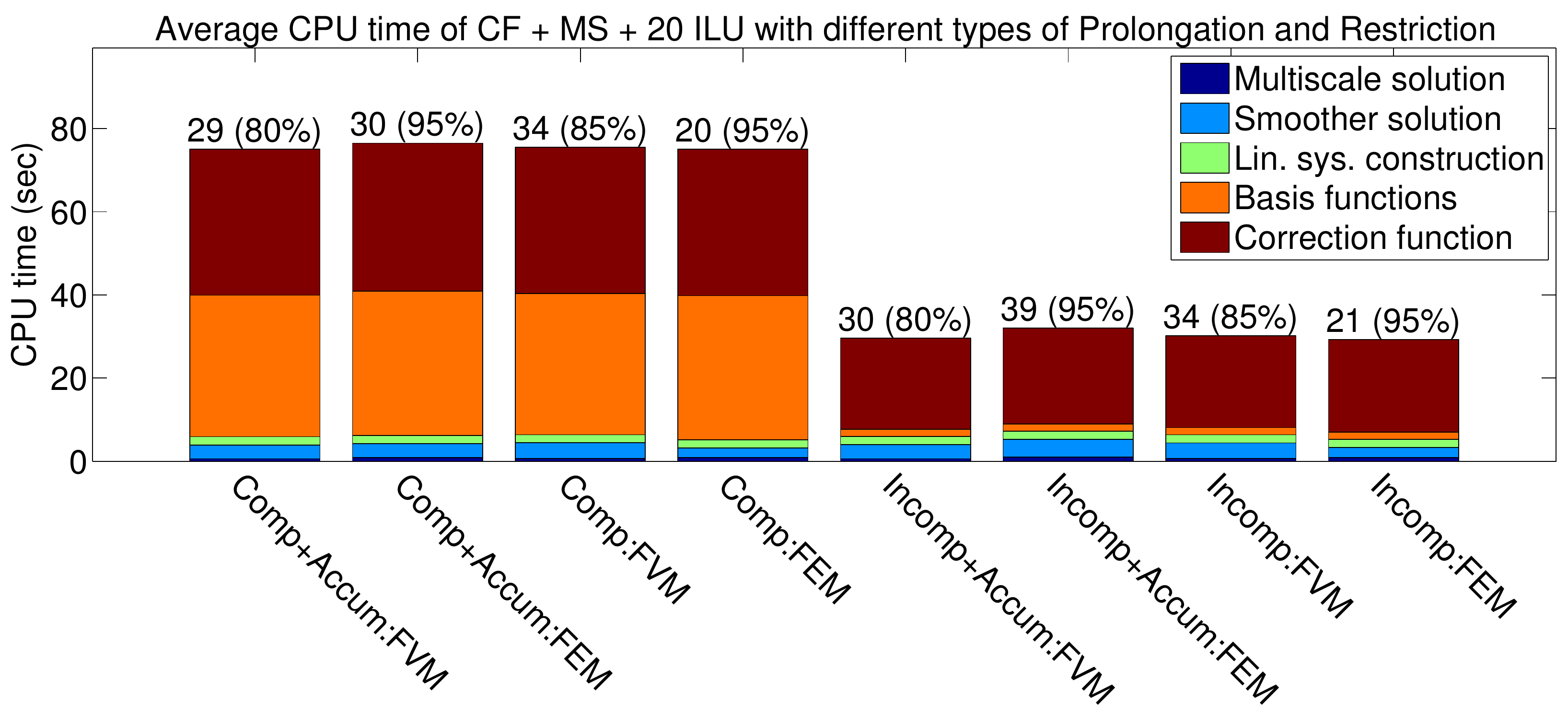}
\caption{Effect of the choice of basis function on the C-AMS performance for the $64^3$ grid-cell problem after a time step of \change{$t^* = 0.4$}. Results are averaged over 20 \change{statistically-equivalent} realizations. The number of iterations is shown on top of each bar. The success percentage is also shown in parentheses. Note that all simulations employ correction functions.}
\label{fig:basis}
\end{figure}

The total CPU time spent in each stage of the solver, as well as the number of iterations (given on top of each bar in Fig.~\ref{fig:basis}), are measured. Also, the success rate of convergence is given inside parentheses beside the average number of iterations. 

The results show that including compressibility in the basis functions does not translate into faster convergence and, thus the additional CPU time required to adaptively update them is not justified. In fact, it is more efficient to use the incompressible (pressure independent) basis functions from Eqs.~\eqref{B3} and \eqref{B4}. Also, the inclusion of the accumulation term and the type of Restriction (MSFE or MSFV) does not play an important role for this patchy test case. Note that none of the choices results in 100\% successful convergence{\change{, even though 20 ILU(0) smoothing iterations have been employed at each iteration}}. This can be attributed to the use of correction functions, as investigated in the next paragraph.

\subsubsection{C-AMS smoothing stage: choice of correction function}

Note that none of the results from the previous test case (Fig.~\ref{fig:basis}) has a $100\%$ success rate. As described in \cite{yixuan-ams}, the CF can be seen as an independent stage, the inclusion of which should be seen as an option and not a necessity for convergence. Fig.~\ref{fig:corr} presents the results of rerunning the previous experiment, this time varying the type of correction function. The plot confirms that eliminating the CF altogether leads to an overall speed-up, and, in addition, a convergence success rate of $100\%$. {\change{As described in \cite{yixuan-ams},} this can be explained by the sensitivity of CF to the heterogeneity of the permeability field, which leads to solver instability. Therefore, the CF should not be considered as candidate for the pre-smoothing stage in an efficient C-AMS procedure. Instead, ILU(0) is performed as post-smoother in order to resolve high-frequency errors. 

\begin{figure} [htb!]
\centering
\includegraphics[width=1\textwidth]{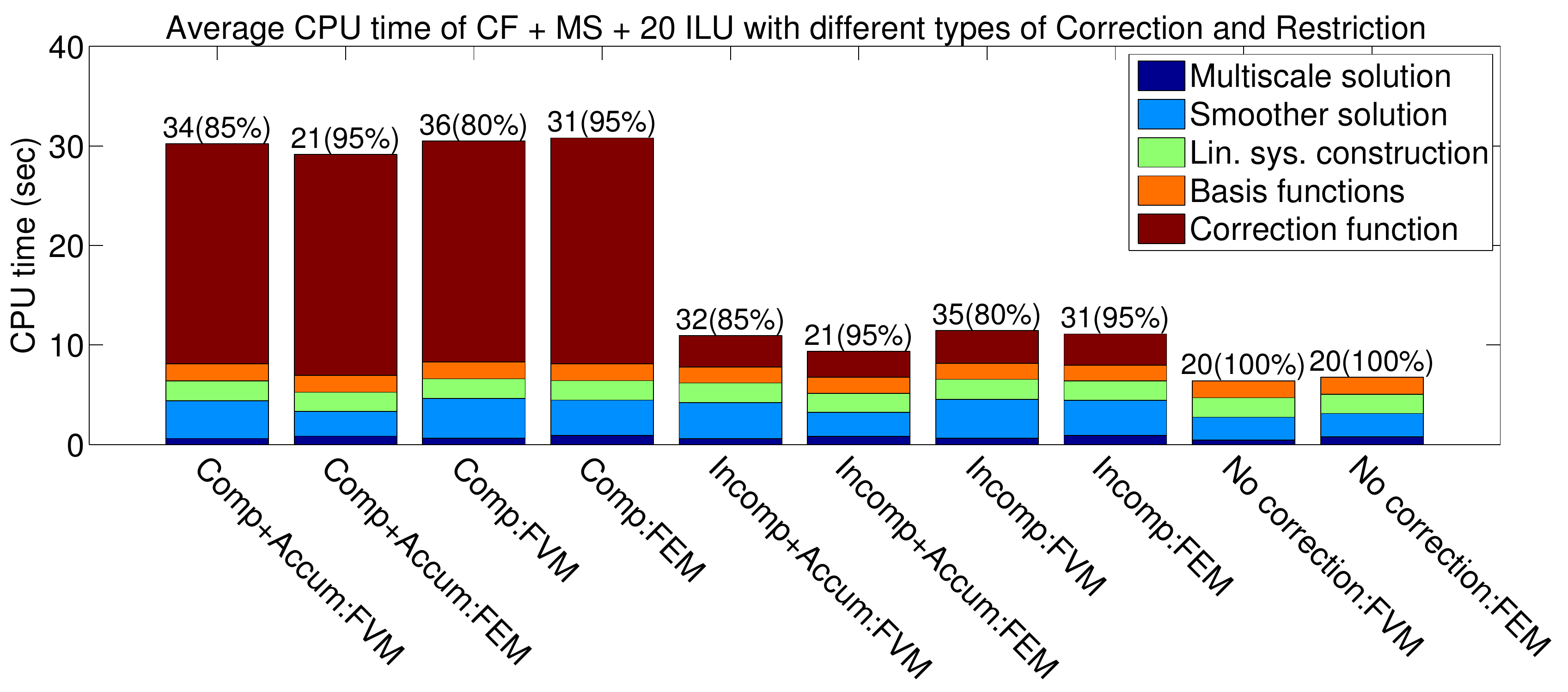}
\caption{Effect of the choice of correction function on the CPU time of the multiscale solution on \change{the permeability Set 1 from Table \ref{tab:permeability set}} after a time step of \change{$t^* = 0.4$}. The number of iterations is shown on top of each bar. Only the last 2 bars on the right correspond to runs in which no correction function was used (i.e., MS + 20 ILU).}
\label{fig:corr}
\end{figure}

\subsubsection{C-AMS smoothing stage: number of smoother iterations}

Another variable in the C-AMS framework is the number of smoothing steps (here, ILU(0)) that should be applied in order to obtain the best trade-off between convergence rate and CPU time. The results of several experiments with the optimum choices (i.e., incompressible basis functions and no incorporation of CF) and various numbers of ILU applications are illustrated in Fig.~\ref{fig:ilu}. It is clear that with this C-AMS setup, an optimum scenario would be found with 5-10 ILU iterations per second-stage call. Note that all C-AMS runs (without correction functions) converged successfully.

\begin{figure} [htb!]
\centerline{\includegraphics[width=0.93\textwidth]{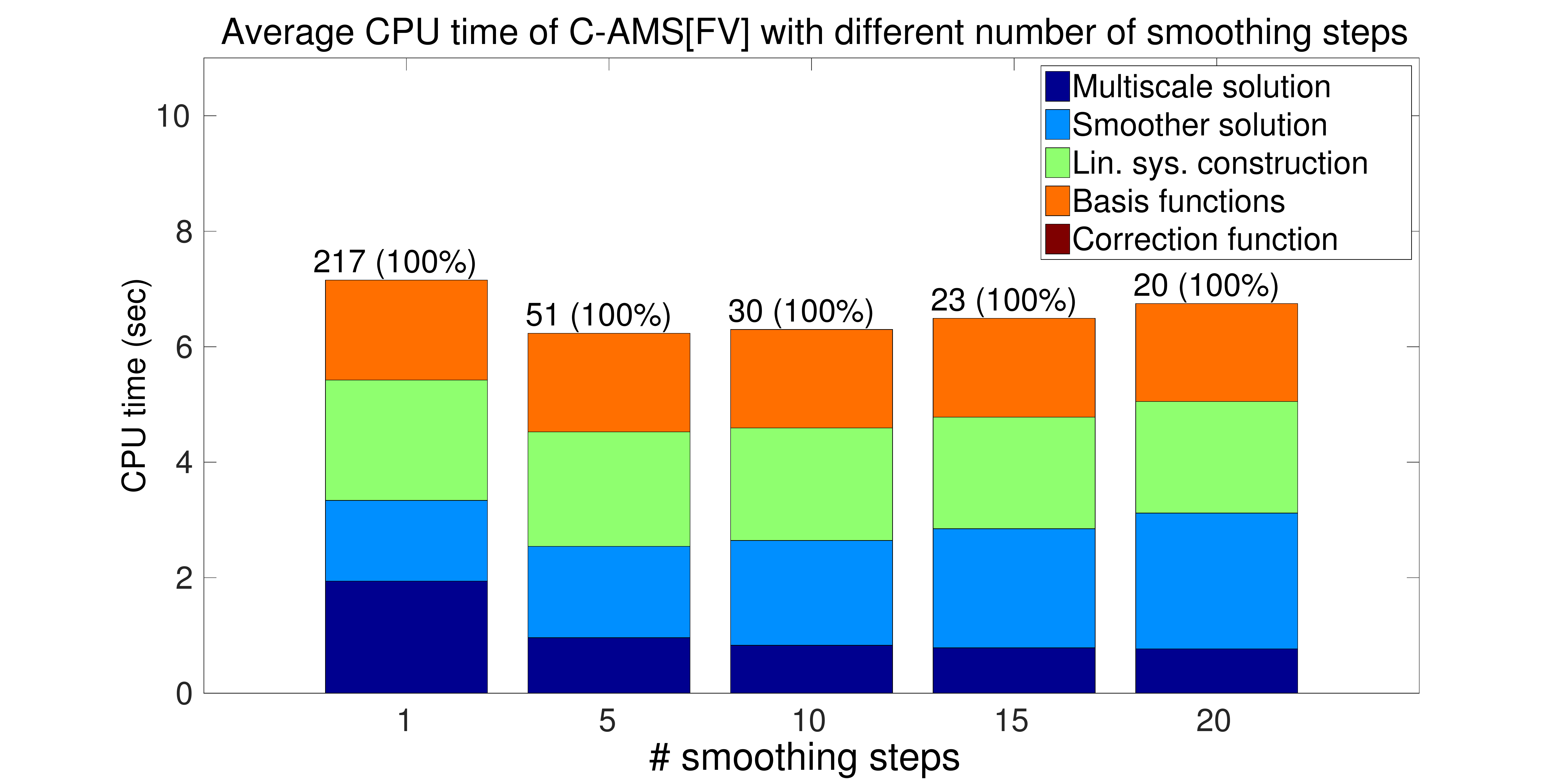}}
\caption{Effect of the number of \change{ILU} smoothing steps on the C-AMS[FV] performance for \change{the permeability Set 1 from Table \ref{tab:permeability set}} \change{(grid aspect ratio is $1$)} after a time step of \change{$t^* = 0.4$}. The number of iterations is shown on top of each bar, with convergence success rate inside parentheses. Note that excluding CF leads to 100\% success rate for all scenarios.}
\label{fig:ilu}
\end{figure}

\subsubsection{C-AMS sensitivity to coarsening ratio: trade-off between size of coarse system and local problem cost}
\label{sec:coarsening}

The coarsening factors used in this paper were found to be optimal after a careful study of the C-AMS sensitivity with the coarsening ratio. As for a thorough study of the new C-AMS solver, it is important to illustrate also its sensitivity with change of coarse-scale system size (and thus the coarsening ratio). This important fact is studied and shown in Figs.~\ref{coarse_64}-\ref{coarse_256} for patchy fields. \change{Not that for the cases studied in this paper, the optimum overall CPU times were obtained with coarse-grid cells with the size of (approximately) the square-root of the domain length in each direction.} 

\begin{figure} [htb!]
\centerline{\includegraphics[width=0.85\textwidth]{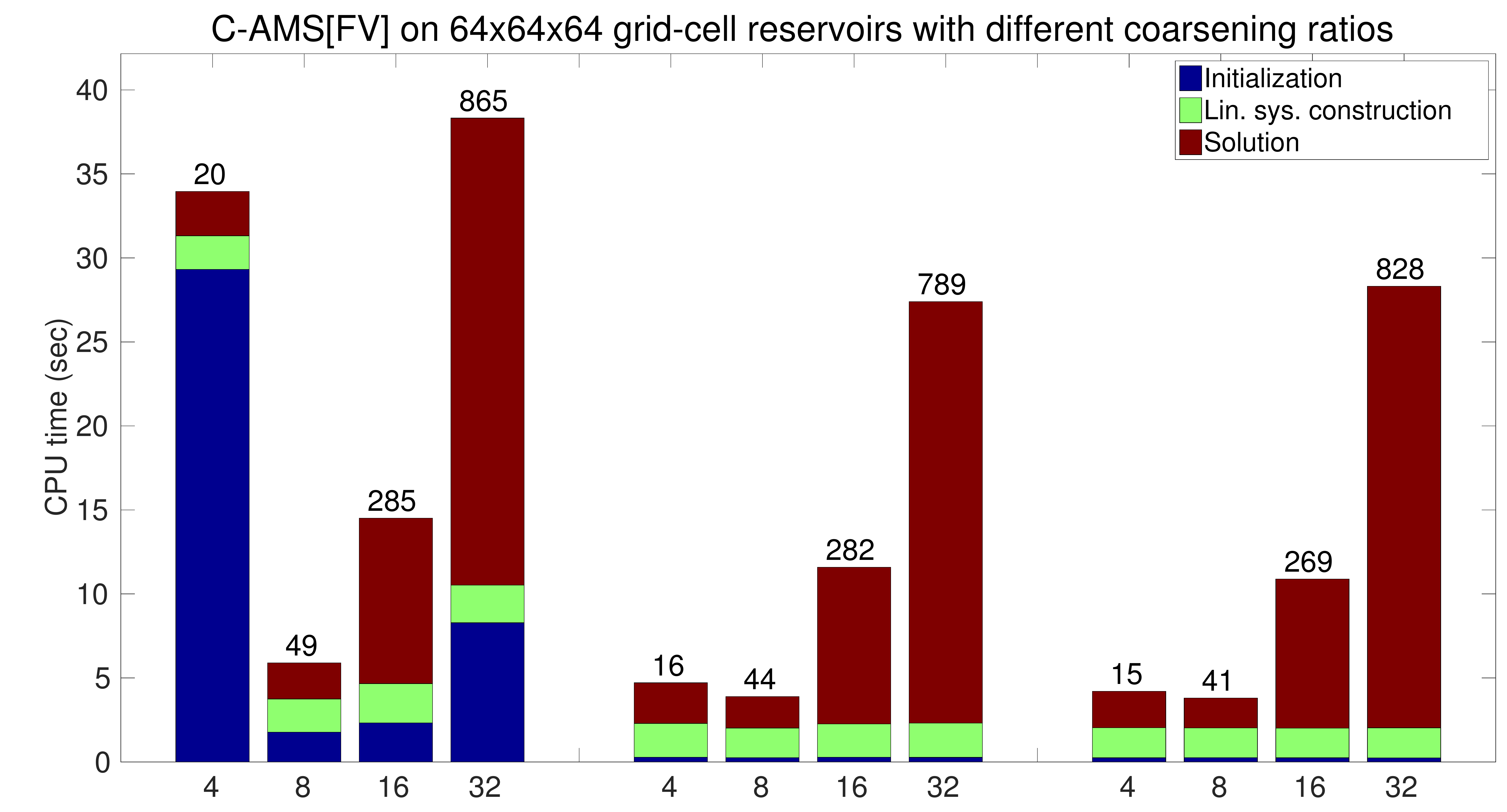}}
\begin{picture}(360,0)
\put(60,5){$t^*:$}
\put(90,5){0.0-0.4}
\put(190,5){0.4-1.0}
\put(295,5){1.0-2.0}
\end{picture}
\caption{Patchy fields: Averaged CPU time (over 20 realizations) of C-AMS[FV] for different coarsening ratios for \change{the permeability Set 1 from Table \ref{tab:permeability set}}. Results support the use of coarsening ratio of $8^3$. \change{A similar behaviour was observed with the FE restriction operator.}}\label{coarse_64}
\end{figure}

\begin{figure} [htb!]
\centerline{\includegraphics[width=.85\textwidth]{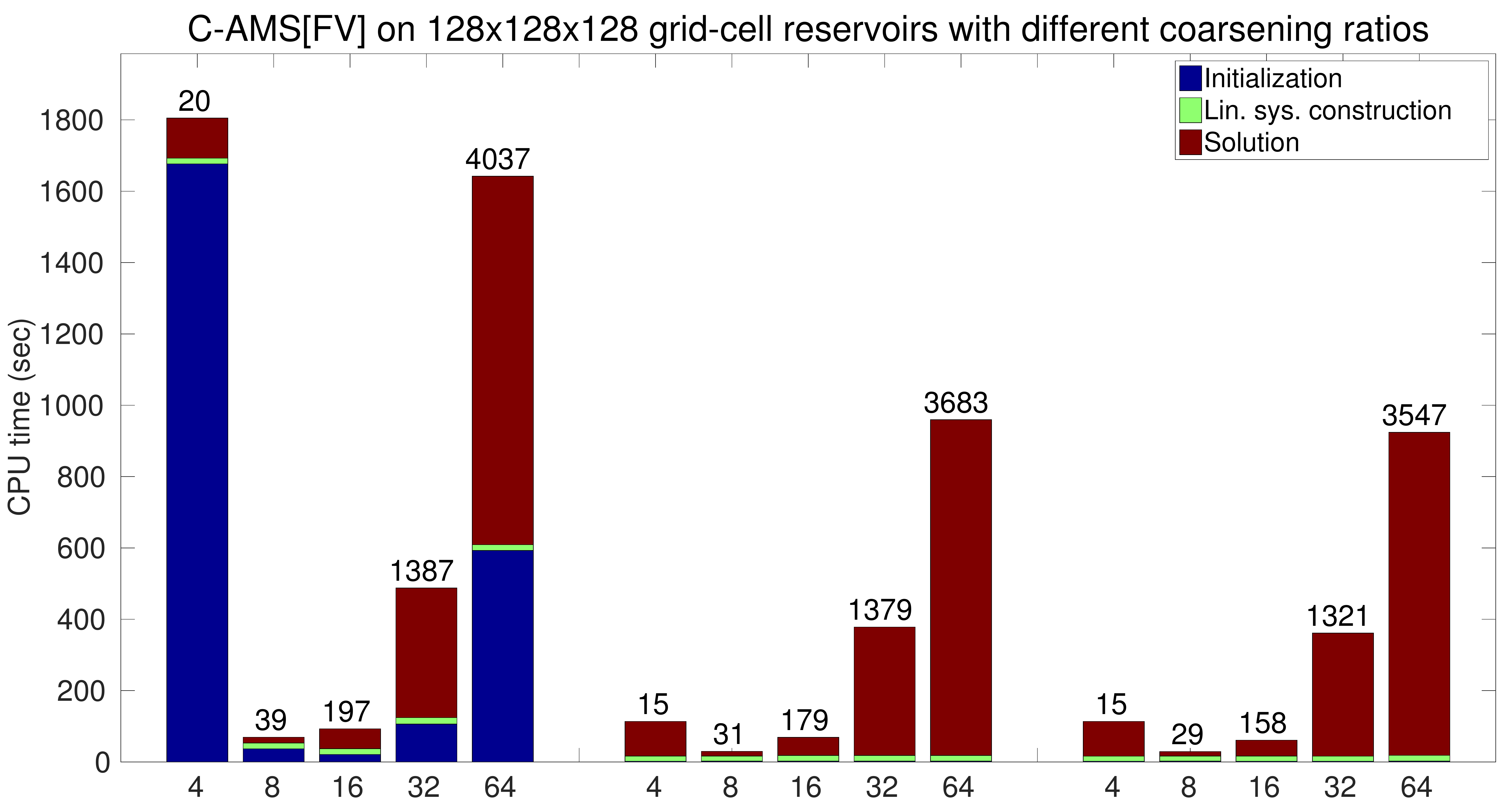}}
\begin{picture}(360,0)
\put(60,5){$t^*:$}
\put(90,5){0.0-0.4}
\put(200,5){0.4-1.0}
\put(300,5){1.0-2.0}
\end{picture}
\caption{Patchy fields: Averaged CPU time (over 20 realizations) comparison of C-AMS [FV] for different coarsening ratios for \change{the permeability Set 2 from Table \ref{tab:permeability set}}. Results support the use of coarsening ratio of $8^3$. \change{A similar behaviour was observed with the FE restriction operator.}}\label{coarse_128}
\end{figure}

\begin{figure} [htb!]
\centerline{\includegraphics[width=.85\textwidth]{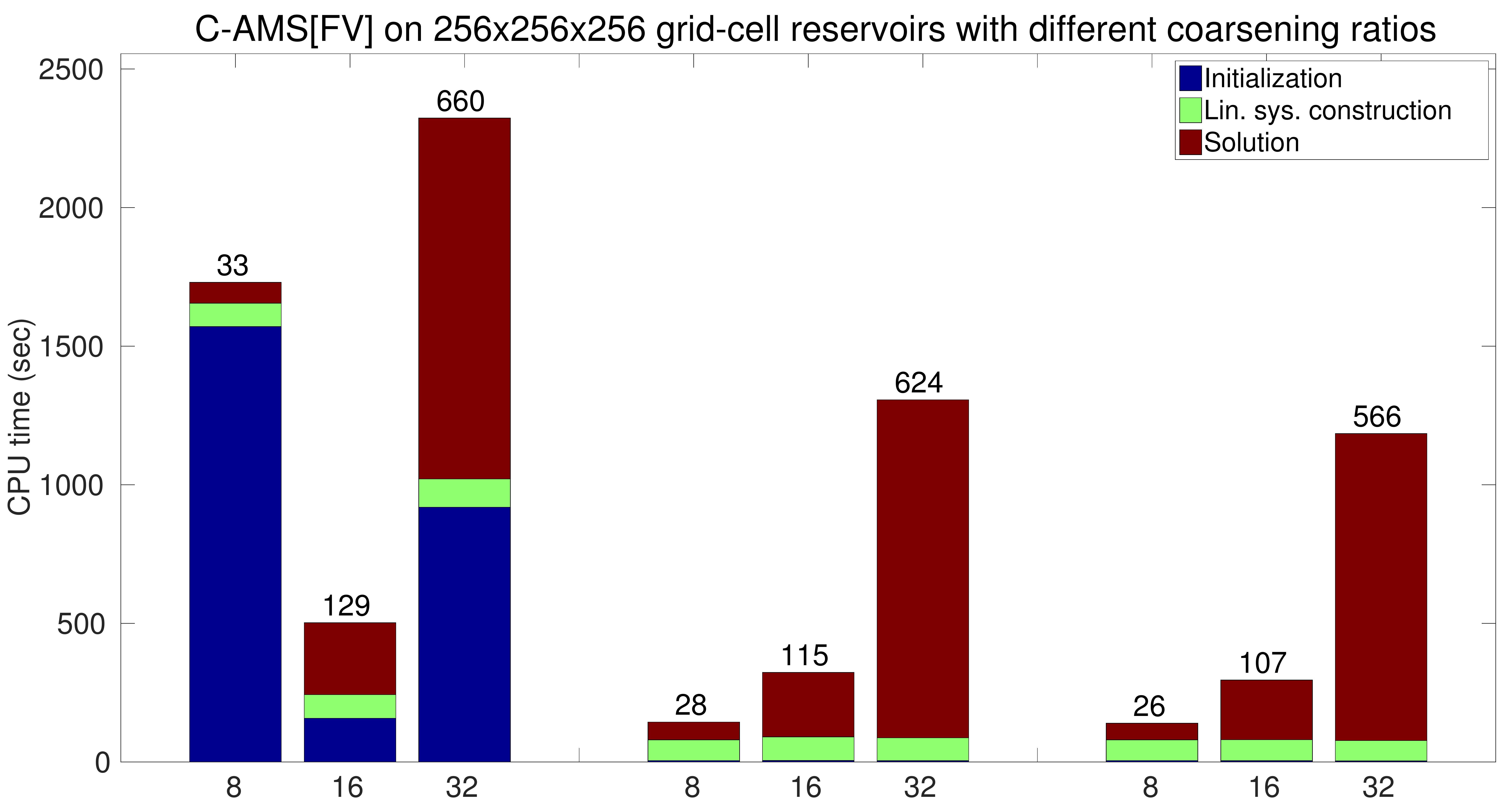}}
\begin{picture}(360,0)
\put(60,5){$t^*:$}
\put(90,5){0.0-0.4}
\put(190,5){0.4-1.0}
\put(295,5){1.0-2.0}
\end{picture}
\caption{Patchy fields: Averaged CPU time (over 20 realizations) comparison of C-AMS[FV] for different coarsening ratios for \change{the permeability Set 3 from Table \ref{tab:permeability set}}. A coarsening ratio of $16^3$ offers the best balance between initialization (basis function computation) and solution time, while $8^3$ results in a more expensive initialization but faster convergence in subsequent time-steps. \change{A similar behaviour was observed with the FE restriction operator.}}\label{coarse_256}
\end{figure}
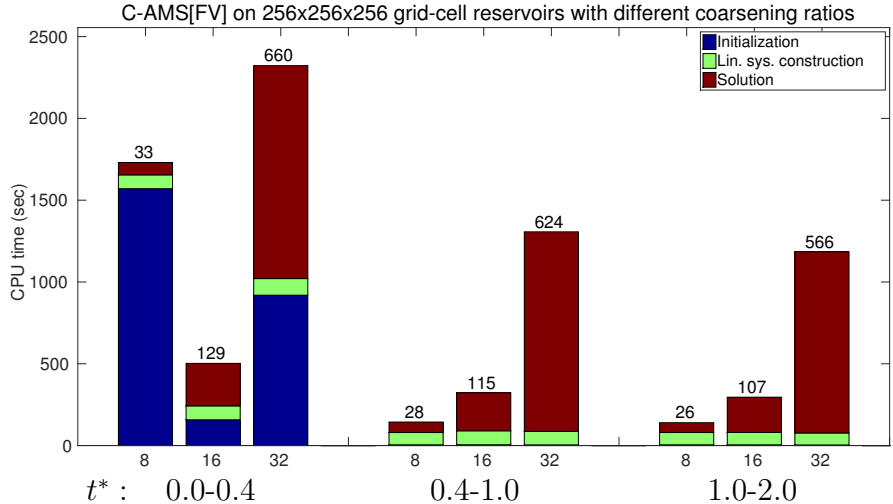



\subsection{C-AMS benchmark versus SAMG}

On the basis of the previously presented studies, the optimal C-AMS strategy includes a global multiscale stage using incompressible basis functions (Eq.~\eqref{B4}), accompanied by 5 iterations of ILU for post-smoothing. \change{In this subsection, C-AMS is compared against SAMG for three sets of different test cases: (1) the heterogeneous domains of different sizes from Table \ref{tab:permeability set}; (2): permeability Set 1 from Table \ref{tab:permeability set} with stretched grids and line-source terms; and, (3): permeability Set 1 from Table \ref{tab:permeability set} with different $ln(k)$ variances (i.e., permeability contrasts).}

\change{In all the presentd experiments, SAMG is called to perform a single V-cycle, repeatedly in a Richardson loop. Its adaptivity is controlled manually, i.e., at the beginning of each Newton-Raphson outer iteration, SAMG is allowed to update its Galerkin operators. On the other hand, during linear iterations, SAMG is instructed to reuse its previous grids and operators. For the test cases considered here, this approach was found more efficient (by a factor in excess of 2) than the \textit{automatic solver control} described in \cite{SAMG}, In all other aspects, SAMG has been used as a black-box commercial solver.} 

\newpage
\change{\subsubsection{Test case 1: heterogeneous domains of different sizes from Table \ref{tab:permeability set}}}
In this subsection, C-AMS is compared against the SAMG algebraic multigrid solver for both patchy and layered permeability fields \change{of Table \ref{tab:permeability set}} over 3 consecutive time steps. The time-lapse pressure solution for one patchy and one layered sample are shown in Fig.~\ref{fig:sol64}, illustrating the propagation of the signal from the western face through the entire domain.

\begin{figure} [htb!]
\centering
\includegraphics[height=0.370\textwidth]{5.pdf}\hspace{2cm}
\includegraphics[height=0.370\textwidth]{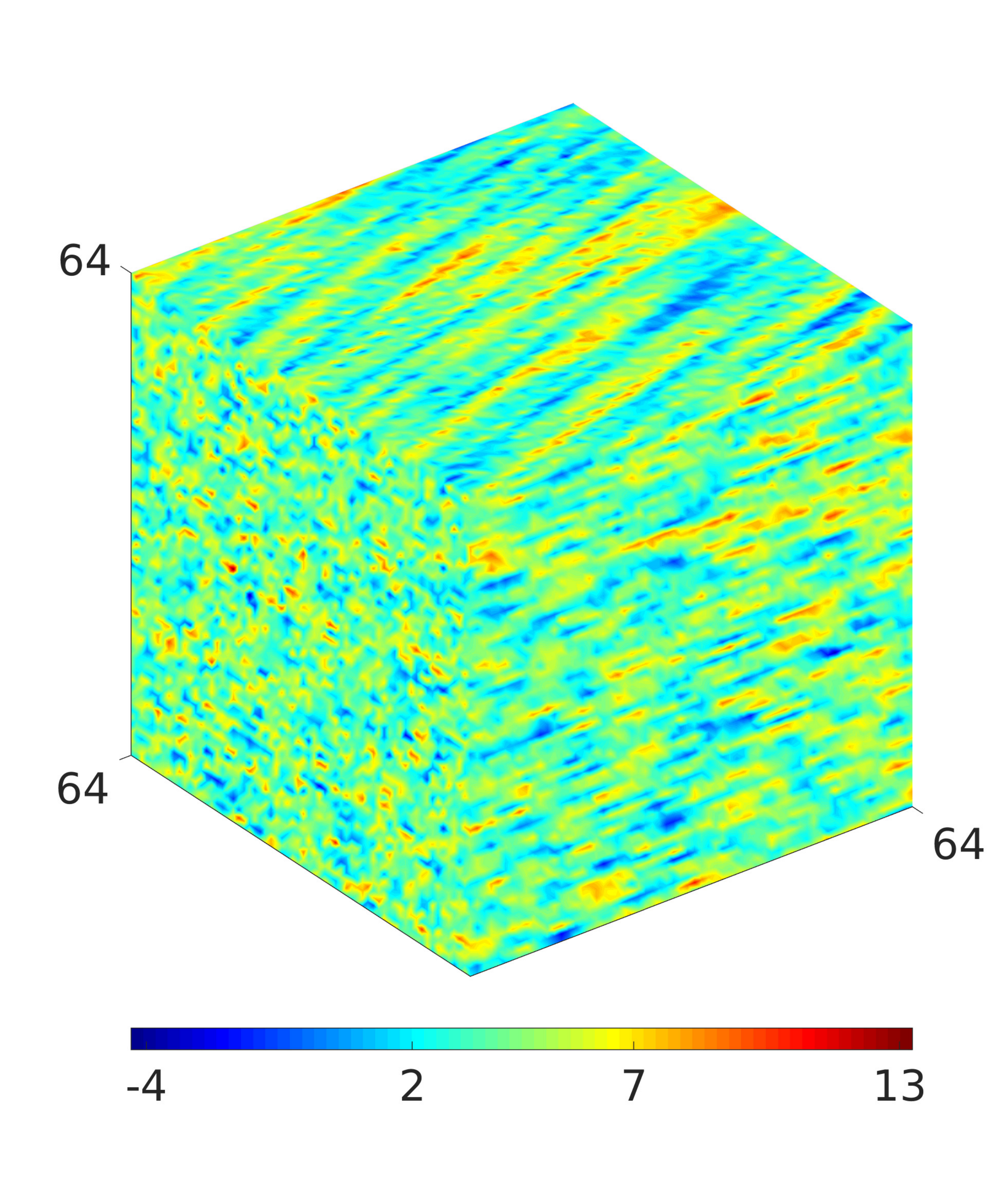}\\
\includegraphics[height=0.370\textwidth]{6.pdf}\hspace{2cm}
\includegraphics[height=0.370\textwidth]{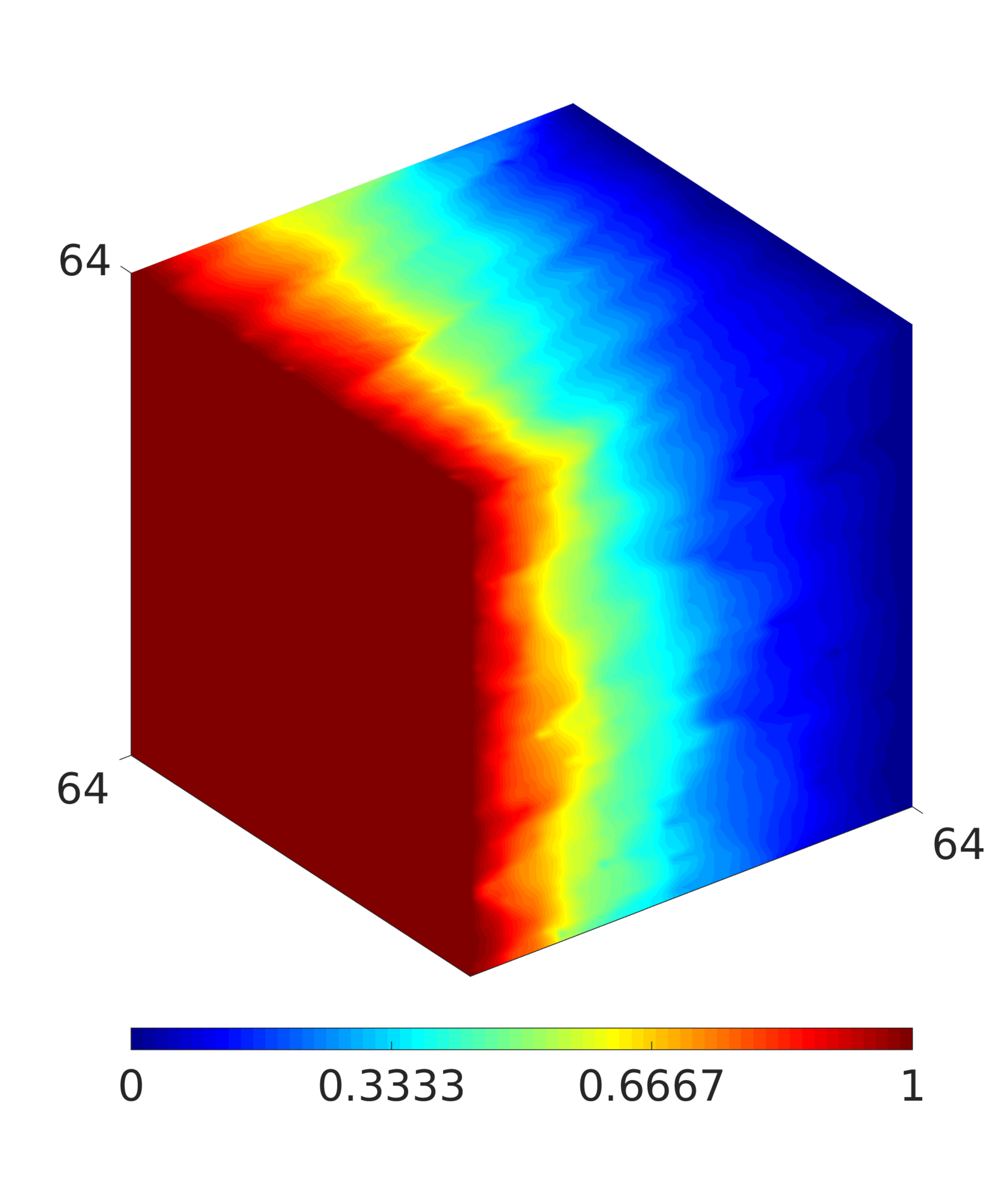}\\
\includegraphics[height=0.370\textwidth]{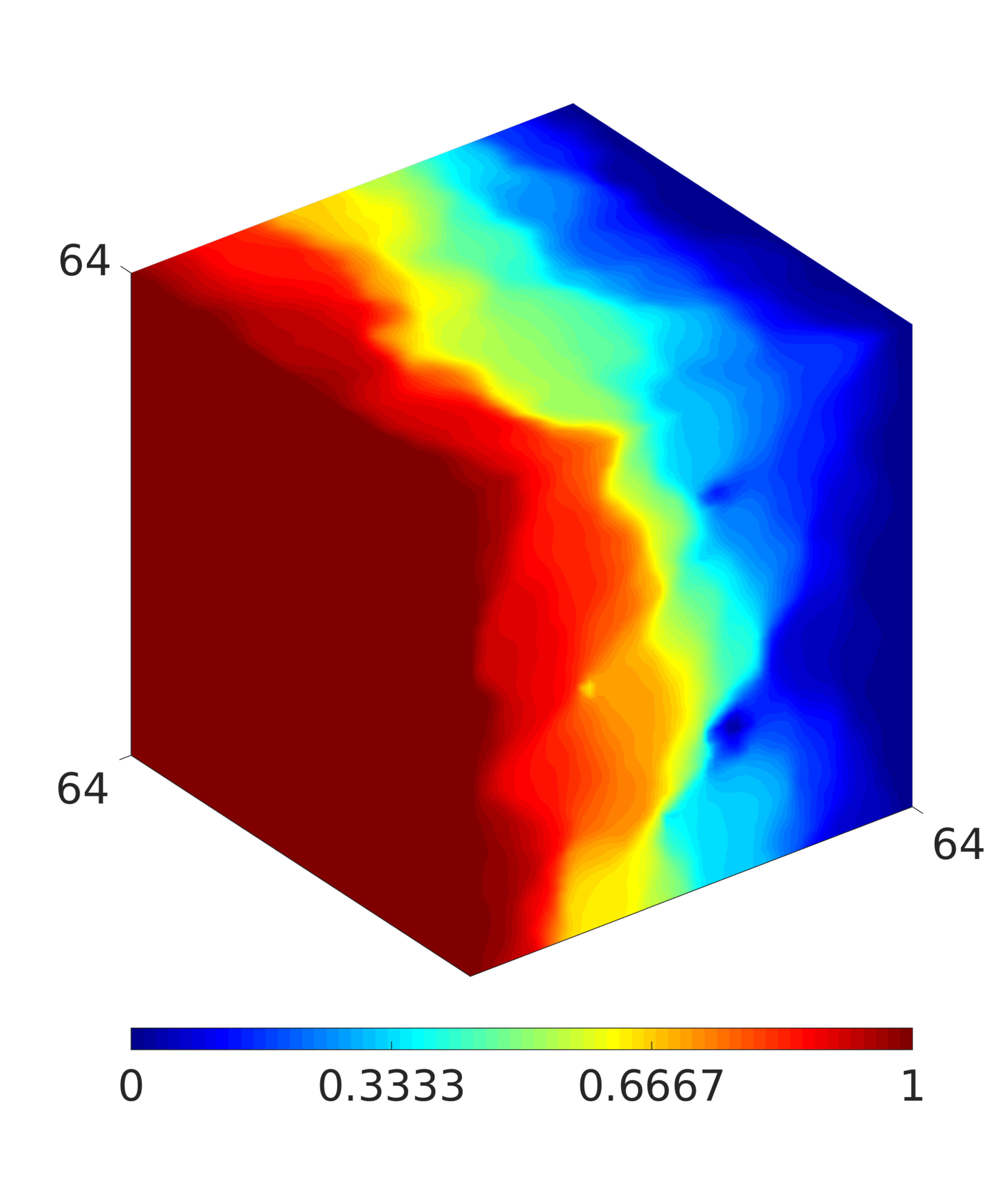}\hspace{2cm}
\includegraphics[height=0.370\textwidth]{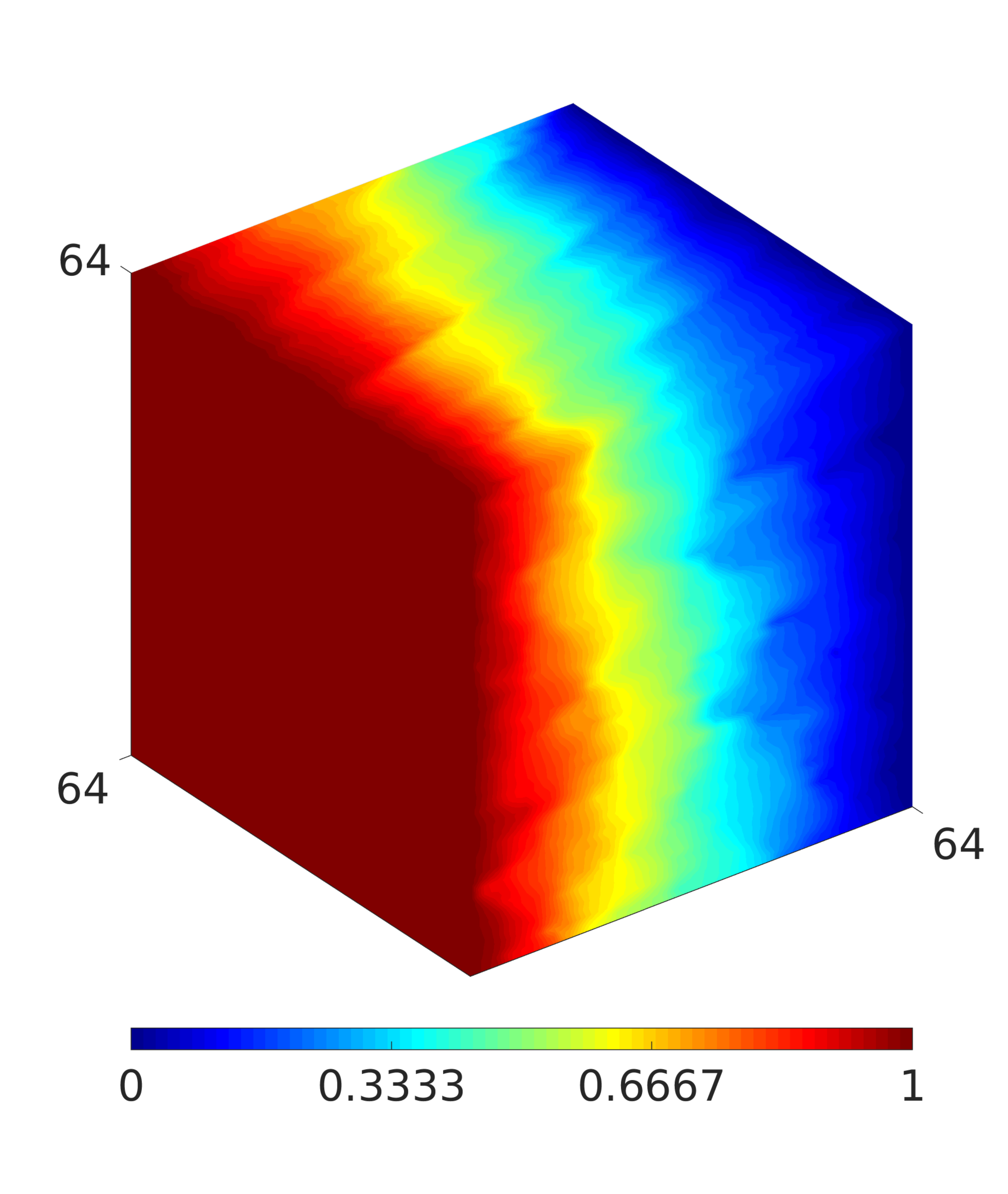} \\
\includegraphics[height=0.370\textwidth]{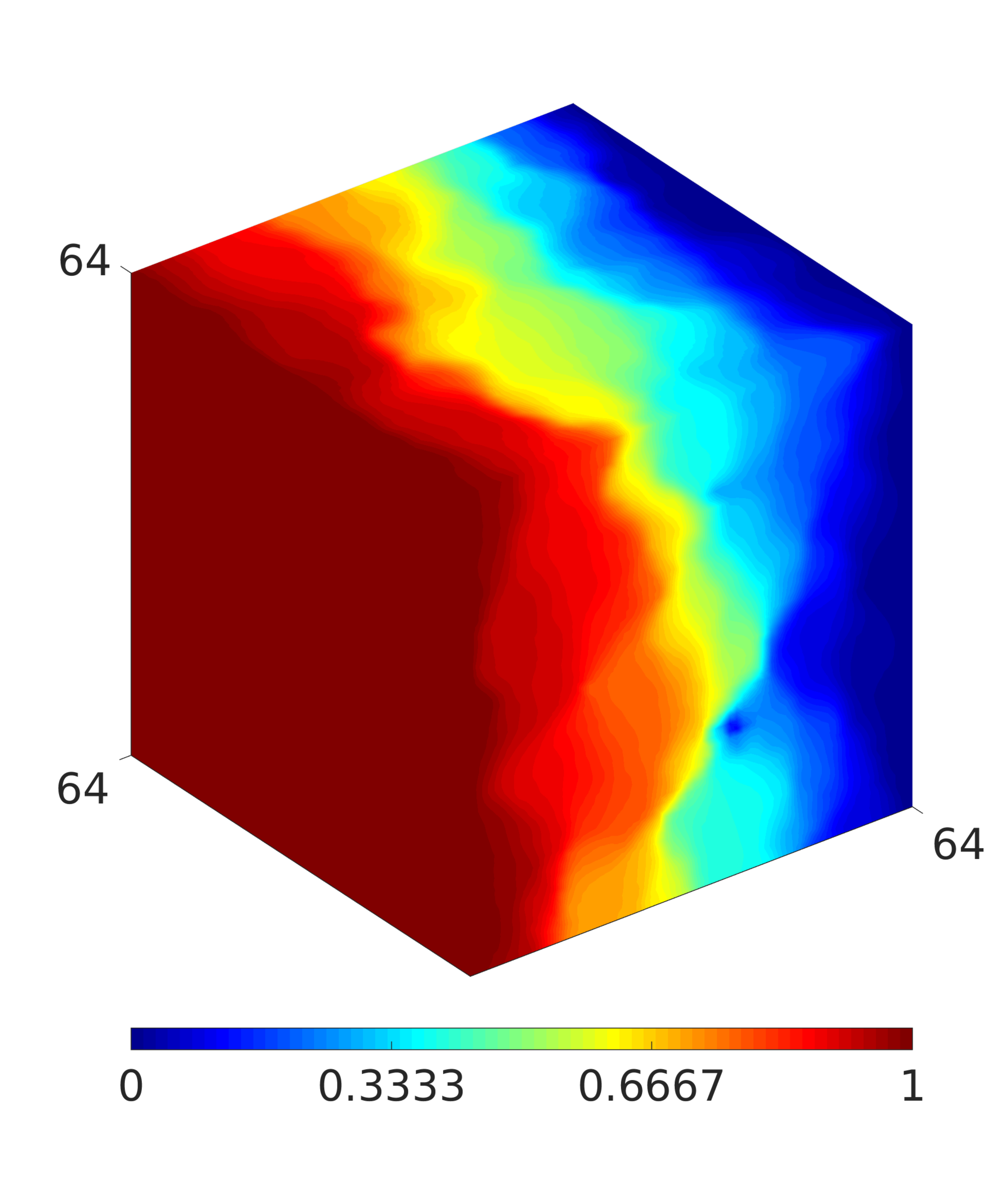}\hspace{2cm}
\includegraphics[height=0.370\textwidth]{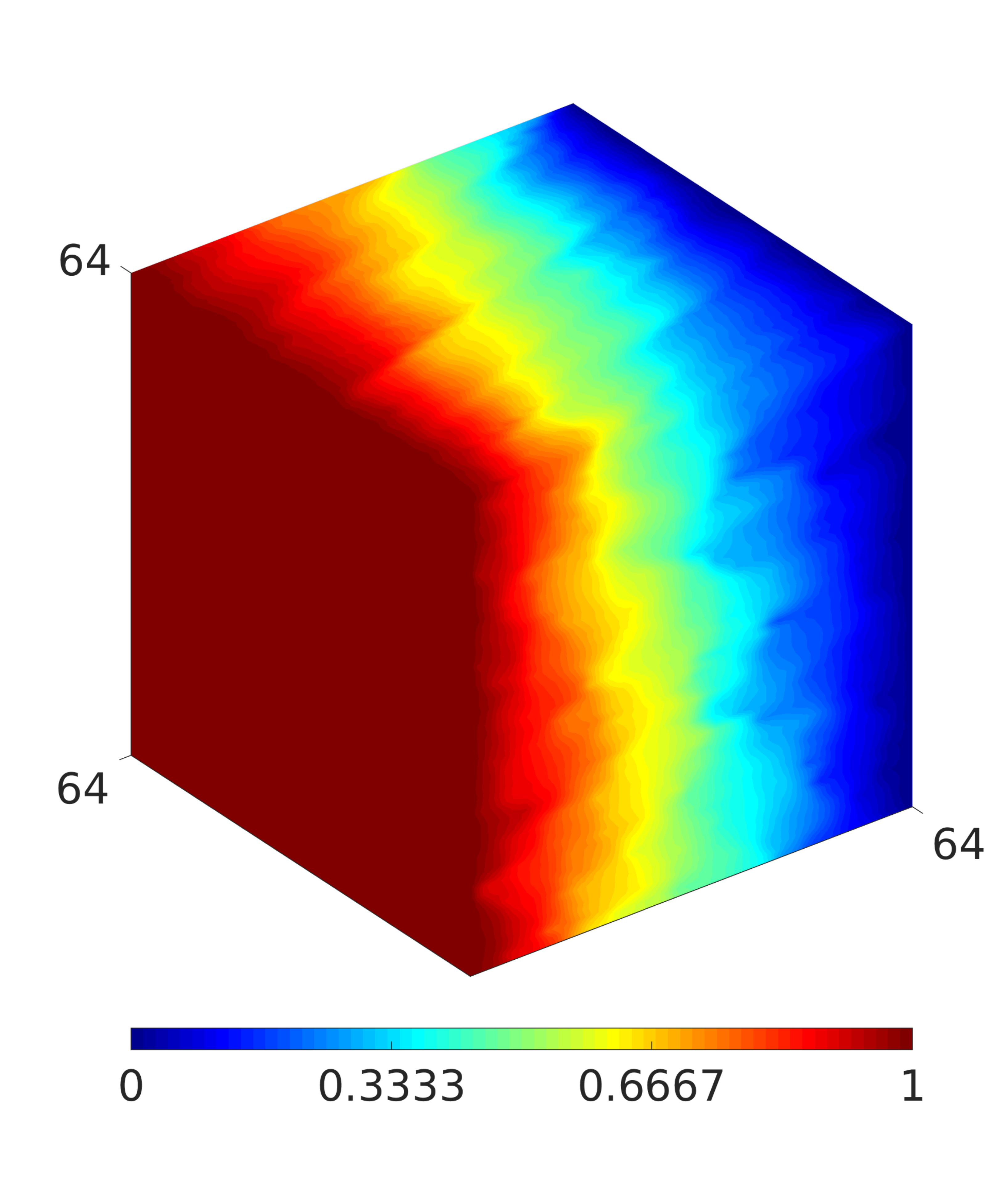}
\caption{Pressure solution on \change{one of the} realizations of \change{permeability Sets 1 (left) and 4 (right) from Table \ref{tab:permeability set} at $t^*=0.4, 1.0$}, and \change{$2.0$} from top to bottom, respectively.}
\label{fig:sol64}
\end{figure}

Figs.~\ref{fig:cpu64} and \ref{fig:cpu128} show the number of iterations and CPU time at 3 consecutive non-dimensional times for different problem \change{Sets 1, 2, 4, and 5 from Table \ref{tab:permeability set}}. \change{Note that C-AMS with FV-based restriction operator did not converge in some of the test cases, while the FE-based variant achieved 100\% success rate due to its SPD property. Therefore, an ideal solution strategy would use MSFE to converge to the desired level of accuracy and then employ a single MSFV sweep, in order to ensure mass conservation  \cite{hadi-erest-spej}.}

\begin{figure} [htb!]
\centering
\centerline{%
\includegraphics[height=0.48\textwidth]{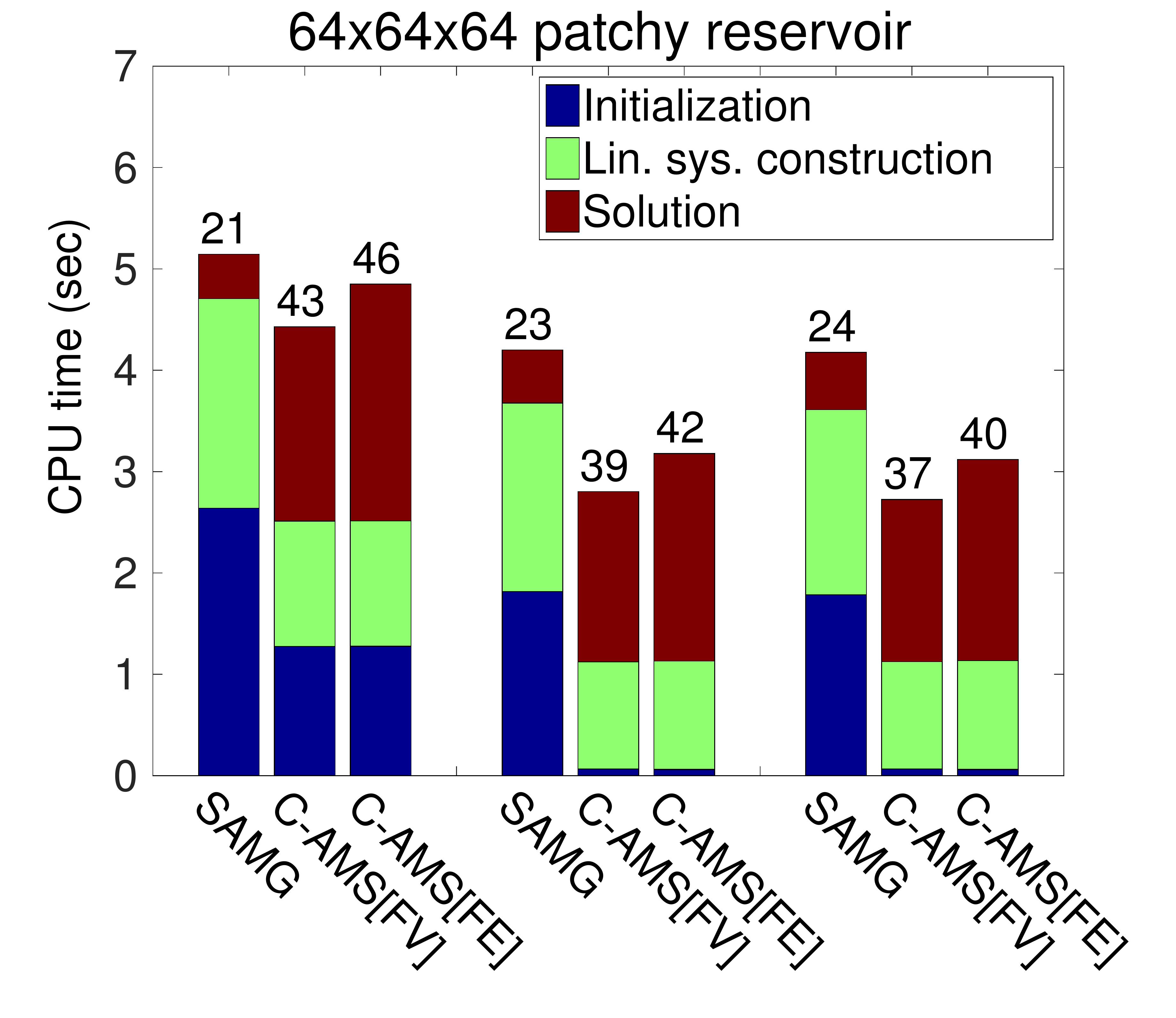}
\includegraphics[height=0.48\textwidth]{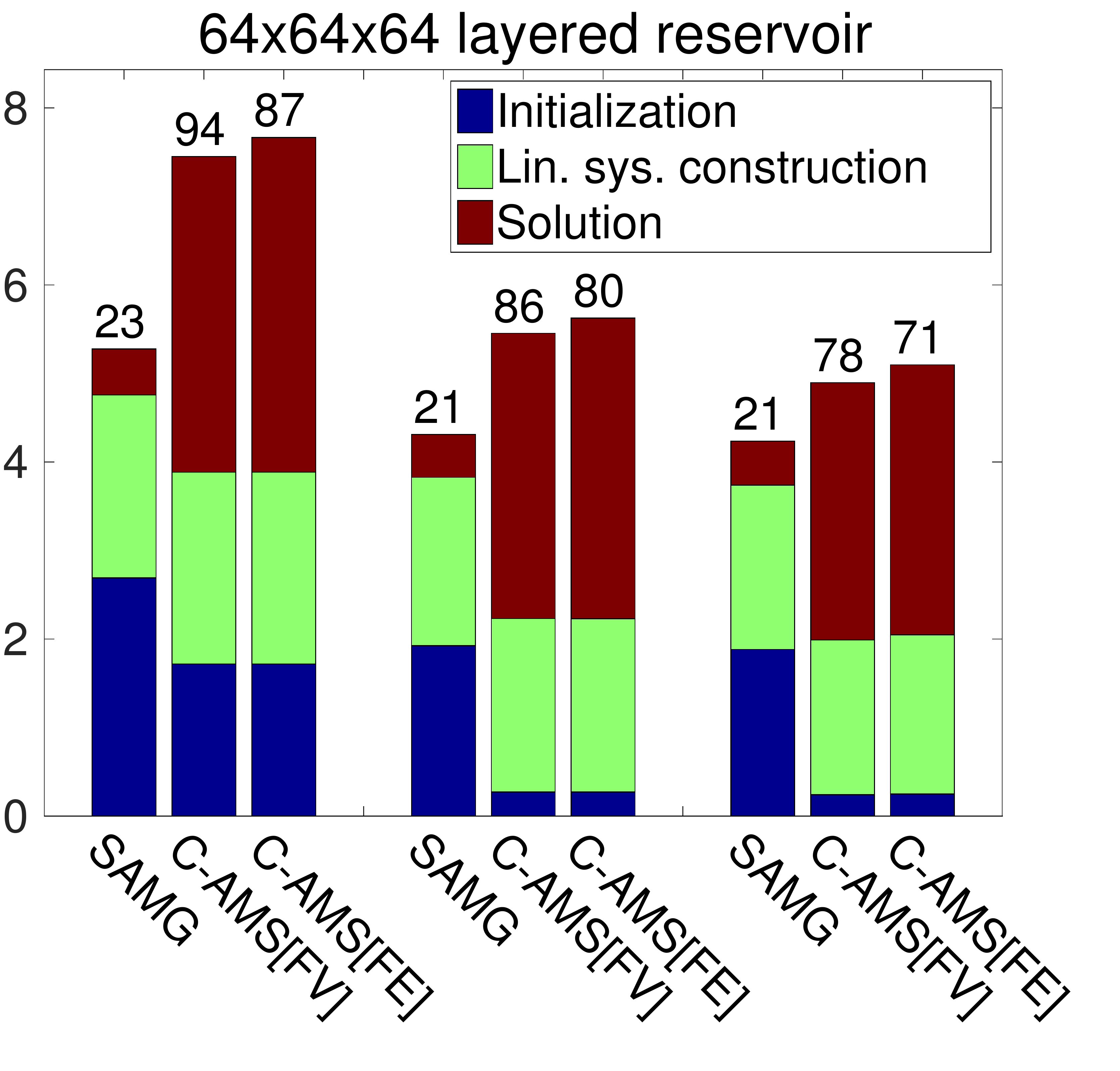}
}
\begin{picture}(360,0)
\put(-5,10){$t^*:$}
\put(20,10){0.0-0.4}
\put(70,10){0.4-1.0}
\put(120,10){1.0-2.0}
\put(210,10){0.0-0.4}
\put(275,10){0.4-1.0}
\put(325,10){1.0-2.0}
\end{picture}
\caption{Averaged CPU time (over 20 realizations) comparison between the C-AMS and SAMG solvers on \change{permeability Sets 1}  (a) and \change{4} (b) \change{from Table \ref{tab:permeability set}} over 3 successive time-steps. The coarsening ratio of C-AMS is $8^3$. Moreover, C-AMS employs 5 ILU smoothing steps per iteration. The number of iterations is given on top of each bar.}\label{fig:cpu64}
\end{figure}

\begin{figure} [htb!]
\centering
\centerline{%
\includegraphics[height=0.48\textwidth]{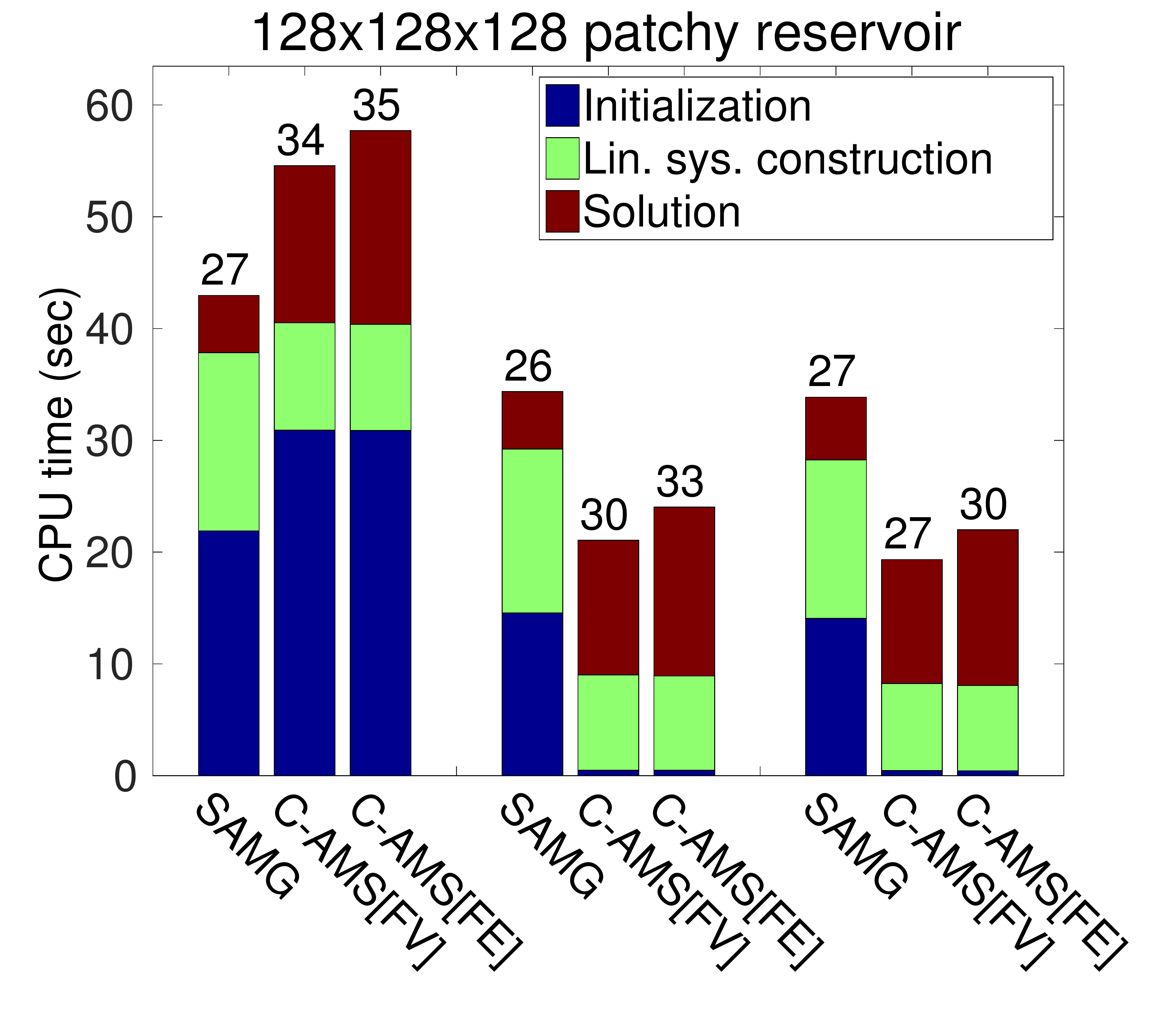}
\includegraphics[height=0.48\textwidth]{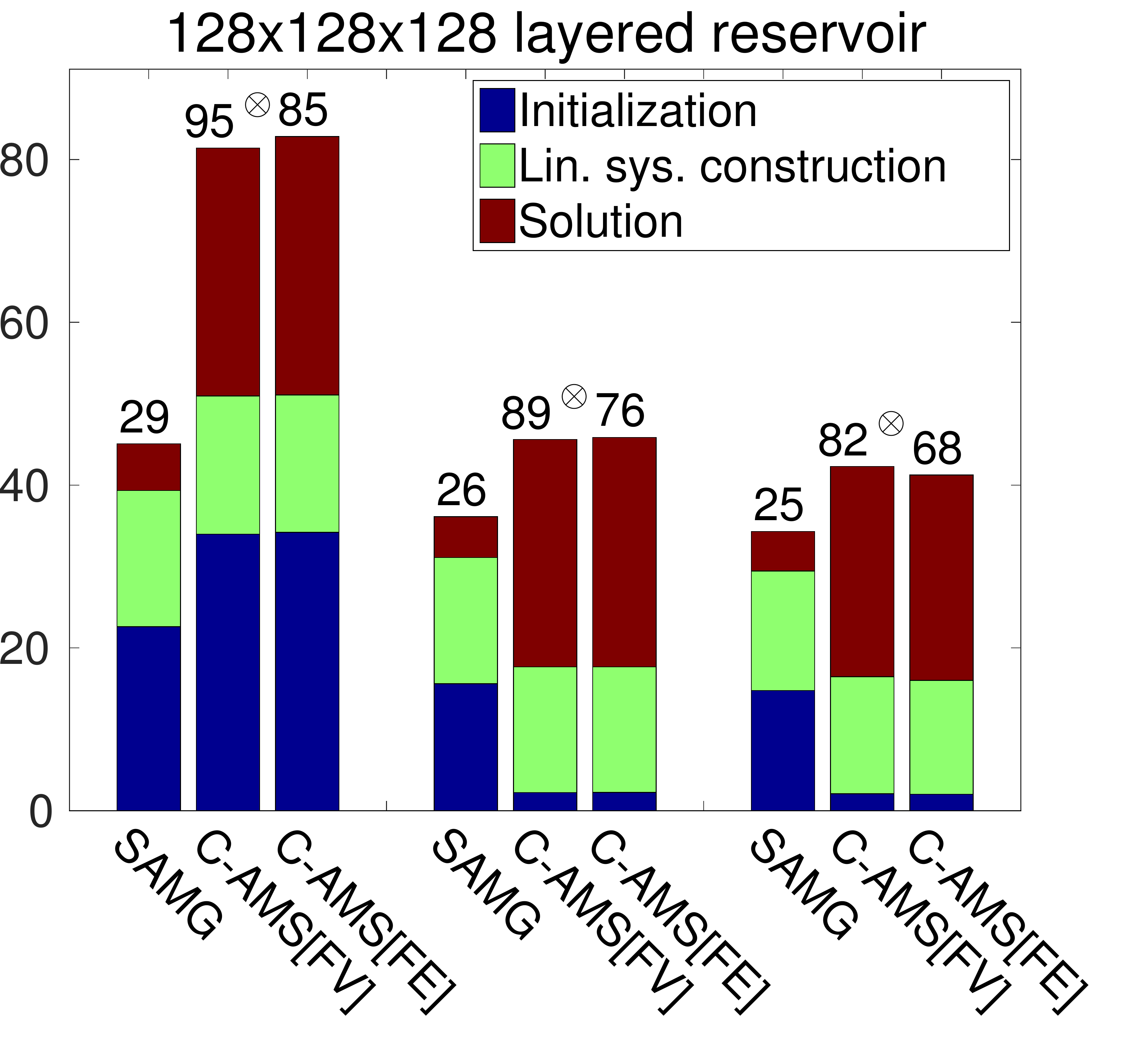}
}
\begin{picture}(360,0)
\put(-5,10){$t^*:$}
\put(20,10){0.0-0.4}
\put(70,10){0.4-1.0}
\put(120,10){1.0-2.0}
\put(210,10){0.0-0.4}
\put(275,10){0.4-1.0}
\put(325,10){1.0-2.0}
\end{picture}
\caption{Averaged CPU time (over 20 realizations) comparison between the C-AMS and SAMG solvers on \change{permeability Sets 2} (a) and \change{5} (b)  \change{from Table \ref{tab:permeability set}} over 3 successive time steps. C-AMS employs the coarsening ratio of $8^3$, along with 5 ILU smoothing steps per iteration. The number of iterations is given on top of each bar. \change{The $\otimes$ symbol signifies convergence success rate of 72\% when FV-based restriction operator was employed.}}\label{fig:cpu128}
\end{figure}

In addition, Figs.~\ref{fig:cpu256} illustrate CPU time (vertical axis) and the total number of iterations (on top of each column), \change{for permeability Sets 3 and 6 from Table \ref{tab:permeability set}, with $8^3$ and $16^3$ coarsening ratios}.

\begin{figure} [htb!]
\centering
\includegraphics[height=0.45\textwidth]{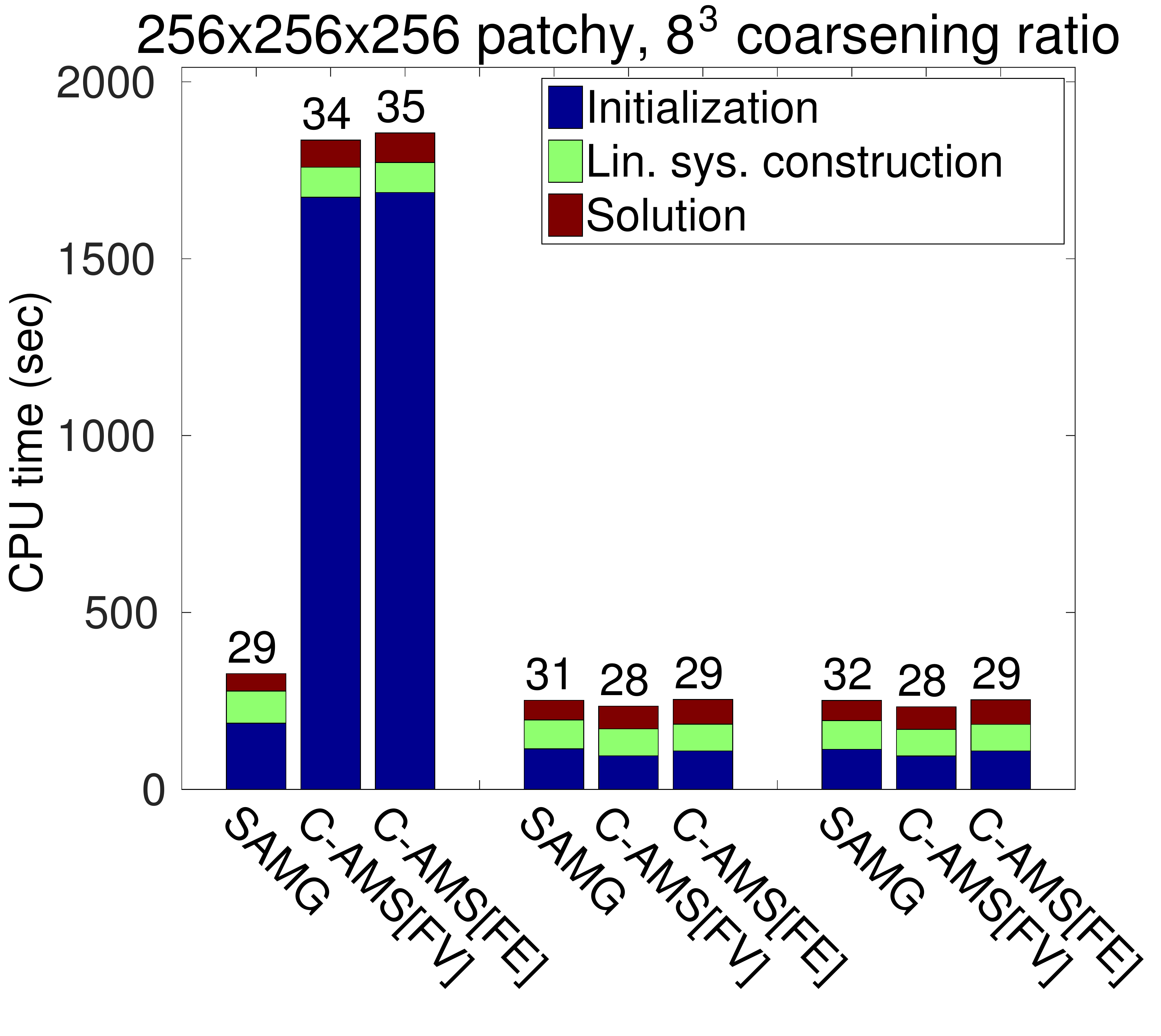}
\includegraphics[height=0.45\textwidth]{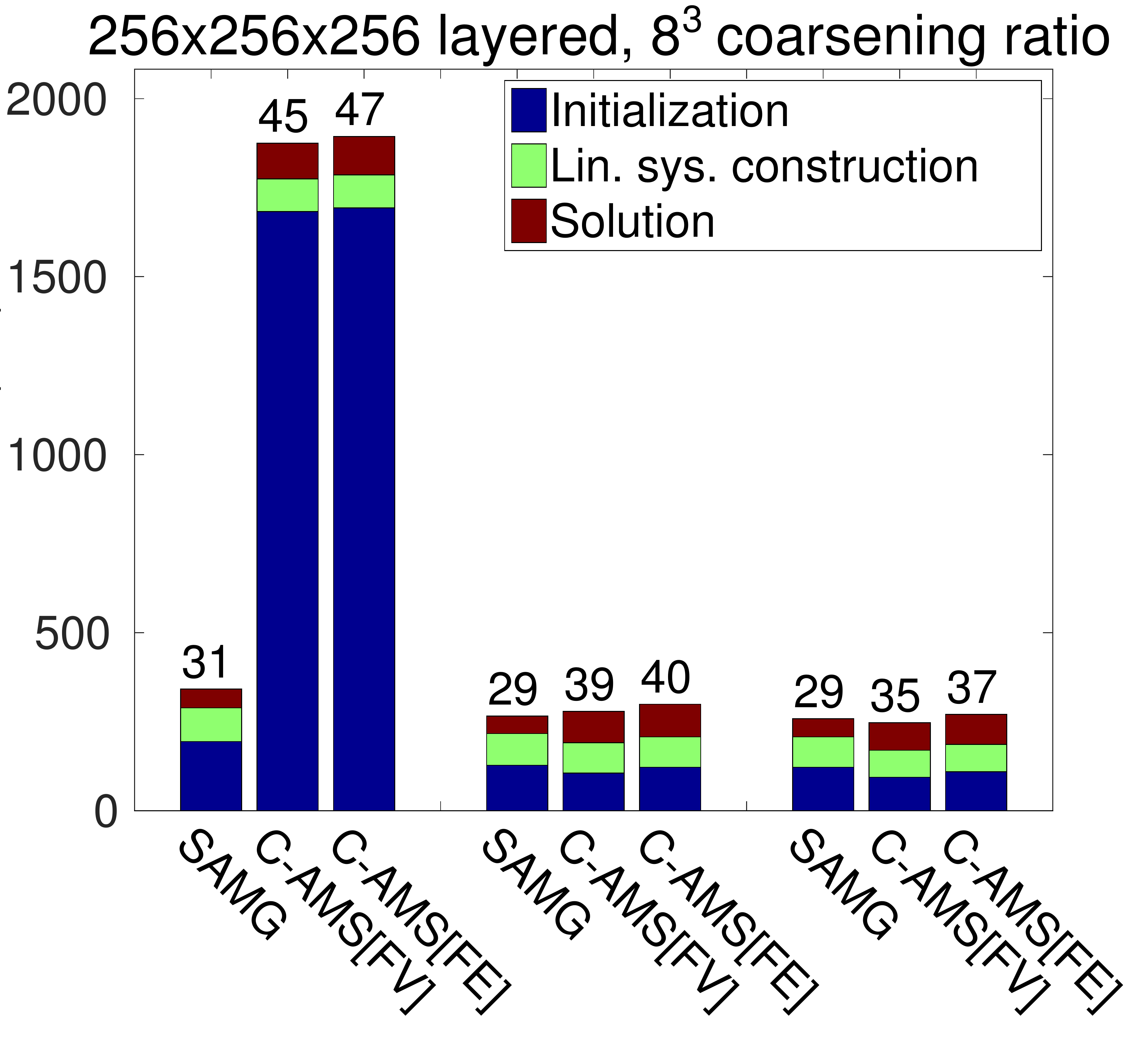}\\
\vspace{3mm}
\includegraphics[height=0.45\textwidth]{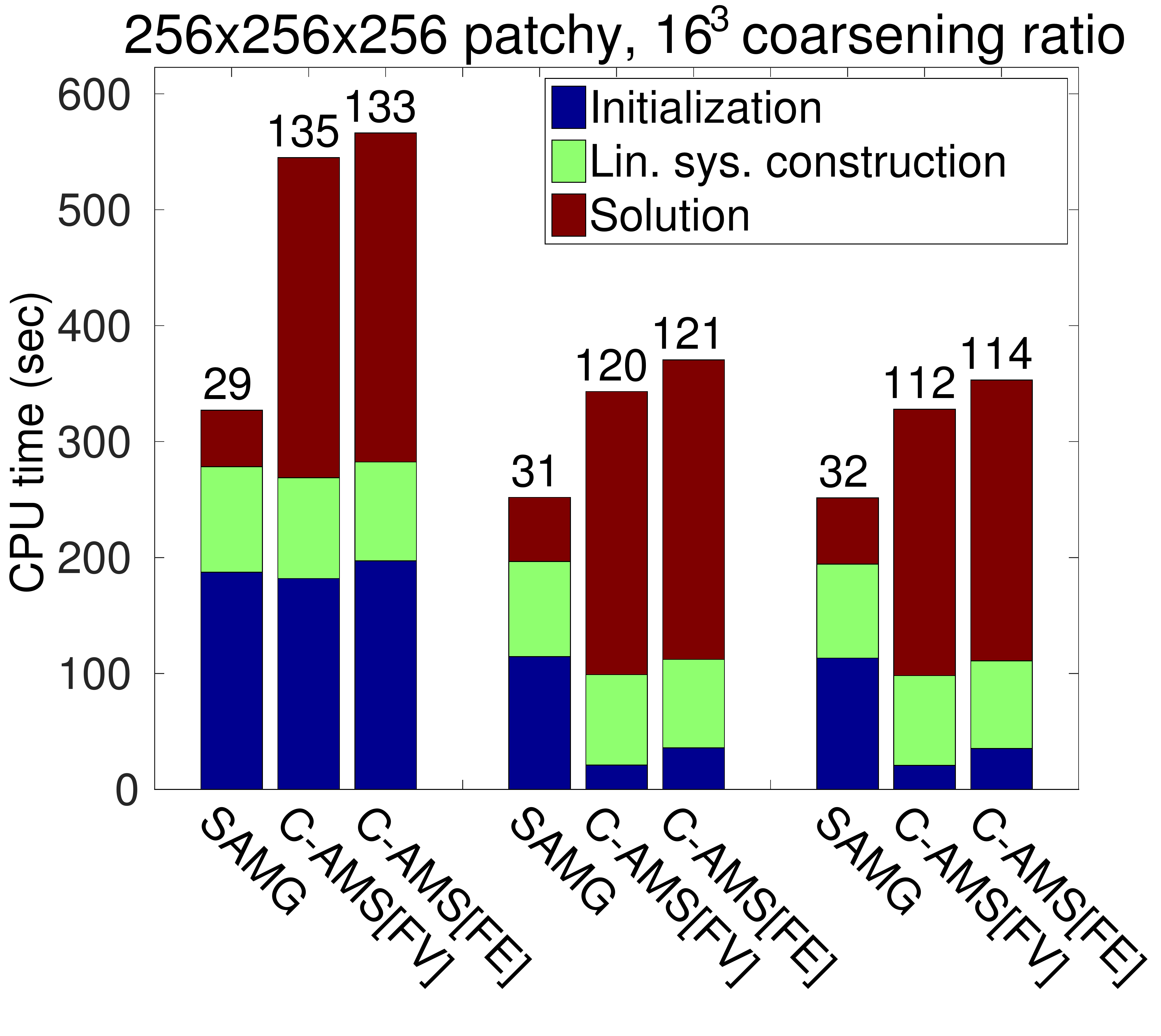}
\includegraphics[height=0.45\textwidth]{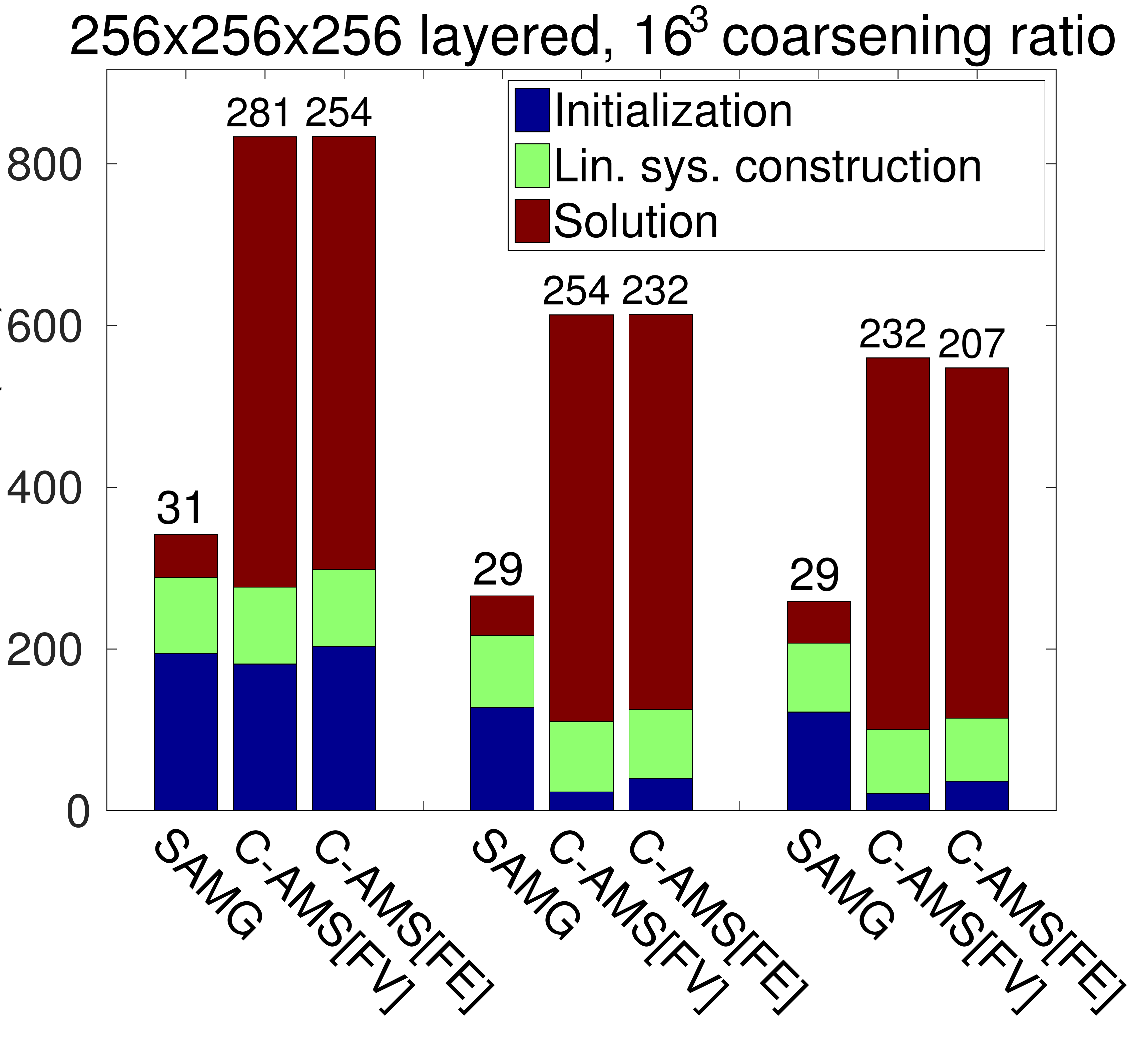}
\begin{picture}(360,0)
\put(-5,10){$t^*:$}
\put(20,10){0.0-0.4}
\put(78,10){0.4-1.0}
\put(130,10){1.0-2.0}
\put(213,10){0.0-0.4}
\put(275,10){0.4-1.0}
\put(325,10){1.0-2.0}
\end{picture}
\caption{Averaged CPU time (over 20 realizations) comparison between the C-AMS and SAMG solvers on \change{permeability Sets 3} (left column) and \change{6} (right column) \change{from Table \ref{tab:permeability set}} over 3 successive time-steps. Different coarsening \change{ratios} of $8^3$ (top row) and $16^3$ (bottom row) are considered for C-AMS. Moreover, C-AMS employs 5 ILU smoothing steps per iteration. The number of iterations is given on top of each bar.}\label{fig:cpu256}
\end{figure}

Note that, except for the first time-step, when all the basis functions are fully computed, C-AMS has a slight edge over SAMG, mainly due to its adaptivity and relatively inexpensive iterations. The initialization cost of C-AMS is particularly high in the $256^3$ case, due to the large number of linear systems (solved with a direct solver) needed for the basis functions. \change{It is clear from} Fig.~\ref{fig:cpu256} that with larger primal-coarse blocks C-AMS requires less setup time, but more iterations to converge. Note that all performance studies presented in this paper are for single-process computations. 

Since reservoir simulators are typically run for many time-steps, the high initialization time of C-AMS is outweighed by the efficiency gained in subsequent steps. Moreover, given the local support of the basis functions, this initialization can be greatly improved through parallel processing. Furthermore, only a few multiscale iteration may prove necessary to obtain an accurate approximation of the pressure solution in each time step for multi-phase flow problems.

\clearpage
\subsubsection{\change{Test Case 2: stretched grids with line-source terms}}
\change{To study the effect of anisotropic permeability fields along with radial injection flow pattern, the permeability Set 1 from Table \ref{tab:permeability set} is considered. The settings are all the same as previous test cases, except the following items. Dirichlet boundary conditions are set at the centers of two vertical sets of fine-scale grid cells: one from (1,1,1) to (1,1,64) and the other from (64,64,1) to (64,64,64) with the values of 1 and 0, respectively. In addition, grid aspect ratios of $\alpha = 1, 5, $ and $10$ are considered (Note that $\Delta x/\alpha = \Delta y = \Delta z$). The non-dimensional time is calculated using $\alpha L$ as the characteristic length. Figure \ref{stretched_p} illustrates the pressure solutions for one of the permeability realizations after the first time step $t^* = 0.4$.}

\begin{figure} [htb!]
\centering
\includegraphics[width=0.32\textwidth]{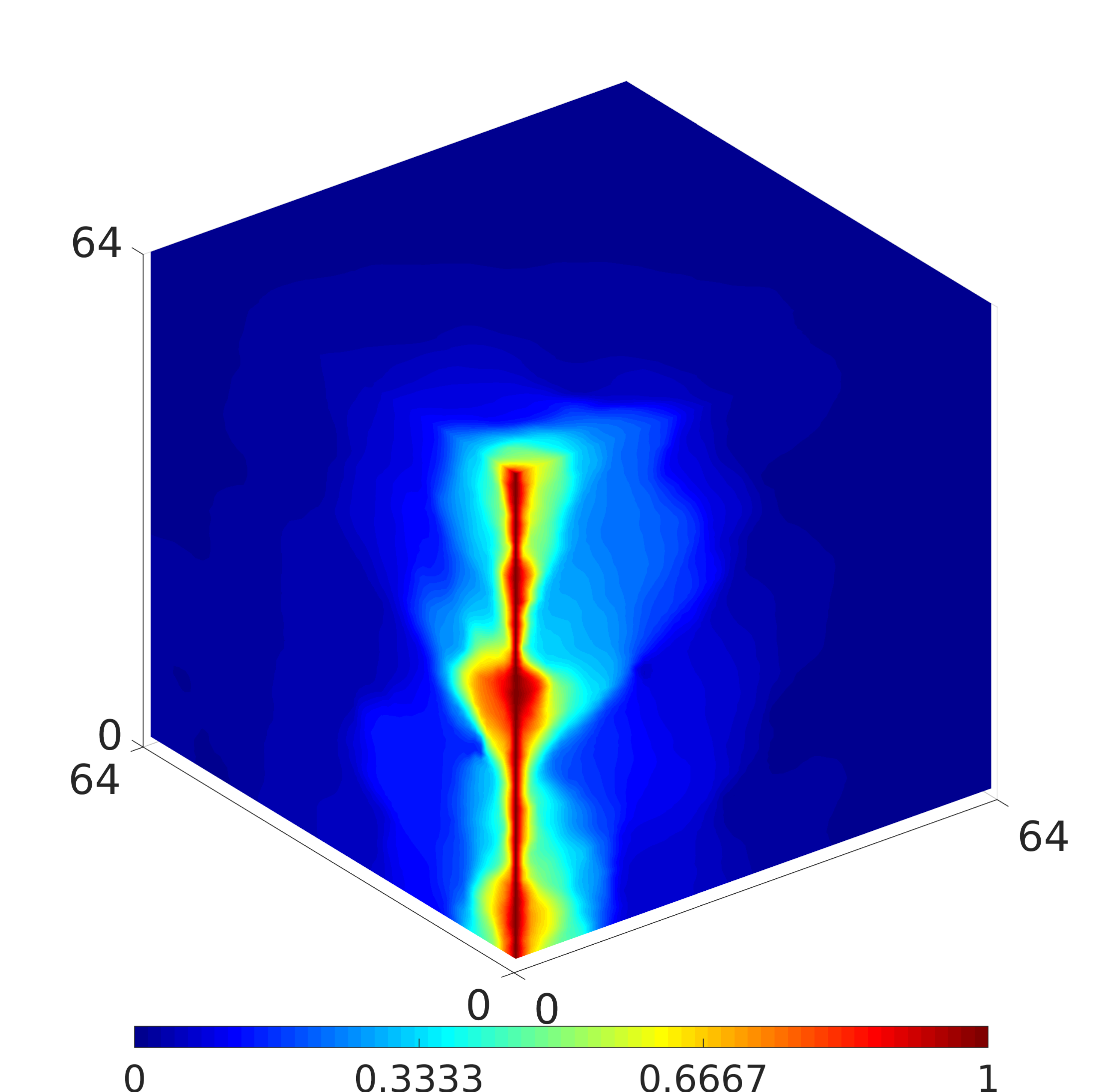}
\includegraphics[width=0.32\textwidth]{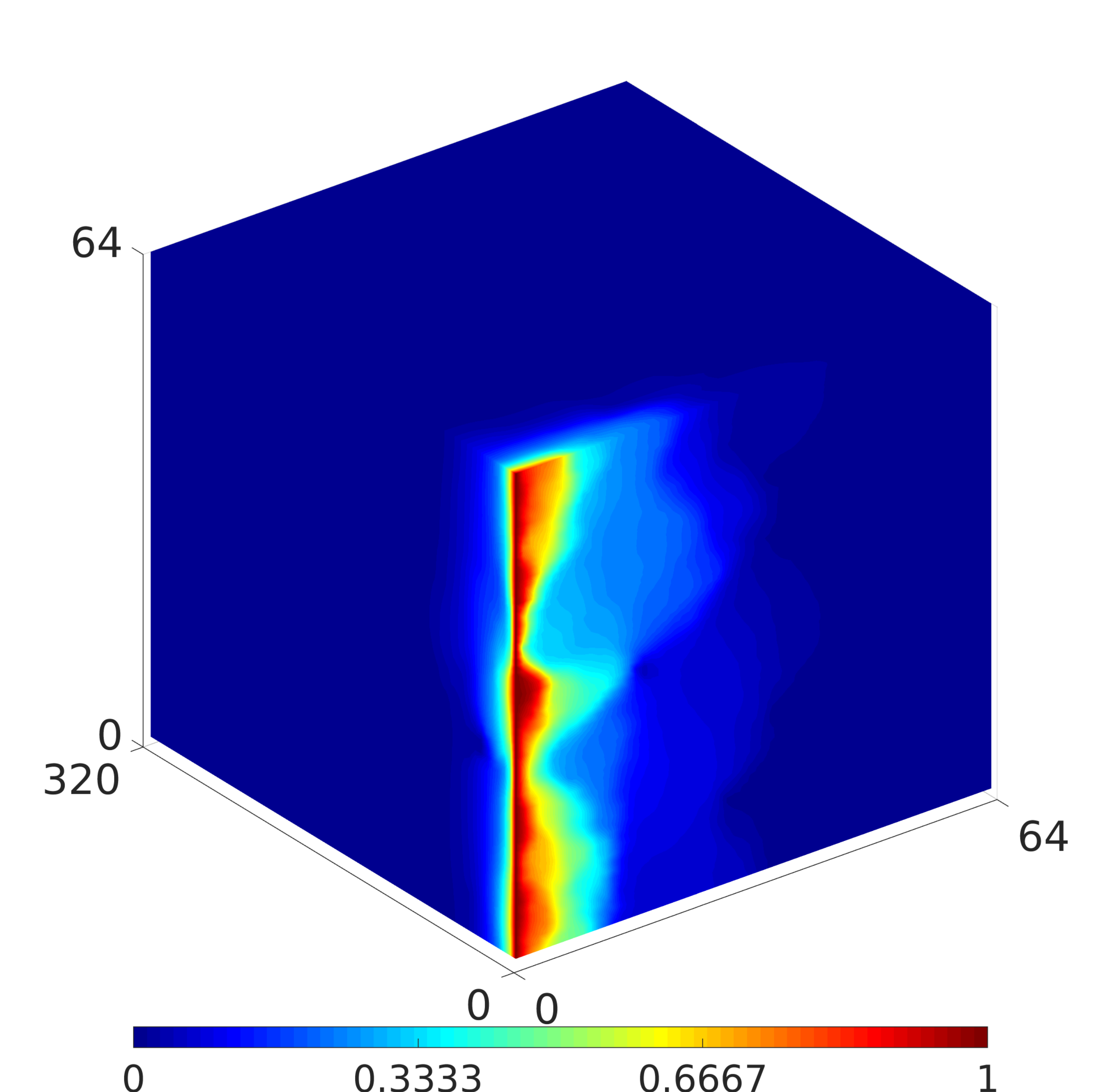}
\includegraphics[width=0.32\textwidth]{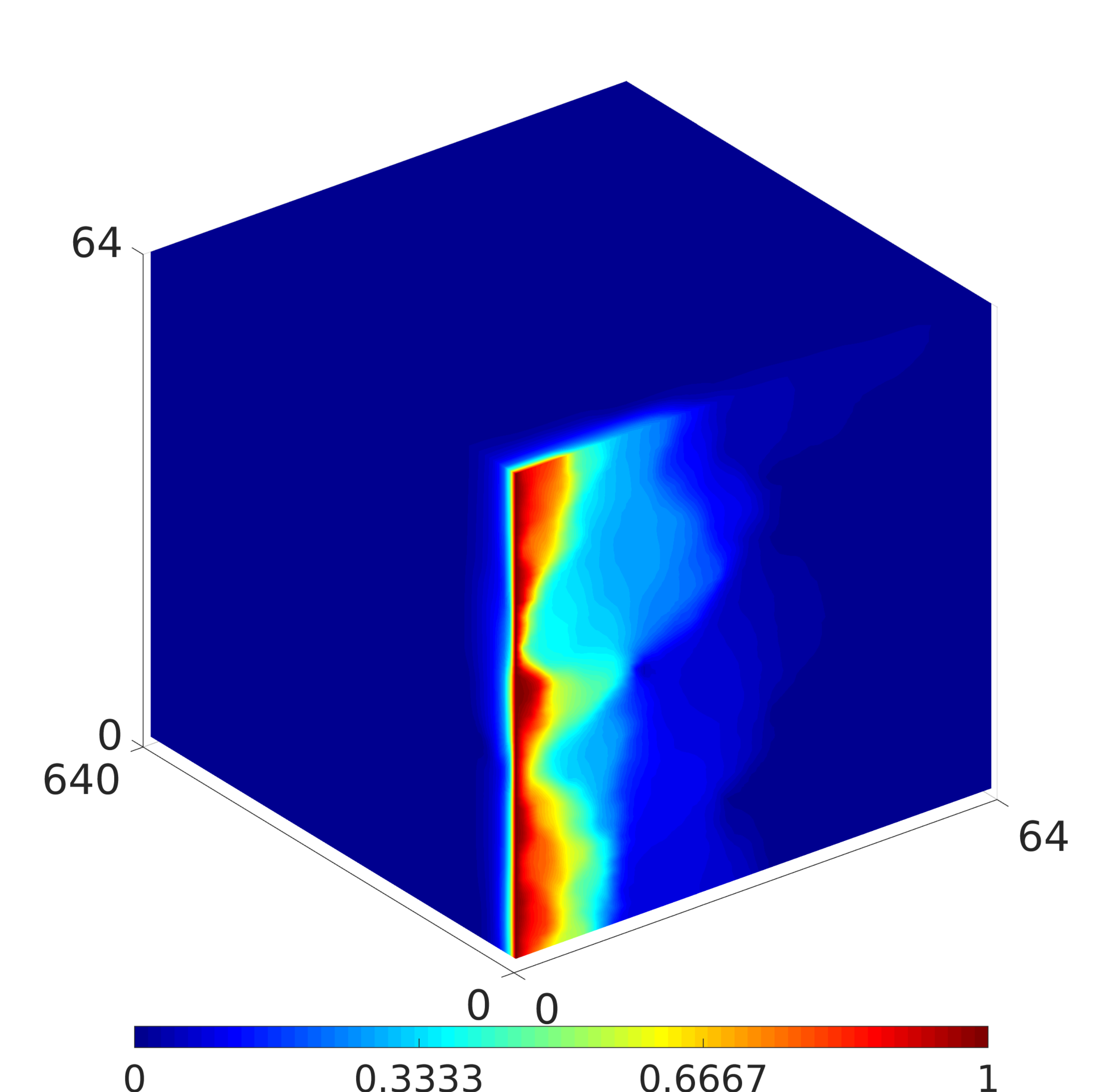}
\caption{\change{Converged pressure solution for one of the realizations of permeability Sets 1 with grid aspect ratio $\alpha = 1, 5,$ and $10$, respectively from left to right, after one time step $t^* = 0.4$. Dirichlet boundary conditions are set at the centers of two vertical sets of fine-scale grid cells: one from (1,1,1) to (1,1,64) and the other from (64,64,1) to (64,64,64) with the values of 1 and 0, respectively.}}\label{stretched_p}
\end{figure}

\change{The performance of C-AMS[FE] and SAMG are presented in Fig. \ref{stretched_cpu}. In contrast to C-AMS[FE], the C-AMS[FV] (not shown) did not lead to $100\%$ convergence success. However, for those C-AMS[FV] successful runs, similar CPU times as in C-AMS[FE] were observed.}

\begin{figure} [htb!]
\centering
\includegraphics[width=0.48\textwidth]{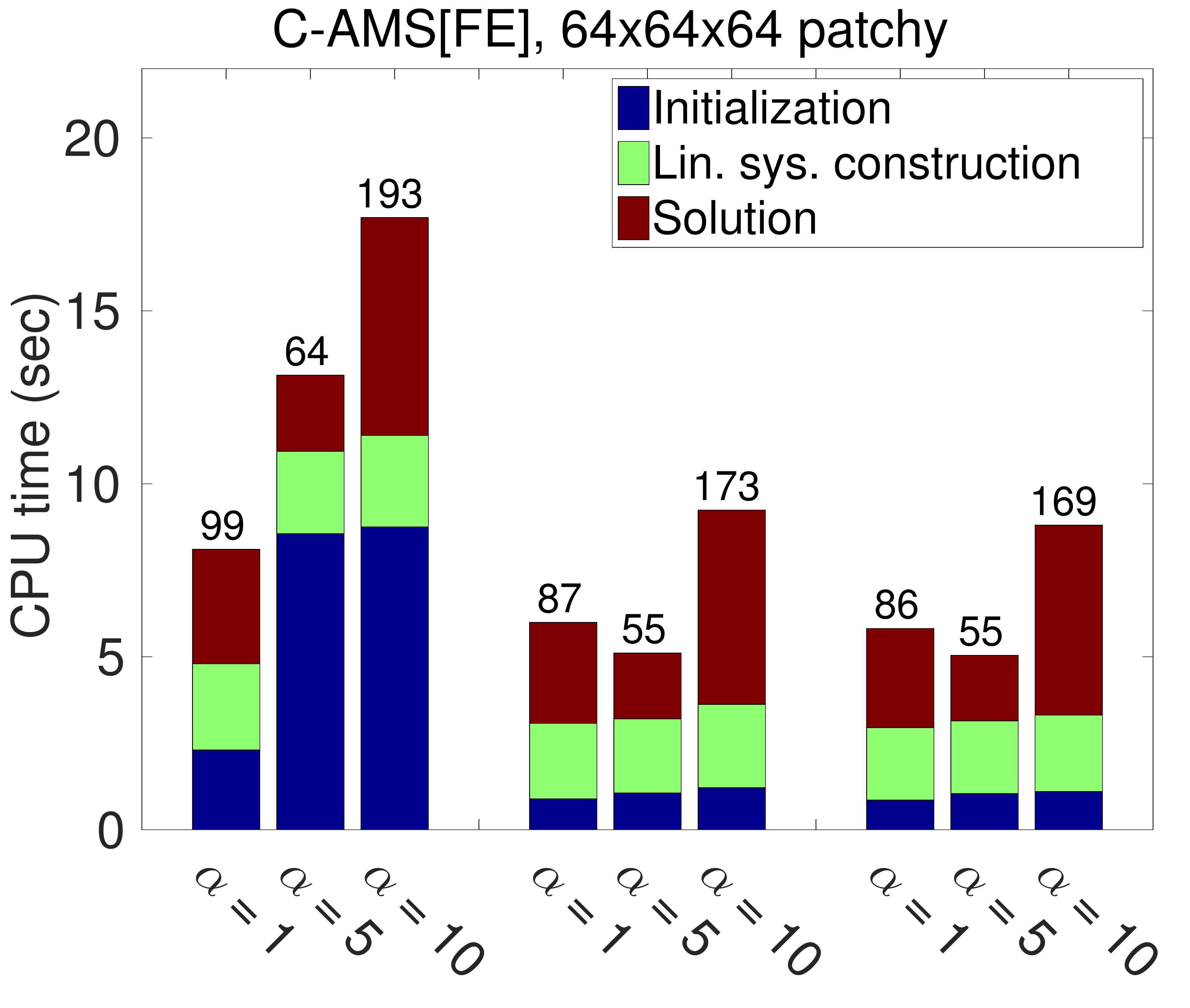}
\includegraphics[width=0.48\textwidth]{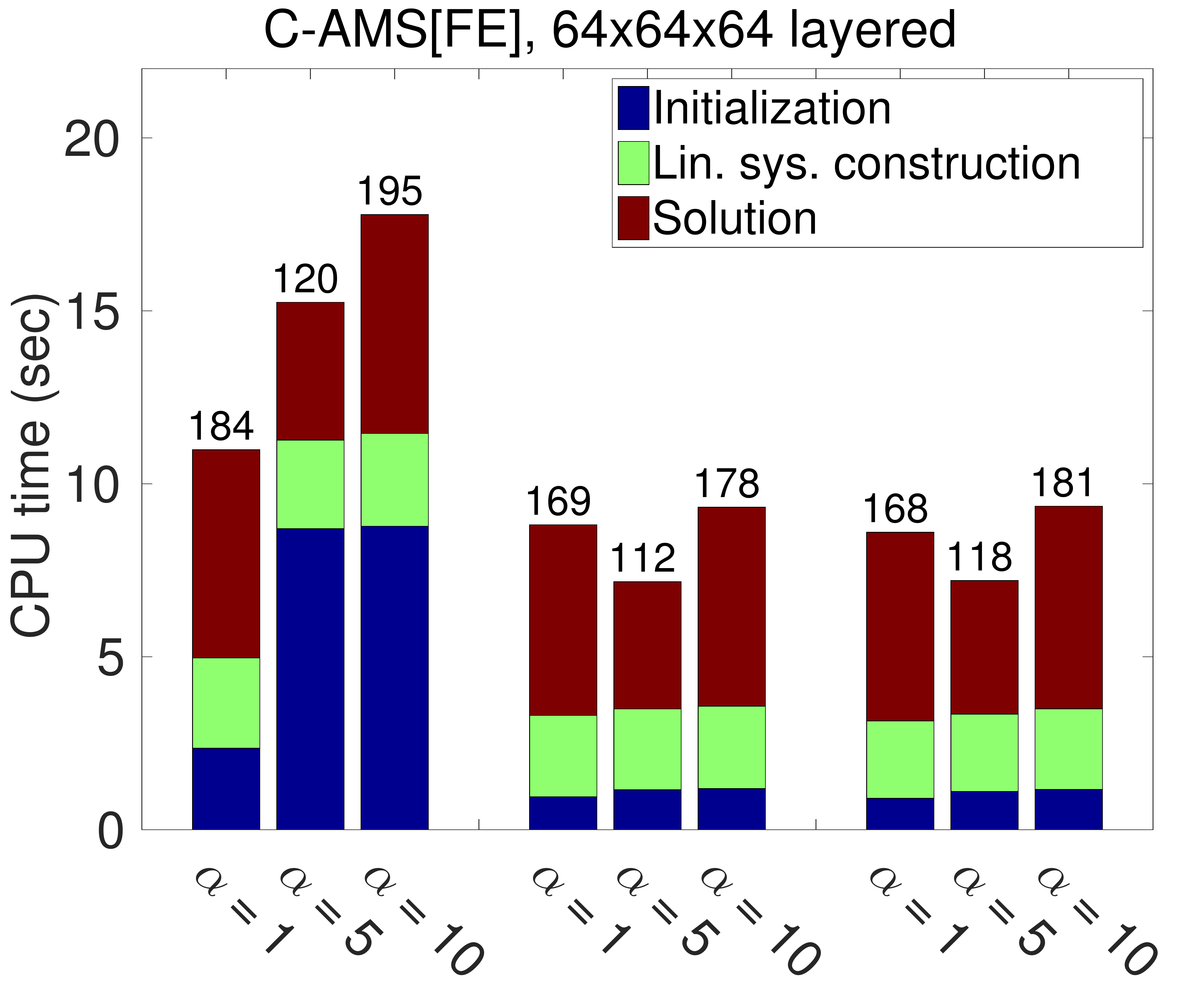}\\
\includegraphics[width=0.48\textwidth]{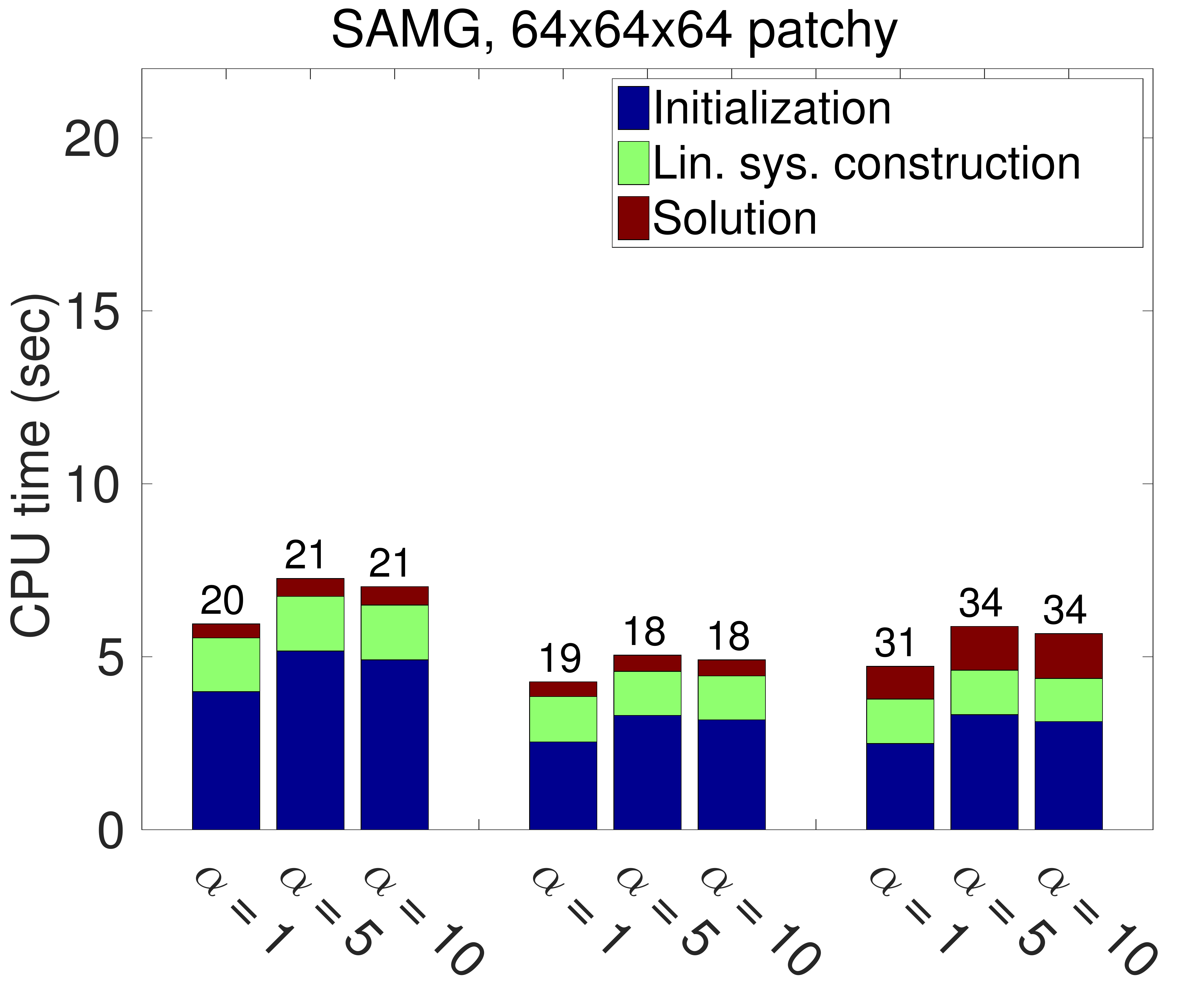}
\includegraphics[width=0.48\textwidth]{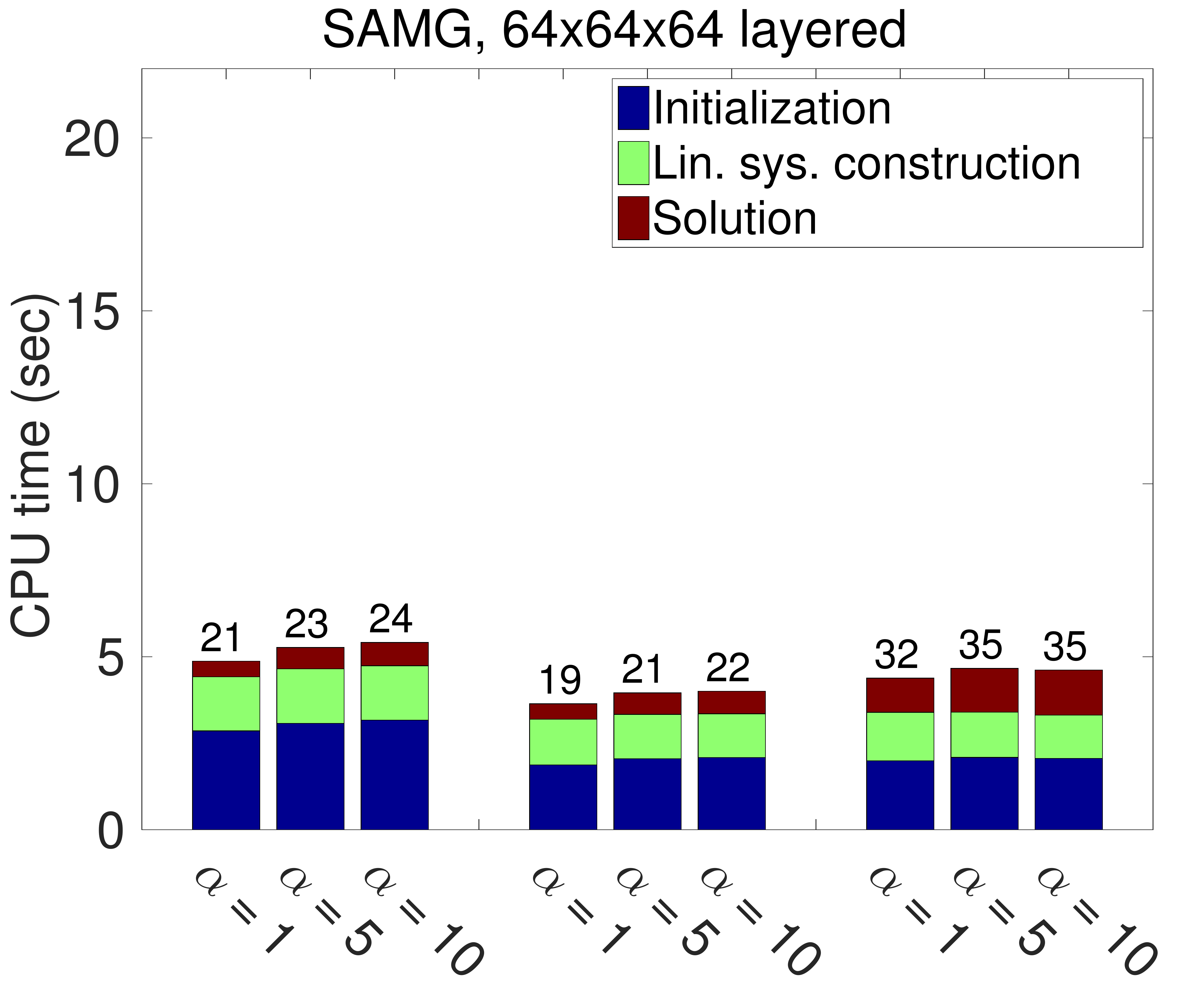}
\begin{picture}(360,0)
\put(-5,0){$t^*:$}
\put(20,0){0.0-0.4}
\put(76,0){0.4-1.0}
\put(132,0){1.0-2.0}
\put(214,0){0.0-0.4}
\put(270,0){0.4-1.0}
\put(320,0){1.0-2.0}
\end{picture}
\caption{\change{Performance of C-AMS (top) and SAMG (bottom) for permeability Set 1 from Table \ref{tab:permeability set} for different grid aspect rations $\alpha$ in $\Delta x/\alpha = \Delta y = \Delta z$ for three successive time steps. Pressure solutions for the first time step is shown for one of the realizations in Fig. \ref{stretched_p}.}}\label{stretched_cpu}
\end{figure}
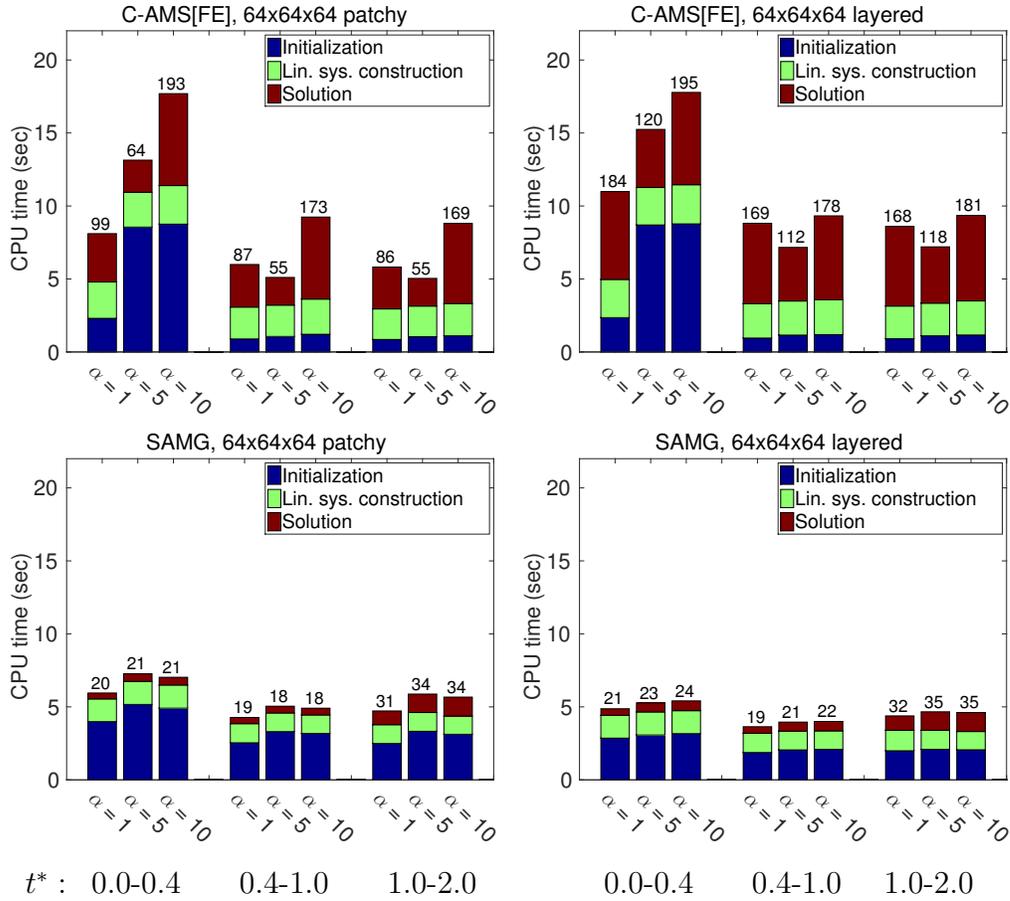

\change{Results shown in Fig. \ref{stretched_cpu} are obtained with the C-AMS coarsening ratios of $8\times8\times8$, $2\times8\times8$, and $2\times8\times8$ for the cases of $\alpha = 1, 5$, and $10$, respectively. Note that as shown in Fig. \ref{stretched_p}, the anisotropic transmissibility (caused by stretched grid effect) would further motivate the use of enhanced coarse-grid geometries for C-AMS. Such a strategy is well developed in algebraic multigrid community, and is the subject of our future studies.}

\clearpage
\subsubsection{\change{Test Case 3: effect of permeability contrast}}
\change{To study the effect of permeability contrast, permeability Set 1 from Table \ref{tab:permeability set} is considered with different $ln(k)$ variances of $\sigma = 2, 4,$ and $8$. Note that the so-far studied cases were for variance $4$, as described in Table \ref{tab:permeability set}. The settings are all the same as the default test cases, i.e., Dirichlet conditions are set at the east and west faces with no-flow condition everywhere else. Figure \ref{cpu_variance} illustrates the performances of C-AMS[FE] and SAMG for this test case. Note that the C-AMS requires more iterations when the permeability contrast is increased. To improve its performance, one can consider enriched multiscale strategies which are based on local spectral analysis \cite{Yalchin-enriched1}, and modified permeability field (with less contrast) for calculation of basis functions \cite{giuseppe-ib08}. Note that the success rates of C-AMS[FV] (not shown) were  90\% (patchy, $\sigma = 2$), 95\% (patchy, $\sigma = 8$) and 40\% (layered, $\sigma = 8$). For the successful runs, the CPU times of C-AMS[FV] were comparable with C-AMS[FE].}

\begin{figure} [htb!]
\centering
\includegraphics[width=0.48\textwidth]{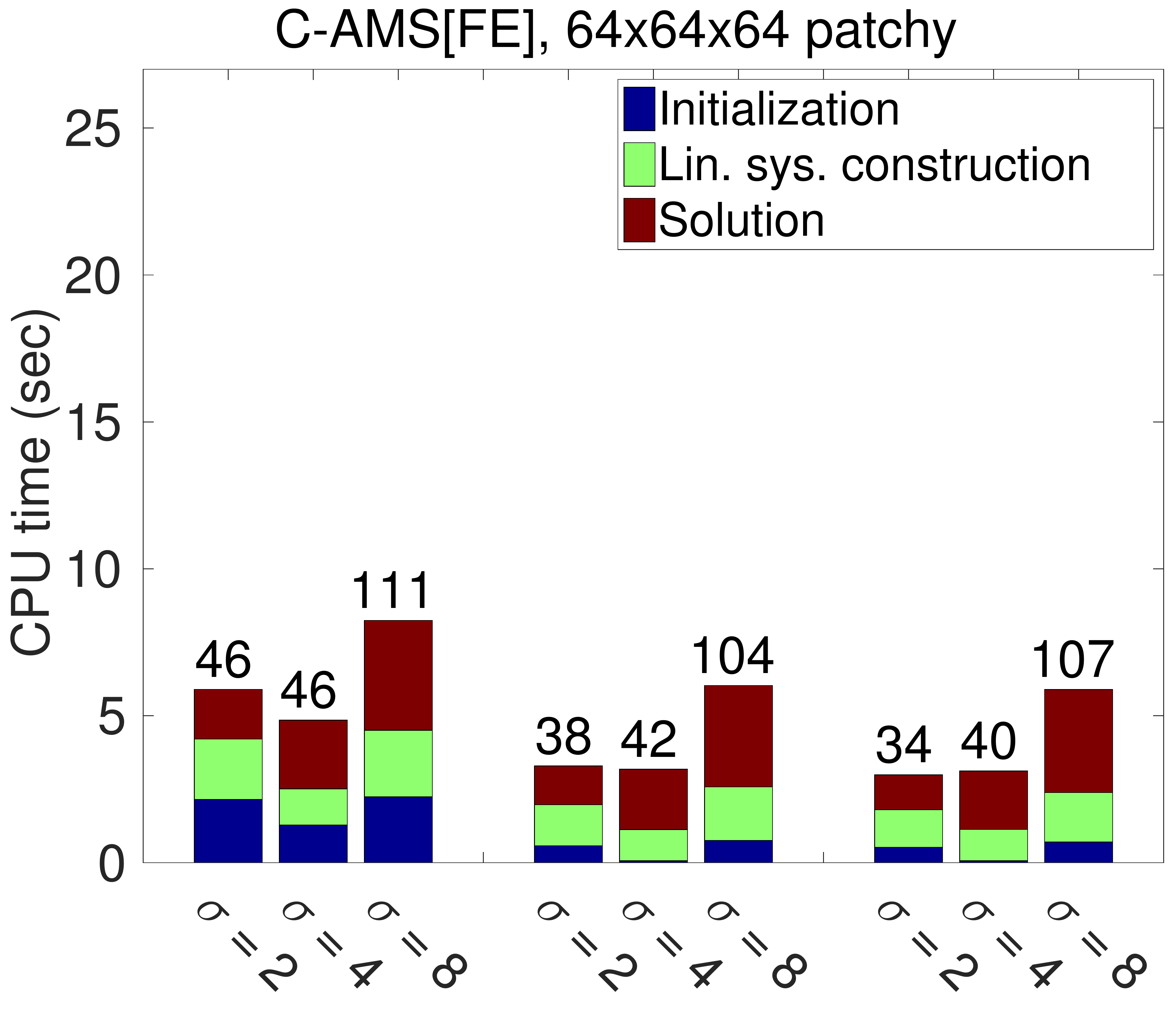}
\includegraphics[width=0.48\textwidth]{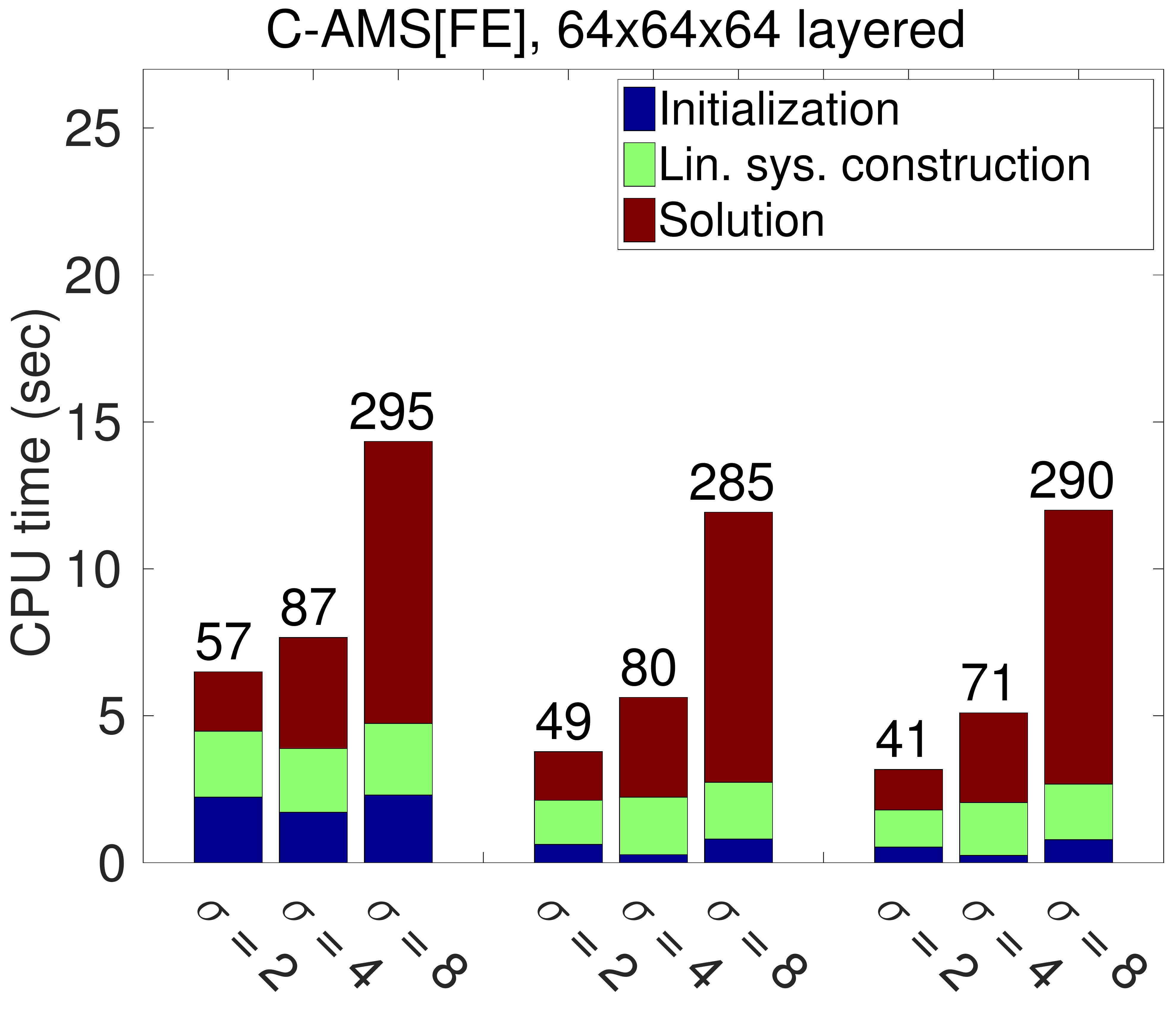}\vspace{2mm}
\includegraphics[width=0.48\textwidth]{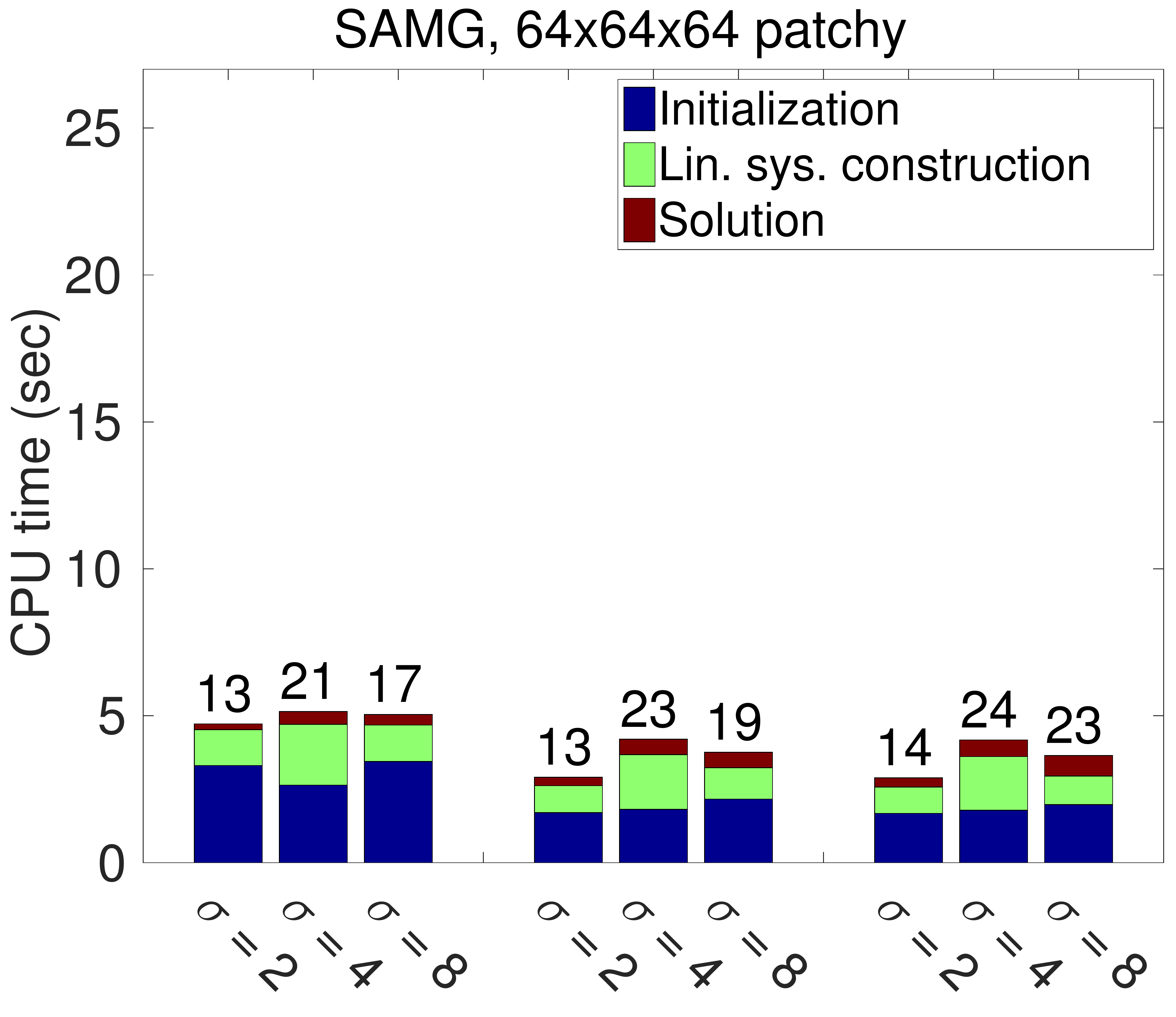}
\includegraphics[width=0.48\textwidth]{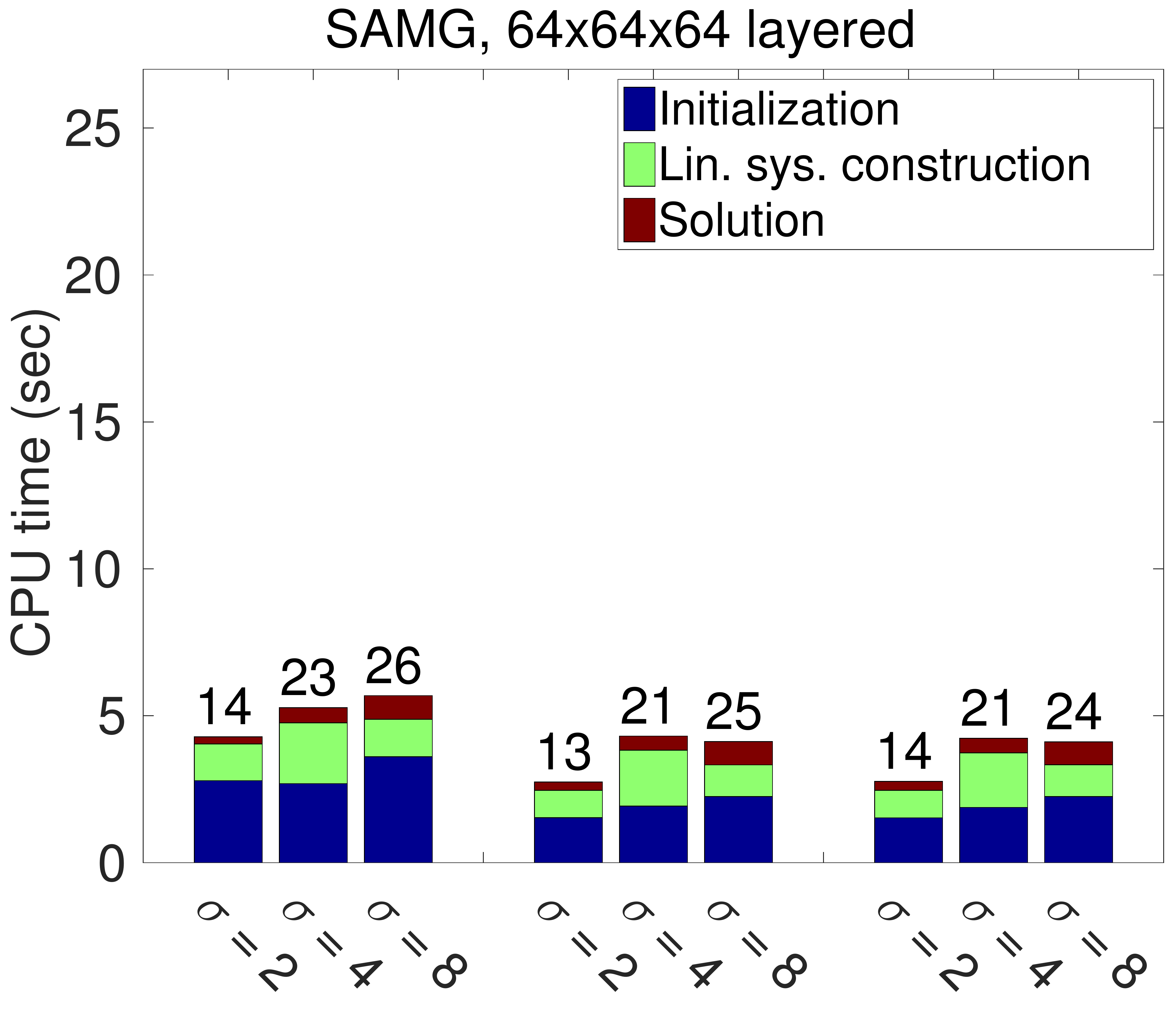}
\begin{picture}(360,0)
\put(-5,0){$t^*:$}
\put(23,0){0.0-0.4}
\put(78,0){0.4-1.0}
\put(132,0){1.0-2.0}
\put(215,0){0.0-0.4}
\put(270,0){0.4-1.0}
\put(328,0){1.0-2.0}
\end{picture}
\caption{\change{Averaged CPU time comparison between C-AMS (top) and SAMG (bottom) for permeability Set 1 from Table \ref{tab:permeability set} for different $ln(k)$ variances of $\sigma = 2, 4$ and $5$.}}\label{cpu_variance}
\end{figure}
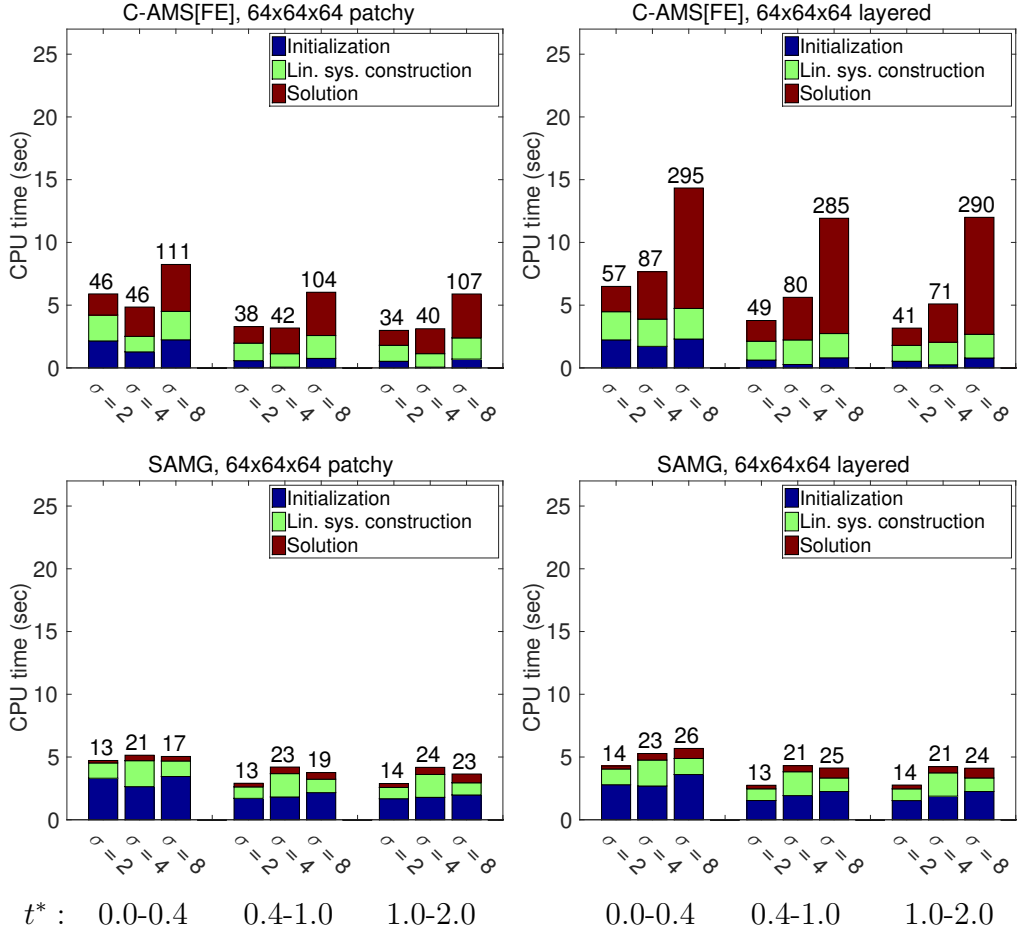


\section{Conclusions}

Algebraic Multiscale Solver for Compressible flows (C-AMS) in heterogeneous porous media \change{was introduced}. Its algebraic formulation benefits from adaptivity, both in terms of the infrequent updating of the linearized system and from the selective update of the basis functions used to construct the prolongation operator.

Extensive numerical experiments on heterogeneous patchy and layered reservoirs revealed that the most efficient strategy is to use basis functions with incompressible advection terms, paired with 5 iterations of ILU(0) for post-smoothing. 

Finally, \change{several} benchmark stud\change{ies were} presented, where the \change{developed} C-AMS research \change{similator} was compared with an industrial-grade multigrid solver, i.e., SAMG. The results show that C-AMS is a competitive solver, especially in experiments that involve the simulation of a large number of time steps. The only drawback is the relatively high initialization time, which can be reduced by choosing an appropriate coarsening strategy or by running the basis function updates in parallel \cite{Abdul-spej-parallel}. Moreover, due to its conservative property, C-AMS requires only a few iterations per time step to obtain a good quality approximation of the pressure solution for practical purposes. Systematic error estimate analyses for 3D multiphase simulations are a subject of ongoing research and, in addition, the C-AMS \change{performance} can be further extended by enrichment of the multiscale operators \cite{Yalchin-enriched1,Nataf-enriched-14,Davide14}, \change{and enriched coarse grid geometries on the basis of the underlying fine-scale transmissibility. Both are subjects of our future studies.}

\section{Acknowledgements}

We would like to acknowledge the financial support from the Chevron/Schlumberger Intersect Alliance Technology and Schlumberger Petroleum Services CV during Matei \c{T}ene's scientific visit at TU Delft between November - February 2014. Since March 2014, Matei \c{T}ene is a PhD Research Assistant at TU Delft, sponsored by PI/ADNOC. The authors also thank Prof. Hamdi Tchelepi of Stanford University for the many helpful discussions.

\bibliographystyle{unsrt}
\bibliography{Bib/Hadi}

\end{document}